\input epsf
\title{Normal zones on Zeeman manifolds with trace class heat and
Feynman kernels and well defined zonal Feynman integrals}

\newcommand{\cmt}[1]{\ifhmode\newline\fi{\sf *** \ \ #1 \\}}

\documentclass[11pt]{article}
\usepackage{latexsym}
\usepackage{epsfig}

\newtheorem{theorem}{Theorem}[section]

\newcommand{\R}{{\mathbf{R}}}
\newcommand{\Z}{{\mathbf{Z}}}

\newcommand{\heading}[1]{\vspace{1ex}\par\noindent{\bf #1}}

\def\:{\colon}

\long\def\onefigure#1#2{
\begin{figure*}[tbh]
\begin{center}
#1
\end{center}
\caption{#2}
\end{figure*}
} 

\newcommand{\iipefig}[1]  
{\smallskip\begin{center}\def\IPEfile{#1.ipe}\input{\IPEfile}\end{center}\smallskip}
\newcommand{\lipefig}[2]  
{\onefigure{\def\IPEfile{frh-#1.ipe}\input{\IPEfile}}{\label{f:#1} #2} }

\author{
{\sc Zolt\'an Imre Szab\'o}\thanks{Lehman College of CUNY, Dept. Math.,
Bronx NY 10468 USA and R\'enyi Institute of Mathematics, Budapest, Hungary.
email: zoltan.szabo@lehman.cuny.edu}
}  
\date{}

\begin{document}

\maketitle
\begin{abstract}
The problem of infinities (divergent integrals) 
appearing in quantum field theory
are treated by {\it renormalization} in the current theories. 
By this perturbative tool
the desired finite quantities are produced by differences of
infinities. The most common reason 
for these infinities appear is the infinite
trace of kernels such as the Wiener-Kac kernel $e^{-tH}$,
or, the Dirac-Feynman kernel $e^{-tH\mathbf i}$,
thus, they assign infinite measures to 
physical objects such as
self-mass, self-charge, e.t.c.. 

This paper gives a new non-perturbative approach to this problem. 
Namely, the Hilbert space $\mathcal H$, on which the quantum
Hamilton operator $H$ is acting, is decomposed into $H$-invariant
Zeeman zones on which
both the Wiener-Kac and the Dirac-Feynman kernels become
of the trace class, both defining the 
corresponding zonal measures on the path-spaces rigorously.

The Hamilton operators, $H$, considered in this paper are those
corresponding to electrons orbiting in a
constant magnetic field. This is one of the most important Hamiltonians,
introduced for explaining the 
Zeeman effect. This Hamiltonian 
is identified with the Laplacian of
of certain Riemannian, called Zeeman manifolds. 
The ``particles" modeled by them are 
called {\it Zeeman particles}. 
The spectral Zeeman zone decomposition is 
introduced, by one of the definitions, 
by the spectrum of the magnetic dipole moment operator. 
A zone corresponds to a magnetic state of these particles.

All the important spectral theoretical objects are
explicitly established. They include the spectrum,
the zonal projection operators, the zonal Wiener-Kac and
Dirac-Feynman kernels, and the zonal partition functions.  
\end{abstract}

\section{Introduction}

{\bf Zones on Zeeman manifolds}
\medskip\medskip

\noindent{\bf Infinities in Quantum Electrodynamics (QED).} 
The problem of infinities (divergent integrals), which is 
present in calculations since the early days of quantum field theory 
(Heisenberg-Pauli (1929-30)) or elementary particle physics 
(Oppenheimer (1930) and Waller (1930) in electron theory)
is treated by {\it renormalization} in the current theories. 
This perturbative tool provides
the desired finite quantities by differences of infinities. 
This problem originated from concepts such as {\it point mass}
and {\it point charge} of classical electron
theory, which provided the first warning that a point 
electron will have infinite electromagnetic
self-mass: the mass $e^2/6\pi ac^2$ for a surface distribution
of charge with radius $a$ blows up for $a\to 0$. In quantum
field theory the Hamiltonian of the field is proportional to this 
electromagnetic self-mass. This is why this infinity launched one 
of the deepest crisis-es in the history of physics.

The infinities mostly appear in the form of infinite
trace of kernels such as the Wiener-Kac kernel $e^{-tH}$,
or, the Dirac-Feyn\-man kernel $e^{-tH\mathbf i}$. 
The WK-kernel provides the fundamental solution of the heat equation
while the DF-kernel provides the fundamental
solution of the Schr\"odinger equation.  
The infinite trace assigns infinite measures to 
physical objects such as
self-mass, self-charge, e.t.c.. 
Because of the divergent integrals emerging
in its construction,
also the Feynman measure, which is analogous 
to the well defined Wiener-Kac measure on the path-spaces, 
requires renormalization. 

This paper offers a new non-perturbative approach to this problem.
The main idea in this approach can be briefly described as follows.
In the first step the quantum Hilbert space 
$\mathcal H$ 
(on which the quantum
Hamilton operator $H$ is acting) is decomposed into the direct
sum of $H$-invariant subspaces, called
Zeeman zones. Then all the operator-actions, 
such as the the heat- or Feynman-flows, 
are considered on these invariant subspaces separately. 
It turns out that both the Wiener-Kac and Dirac-Feynman kernels are 
of the trace class on each zone, furthermore, both define the 
corresponding zonal measures on the path-spaces rigorously.
\medskip

\noindent{\bf The Hamilton operator considered in this paper} is
the classical Zeeman operator 
\begin{equation}
\label{Zee_int}
H_Z=-{\hbar^2\over 2\mu}\Delta_{(x,y)} -
{\hbar eB\over 2\mu c\mathbf i}
D_z\bullet
+{e^2B^2\over 8\mu c^2}(x^2+y^2) ,
\end{equation}
of a free charged particle.
At the time it was introduced, the new feature of this operator was 
the orbital angular momentum operator 
$D_z\bullet =x\partial_y-y\partial_x$, which was the forerunner for
an adequate spin-concept. 
This operator is the result of a long agonizing creative 
effort \cite{to}, searching for a Hamilton operator explaining 
the Zeeman effect. (Note that the $D_z\bullet$
commutes with the rest part, $\mathbf O$, of the complete operator.
Thus the spectrum appears on common eigenfunctions, resulting that the
$D_z\bullet$ splits the spectral lines of $\mathbf O$ (Zeeman effect).) 
Pauli, who added a spin angular momentum operator to the orbital one,
developed the non-relativistic spin-concept. The relativistic concept
due to Dirac. 
Actually, the $H_Z$ is the Hamilton operator of an 
electron orbiting about the origin of the 
$(x,y)$-plane in a constant magnetic field $\mathbf K=B\partial_z$.
It had been established by means of the Maxwell equations.

\medskip

\noindent{\bf Mathematical modeling; Zeeman manifolds.}
An interesting feature of Zeeman operators, $H_Z$, is that they 
can be identified 
with the Laplace operators on certain Riemannian manifolds, namely,
with the Laplacians on two step nilpotent Lie groups endowed 
with the natural left invariant metrics. The details are as follows.

A 2-step nilpotent metric Lie group is defined on the
product $\mathbf v\oplus\mathbf z$ of Euclidean spaces,
$\mathbf v=\mathbf R^k$ and
$\mathbf z=\mathbf R^l$, where the
components are called X- and Z-space respectively. The
main object defining the Lie algebra is the linear space $J_{\mathbf z}$
of skew endomorphisms $J_Z,\, Z\in\mathbf z$, acting on the
X-space. The metric, $g$, is the left invariant extension of the
natural Euclidean metric on the Lie algebra. Particular
2-step nilpotent Lie groups are the
Heisenberg-type Lie groups which are defined by endomorphism spaces
satisfying the Clifford condition $J^2_Z=-|Z|^2id$. 
They are attached to Clifford modules.
Each group, $(N,g)$, extends into a solvable group $(SN,g_s)$. 

In this introduction just the relatively simple H-type
groups will be considered. On them the Laplacian appears in the form
\begin{equation}
\label{Delta}
\Delta=\Delta_X+(1+\frac 1{4}|X|^2)\Delta_Z
+\sum_{\alpha =1}^l\partial_\alpha D_\alpha \bullet,
\end{equation}
where $\partial_\alpha$ is partial derivative on the Z-space
$\mathbf z=\mathbf R^l$ and $D_\alpha\bullet$ denotes 
the directional derivative along
the X-field $J_\alpha (X)$ defined for the natural basis element 
$e_\alpha$ in the Z-space. Center periodic H-type groups are 
introduced by
factorizations, $\Gamma_\gamma\backslash H$, with Z-lattices 
$\Gamma_\gamma =\{Z_\gamma\}$ defined on the Z-space. 
On these manifolds the Laplacian appears
in a much more handy form. In fact, in this case the $L^2$ function
space is the direct sum of function
spaces $W_\gamma$ spanned by functions of the form
$\Psi_\gamma (X,Z)
=\psi (X)e^{2\pi\mathbf i\langle Z_\gamma ,Z\rangle}.
$
Each $W_\gamma$ is invariant under the action of the Laplacian, i. e., 
$\Delta \Psi_{\gamma }(X,Z)=\Box_{\gamma}\psi (X)
e^{2\pi\mathbf i\langle Z_\gamma ,Z\rangle}$,  
where operator $\Box_{\gamma }$, acting on 
$L^2(\mathbf v)$, is of the form
\begin{equation}
\label{Box}
\Box_{\gamma }
=\Delta_X + 2\pi\mathbf i D_{\gamma }\bullet -4\pi^2
|Z_\gamma |^2(1 + 
\frac 1 {4} |X|^2).
\end{equation}
  
Notice that (\ref{Zee_int}) is nothing but (\ref{Box}) 
on the 3D-Heisenberg
group. On a $(k+1)$-dimensional Heisenberg group, defined by a complex
structure $J$ acting on the even dimensional Euclidean space  
$\mathbf v=\mathbf R^k$, the Laplacian (\ref{Box}) appears in the form 
\begin{equation}
\label{Box2}
\Box_{\lambda}=\Delta_X +2 \mathbf i 
D_{\lambda }\bullet -
\lambda^2|X|^2 -4\lambda^2.
\end{equation}
Number $k/2$ is interpreted as the number of particles. The single
complex structure $J$ is interpreted such that these particles
are identical, rotating in the same plane defined by the same 
constant magnetic field $B$.

Operators (\ref{Zee_int}) and (\ref{Box2}) are identified by 
$H=-(1/2)\Box_\lambda$ and by the particular choice 
$\mu =\hbar =1, \lambda=eB/2c$ of the constants.
Operator (\ref{Box}) contains also the
constant $-4\pi^2|Z_\alpha |^2=4\lambda^2$, which is proportional
to $\hbar^2$ on the microscopic level, thus, it is usually neglected
in quantum physics. Also note that
the particles described by these Hamiltonians are free ($V=0$).

The Zeeman operator appears as Laplacian on center periodic
2-step nilpotent Lie groups in a more complex form.
These models represent $k/2$ number of charged particles, each of them
is orbiting in its own constant magnetic field. The system
can be in crystal states represented by the endomorphisms
$J_{\gamma}$. The Hamilton operators belonging to these crystal states
are $-{1\over 2}\Box_\gamma$.  This model matches
Dirac's famous multi-time theory. 

The Riemannian manifolds introduced so far are prototypes of a general
Zeeman manifold concept. This general concept
is beyond the scope of this article and will be developed in a 
subsequent paper. This exposition proceeds with considering
2-step nilpotent Lie groups.

\medskip

\noindent{\bf Introducing the zones.}
In \cite{sz5,sz6} the problem of infinities is approached 
by the above mentioned {\it Zeeman zone
decomposition} of the Hilbert space of complex valued 
$L^2$-functions defined on the X-space. 
This {\it Hilbert space}, $\mathcal H$,
is isomorphic to the weighted space defined by the Gauss density 
$d\eta_\lambda (X)=e^{-\lambda |X|^2}dX$. 
This Hilbert space is spanned by the complex valued polynomials. 

Next the Hilbert space is interpreted in this way.
The natural {\it complex Heisenberg group 
representation} on $\mathcal H$ is defined by  
\begin{equation}
\label{rho}
\rho_{\mathbf c} (z_i)(\psi )=(-\partial_{\overline{z}_i}+
\lambda_i z_i\cdot )\psi\quad ,\quad
\rho_{\mathbf c} (\overline{z}_i)(\psi )=\partial_{z_i}\psi , 
\end{equation}
where $\{z_i\}$ is a complex coordinate system on the X-space.
This representation is reducible. In fact, it is irreducible on the
space of holomorphic functions, where it is called Fock representation.
Besides the holomorphic subspace there are infinitely many
other irreducible invariant subspaces. In the literature
only the Fock representation, defined on the space of 
holomorphic functions, is well known. The above representation is
called {\it extended Fock representation}. In the function operator
correspondence, this representation associates operator (\ref{Zee_int})
to the Hamilton function of an electron orbiting in a constant 
magnetic field.

The zones are defined in two different ways. First, they can be defined
by the invariant subspaces of representation (\ref{rho}). The actual 
construction uses Gram-Schmidt orthogonalization. On the complex plane
$\mathbf v=\mathbf C$, which corresponds to the 2D-Zeeman operator 
(\ref{Zee_int}),
the first zone, $\mathcal H^{(0)}$, is the holomorphic zone spanned
by the holomorphic polynomials. To construct the second zone
one considers the function space $G^{(1)}$ consisting of functions 
of the form $\overline zh$, where the $h$ is a holomorphic function.
Then $\mathcal H^{(1)}$ is the orthogonal component of $G^{(1)}$
to the previous zone. E. t. c., one can construct
all the zones, $\mathcal H_\lambda^{(a)}$, by continuing the
Gram-Schmidt orthogonalization applied to function spaces $G^{(a)}$
spanned by functions of the form $\overline z^ah$. 
The zone index $a$ indicates the maximal number
of antiholomorphic coordinates $\overline z_i$ in the polynomials
spanning the zone. 

In the 2D-case the zones are irreducible under 
the action of the extended
Fock representation. In the higher dimensional cases the Gram-Schmidt
process results reducible zones, 
called {\it gross zones.} More precisely,
the holomorphic zone is always irreducible and the gross zones
of higher indexes decompose
into {\it irreducible zones}, which can also be explicitly described.

The second technique defines the very same zones by  
computing the spectrum and the corresponding 
eigenfunctions explicitly.  According to these computations, the
eigenfunctions appear in the form
\begin{equation}
h^{(p,\upsilon)}(X)=
H^{(p,\upsilon)}(X) 
e^{-\lambda |X|^2/2}
\end{equation}
with the corresponding eigenvalues
\begin{equation}
\label{eigv_int}
-((4p+k)\lambda +4k\lambda^2),  
\end{equation}
where $p$ resp. $\upsilon$ are the holomorphic resp. antiholomorphic
degrees of the polynomial $H^{(p,\upsilon )}$.
Numbers $l=p+\upsilon$ and $m=2p-l$ are
called azimuthal and magnetic quantum numbers respectively. The
above function is an eigenfunction also of the magnetic dipole
moment operator with eigenvalue $m$.

Now the zones are created such that the above 
eigenfunction falls into the zone with the zone index $\upsilon$. 
According to the formula $\upsilon ={1\over 2}(l-m)$, 
the zones are determined by the magnetic quantum number $m$.
Thus a zone exhibits the magnetic state of a zonal particle.

Note that the eigenvalues (\ref{eigv_int}) are independent of the 
antiholomorphic index and they depend just on the holomorphic index. 
As a result,
each eigenvalue has infinite multiplicities. 
On the irreducible zones, however,
each multiplicity is $k/2$. (Here we suppose that there is only one 
parameter $\lambda$ involved, meaning that the particles are identical.
If the particles (i. e., the $\lambda_i$'s) 
are properly distinct, then the 
multiplicity is $1$ on each zone.) Moreover, 
two irreducible zones are isospectral.
\medskip\medskip

\noindent{\bf Normal de Broglie Geometry}
\medskip\medskip

\noindent{\bf Introducing the point-spreads by projection 
kernels.} It is remarkable that all the 
important objects such as the projections onto the zones,
the zonal heat- and Feynman-kernels and their well defined trace,
and the several Feynman-Kac type formulas 
can explicitly be computed \cite{sz5,sz6}. 

First the {\it projection operators}, $\delta^{(a)}$, 
onto the zones $\mathcal H^{(a)}$ are established. 
If $\{\varphi_i^{(a)}\}_{i=1}^\infty$ is an 
orthonormal basis in $\mathcal H^{(a)}$, 
then the corresponding projection can be
formally defined as convolution with the kernel 
\begin{equation}
\delta^{(a)}(z,w)=\sum_i \varphi^{(a)}_i(z)
\overline{\varphi^{(a)}_i(w)},
\end{equation}
where $z$ and $w$ represent complex vectors on 
$\mathbf C^{k\over 2}=\mathbf R^k$.
Interestingly enough, these self-adjoint operators are integral
operators having smooth Hermitian integral kernels.
These kernels can be interpreted as restrictions of the 
global Dirac delta distribution, 
$\delta_z(w)=\sum \varphi_i(z)\overline{\varphi}_i(w)$, onto the
zones. They have the following
explicit form
\begin{equation}
\delta_{\lambda z}^{(a)}(w)
= {\lambda^{k/2}\over \pi^{k/2}}L^{((k/2)-1)}_a(\lambda |z-w|^2)
e^{\lambda (z\cdot\overline{w}-{1\over 2}(|z|^2+ |w|^2)},
\end{equation}
where $L_a^{((k/2)-1)}(t)$ is the corresponding Laguerre polynomial.
Among these kernels only the first one, the projection kernel
onto the holomorphic zone, is well known. It is nothing but the
Bergman kernel. The new mathematical feature of these formulas is
that they are explicitly determined regarding each zone and not just
for the holomorphic zone.

These kernels represent one of the most important concepts
in this theory. They can be interpreted such that, 
on a zone, a point particle appears as a spread
described by the above wave-kernel. 
Note that how these kernels, called zonal point-spread, are derived
from the one defined for the holomorphic zone. This holomorphic 
spread is just multiplied by the radial 
Laguerre polynomial corresponding to the zone. These point-spreads
show the most definite similarity to the de Broglie waves
packets (cf. \cite{bo}, pages 61).
These zonal kernels can be interpreted such that 
a point particle concentrated
at a point $Z$ appears on the zone as an object which spreads around 
$Z$ as a wave-package with wave-function 
described by the above kernel explicitly. 

The wave-package interpretation of physical objects started out with the
de Broglie theory. This concept was finalized in the Schr\"odinger
equation. The mathematical formalism did not follow this development,
however, and the Schr\"odinger theory is built up on such mathematical
background which does not exclude the existence of the controversial
point objects. On the contrary, 
an electron must be considered as a point-object in
the Schr\"odinger theory as well 
(cf. Weisskopf's argument on this problem in \cite{schw, sz5}). 
An other demonstration for 
the presence of
point particles in the classical theory is the duality principle, 
stating that objects manifest themselves sometime as waves and sometime
as point particles. The bridge between the two visualizations is
built up in Born's probabilistic theory,
where the probability for that a particle, 
attached to a wave $\xi$, can be found on a domain $D$ is measured by
$\int_D\xi\overline\xi$.

These controversial point-objects, by having infinite self-mass or
self-charge attributed to them 
by the Schr\"odinger equation, launched one
of the deepest crisis's in the history of physics.
In the zonal theory de Broglie's idea is established on a mathematical
level. Although the points are ostracized from this theory, 
the point-spreads 
still bear some reminiscence of the point-particles. For instance, they
are the most compressed wave-packages and all the 
other wave-functions in the 
zone can be expressed as a unique superposition of the point-spreads.
If $\xi$ is a zone-function, the above integral 
measures the probability of that that
the center of a point-spread is on the domain $D$. This 
interpretation restores, in some extend, the duality principle in the
zonal theory.

Function 
$
\delta^{(a)}_{\lambda Z}
\overline{\delta}^{(a)}_{\lambda Z}
$
is called the density of the spread around $Z$. By this reason, 
function 
$
\delta^{(a)}_{\lambda Z}
$
is called spread-amplitude. Both the spread-amplitude and 
spread-density generate well defined measures on the path-space
consisting of continuous curves connecting two arbitrary points.
Both measures can be constructed by the method applied in constructing
the Wiener measure.
 
The point-spread 
concept bears some remote reminiscence of Heisenberg's 
suggestion (1938) for the existence of 
a fundamental length $L$, analogously
to $h$, such that field theory was valid only for distances larger than
$L$ and so divergent integrals would be cut off at that distance.
This idea has never became an effective theory, however. 
Other distant relatives of the point-spread concept are the 
smeared operators, i. e. those suitably averaged over small 
regions of space-time, considered by Bohr and Rosenfeld in quantum 
field theory. There are also other theories where an electron is 
considered to be extended. Most of them fail on lacking the
explanation for the question: 
Why does an extended electron not blow up? The zonal theory is checked
against this problem in \cite{sz5}, section 
(F) ``Linking to the blackbody radiation; Solid zonal particles''.
\medskip

\noindent{\bf Global Wiener-Kac and Feynman flows.}
Both definitions imply that the zones are invariant 
under the action of the Hamilton operator,
therefore, kernels such as the {\it heat (Wiener-Kac) and Feynman 
kernels} can be restricted onto them. The 
zonal kernels are defined by these restrictions. 

Since the spectrum is discrete, also the global kernels, 
defined for the total space $\mathcal H$, can be introduced 
by the trace formula using an orthonormal basis
consisting of eigenfunctions on the whole space $\mathcal H$.
Despite of the infinite multiplicities on the global setting,
both global kernels are well defined smooth functions. 
If the Zeeman operator $H_Z$ is non-degenerated and the
distinct non-zero parameters $\{\lambda_i\}$,  $i=1,\dots ,r$,
are defined on $k_i$-dimensional subspaces, then for the Wiener-Kac 
kernel we have
\begin{eqnarray}
\label{d_1gamm}
d_{1\gamma} (t,X,Y)=e^{-tH_Z}(t,X,Y)=\\
=\prod \big({\lambda_i\over 2\pi sinh(\lambda_i t)}\big)^{k_i/2}e^{
-\sum{\lambda_i}({1\over 2}coth(\lambda_i t)
|X_i-Y_i|^2+\mathbf i
\langle X_i,J (Y_i)\rangle}.\nonumber
\end{eqnarray}
This kernel satisfies the Chapman-Kolmogorov identity as well as the
limit property $\lim_{t\to 0_+}d_1(t,X,Y)=\delta(X,Y)$, however, 
it is not of the trace class. Thus functions such as 
the partition function or the zeta function can not be defined in the
standard way. Note that by regularization (renormalization) only
well defined relative(!) partition and zeta functions are introduced..

The explicit form of the global Feynman-Dirac kernel is
\begin{eqnarray}
d_{\mathbf i}(t,X,Y)=e^{-t\mathbf iH_Z}(t,X,Y)=\\
=\prod\big({\lambda_i\over 2\pi
\mathbf i sin(\lambda_i t)}\big)^{k_i/2}e^
{\mathbf i\sum{\lambda_i}\{
{1\over 2} cot(\lambda_i t)|X_i-Y_i|^2-\langle X_i,J(Y_i)\rangle \}}.
\nonumber
\end{eqnarray}
Since for fixed $t$ and $X$ function $d_{\mathbf i}(t,X,Y)$ is neither
$L^1$- nor $L^2$-function of variable $Y$, 
the integral required
for the Chapman-Kolmogo\-rov identity is not defined for this kernel. 
It is not of the trace class either. Nevertheless, it satisfies the 
above limit property. Thus the
constructions with the global Feynman-Dirac kernel lead to divergent 
integrals in the very first step.

It is well known in the history that 
Kac, who tried to understand Feynman, was able to introduce a 
rigorously defined
measure on the path-spaces only by the kernel $e^{-tH}$. 
This measure was, actually, established earlier by Wiener for the 
Euclidean Laplacian $\Delta_X$. 
Note that the heat kernel involves a Gauss density which makes this
constructions possible. Whereas, the Feynman kernel does not involve
such term. This is why no well defined constructions can be
carried out with this kernel. One can strait out all this
difficulties, however, by considering these constructions 
on the zones separately.
\medskip

\noindent{\bf Zonal Wiener-Kac and Feynman flows.}
Also the zonal WK- resp. FD-kernels are well defined smooth functions.
The gross zonal
Wiener-Kac kernels are of the trace class, which can be
described, along with their partition functions, by the following
explicit formulas. 
\begin{eqnarray}
\label{d_1^a}
d_{1}^{(0)}(t,X,Y)=
\prod\big({\lambda_i e^{-\lambda_i t}\over \pi}\big)^{k_i\over 2}
e^{\sum \lambda_i(-{1\over 2} (|X_i|^2+|Y_i|^2)+ e^{-2\lambda_i t}
\langle X_i,Y_i+\mathbf iJ(Y_i)\rangle )},
\\
d_{1}^{(a)}(t,X,Y)=
(L^{({k\over 2}-1)}_a(\sum\lambda_i |X_i-Y_i|^2)+LT_1^{(a)}
(t,X,Y))d_1^{(0)}(t,X,Y),
\end{eqnarray}
where 
$LT_1^{(1)}$ 
is of the form 
\begin{eqnarray}
LT_1^{(1)}(t,X,Y))=
(1-e^{-2t})
lt_1^{(1)}(t,X,Y))=
\\
(1-e^{-2t})
(|X|^2+|Y|^2-1-
(1+e^{-2t})\langle X,Y+\mathbf iJ(Y)\rangle ) \nonumber
\end{eqnarray}
and for the general terms,
$LT_1^{(a)}$,
recursion formula can be established. Furthermore,
\begin{eqnarray}
\mathcal Z_1^{(a)}(t)=
Trd_{1}^{(a)}(t)
={a+(k/2)-1\choose a}\prod {e^{-{k_i\lambda_i t\over 2}}\over
(1-e^{-2\lambda_i t})^{k_i\over 2}}=
TrD_{1}^{(a)}(t),
\end{eqnarray}
where 
$
D_{1}^{(a)}(t,X,Y)=
L^{({k\over 2}-1)}_a(\sum\lambda_i|X_i-Y_i|^2)
d_{1}^{(0)}(t,X,Y)
$
is the dominant zonal kernel. The remaining long term kernel
in the WK-kernel vanishes for 
$\lim_{t\to 0_+}$ and is of the 0 trace
class.
The zonal WK-kernels satisfy the Chapman-Kolmogorov identity 
along with the limit property
$\lim_{t\to 0_+}d_1^{(a)}=\delta^{(a)}$. 

Similar statements can be established regarding the zonal
DF-flow. The gross zonal 
Dirac-Feynman kernels are of the trace class which, together with their
partition functions, can be described 
by the following explicit formulas.

\begin{eqnarray}
\label{d_i^a}
d_{\mathbf i}^{(0)}(t,X,Y)=
\prod\big({\lambda_i e^{-\lambda_i t\mathbf i}\over \pi}
\big)^{k_i\over 2}
e^{\sum \lambda_i(-{1\over 2} (|X_i|^2+|Y_i|^2)+ 
e^{-2\lambda_i t\mathbf i}
\langle X_i,Y_i+\mathbf iJ(Y_i)\rangle )},
\\
d_{1}^{(a)}(t,X,Y)=
(L^{({k\over 2}-1)}_a(\sum\lambda_i |X_i-Y_i|^2)+LT_{\mathbf i}^{(1)}
(t,X,Y))d_1^{(0)}(t,X,Y),
\end{eqnarray}
where 
$LT_{\mathbf i}^{(1)}$ 
is described by 
\begin{eqnarray}
LT_{\mathbf i}^{(1)}(t,X,Y))=
(1-e^{-2t\mathbf i})
lt_{\mathbf i}^{(1)}(t,X,Y))=
\\
(1-e^{-2t\mathbf i})
(|X|^2+|Y|^2-1-
(1+e^{-2t\mathbf i})\langle X,Y+\mathbf iJ(Y)\rangle )\nonumber
\end{eqnarray}
and a
general long term, 
$LT_{\mathbf i}^{(a)}$,
can be defined recursively. Furthermore, 
\begin{eqnarray}
\mathcal Z_{\mathbf i}^{(a)}(t)=
Trd_{\mathbf i}^{(a)}(t)
={a+(k/2)-1\choose a}\prod {e^{-{k_i\lambda_i t\mathbf i\over 2}}\over
(1-e^{-2\lambda_i t\mathbf i})^{k_i\over 2}}=
TrD_{\mathbf i}^{(a)}(t),
\end{eqnarray}
where 
$
D_{\mathbf i}^{(a)}(t,X,Y)=
L^{({k\over 2}-1)}_a(\sum\lambda_i|X_i-Y_i|^2)
d_{\mathbf i}^{(0)}(t,X,Y)
$
is the dominant kernel. The remaining longterm term in the zonal
DF-kernel is of the 0 trace class.

The zonal DF-kernels are zonal fundamental solutions 
of the Schr\"o\-dinger equation
$(\partial_t+\mathbf i(H_{Z})_{X})d^{(a)}_{\mathbf i\gamma}(t,X,Y)=0$,
satisfying the Chapman-Kolmogorov identity as well as 
the limit property $\lim_{t\to 0_+}d_{\mathbf i}^{(a)}=\delta^{(a)}$.
 
On each zone,
both the WK- and the FD-kernels are of the 
trace class, moreover, both define rigorous
complex zonal measures, 
the zonal Wiener-Kac measure 
$dw_{1xy}^{T(a)}(\omega)$ 
and the zonal Feynman measure 
$dw_{\mathbf ixy}^{T(a)}(\omega)$, 
on the space of continuous curves $\omega :[0,T]\to \mathbf R^k$ 
connecting two points $x$ and $y$. The existence of zonal WK-measure
is not surprising, since this measure exists even for the global 
setting. However, the trace class property is a new feature indeed. In
case of the zonal Feynman measure both the trace class property 
and the existence are new features. Note that the zonal DF-kernels
involve a Gauss density which makes these constructions well defined.
\newpage

\begin{center} 
\large{\bf PART ONE}
\end{center} 
\begin{center} 
 \large{\bf ZONES ON ZEEMAN MANIFOLDS}
\end{center} 
\medskip

\section{Zeeman operators}

\noindent{\bf Zeeman operators defined physically.}
There is a formal correspondence between Quantum Theory and
Classical Theory, formulated by Bohr (1923) 
as Correspondence Principle. There are many aspects of this 
complex principle.
In short, Quantum Theory substitutes
the continuous functions describing the physical systems in Classical
Physics by operators
acting on a complex Hilbert space $\mathcal H$. The discontinuities
experienced on the quantum (microscopic) level appear in the discrete
spectrum of these operators.
This replacement, called {\it function-operator correspondence}, 
actually means derivation of the operators from the functions of
the classical theory.  
For instance, the Hamilton operators are derived from the Hamilton
functions describing the total energy of the physical systems. 
The elaboration of a correct correspondence is guided by the 
requirement (principle): Quantum Theory must approach Classical
Theory asymptotically in the limit of large quantum numbers. 
 
It is obvious that 
the function-operator correspondence 
should be defined, on the first place, for the canonical
coordinates $(q_i,p_i=\dot{q}_i)$. In a time independent
physical system the other functions depend only on these coordinates,
therefore, one should just extend this correspondence
to the other functions, such as the Hamilton functions, 
by methods known both in physics and geometric
quantization. (With some caution, the Taylor expansion is 
a natural tool for this extension.) 

The Heisenberg Lie algebras can be defined
by restricting the Poisson brackets
\begin{equation}
\{f,g\}=\sum
\partial_{q_i}(f)\partial_{p_i}(g) -
\partial_{p_i}(f)\partial_{q_i}(g)
\end{equation}
onto the linear space of functions spanned by the functions 
$\{q_i,p_i,1\}$. 
The exponential map maps
this algebra to the Heisenberg group. 

In the quantization process one considers a
unitary representation of the Heisenberg group on a 
complex Hilbert space. 
In case of the {\it real Heisenberg groups}, 
this Hilbert space, $\mathcal H$, is the 
$L^2$-space of complex valued functions depending on the position 
coordinates $q_i$ only.
I. e., $\mathcal H=L_{\mathbf C}^2(\mathbf v_q)$, 
where, on the linear space $\mathbf v_q$ of position vectors, 
the complex inner product is defined by 
$\langle f,g\rangle =\int f\overline gdq$.
{\it Representation $\rho :f\to f_\rho$ of the real Heisenberg 
group} is defined by
\begin{equation}
\label{hrep}
\rho (q_i)(\psi )=q_i\psi\quad ,\quad
\rho (p_i)(\psi )={\hbar\over \mathbf i}\partial_{q_i}
\psi\quad ,\quad \rho (1)=id.
\end{equation}

The associated operators,
$\tilde{f}=\rho (f)$, satisfy the Heisenberg commutative relations:
\begin{equation}
\label{comrel}
[\tilde{q}_i,\tilde{q_j}]=0\quad ,\quad
[\tilde{q}_i,\tilde{p_j}]=\hbar\mathbf i\delta_{ij}\quad ,\quad
[\tilde{p}_i,\tilde{p_j}]=0,
\end{equation}
therefore, the $\rho$ is a Lie algebra representation and 
the exponential map $e^{\mathbf i\rho}$ 
defines a unitary representation of
the Heisenberg group. Representation $\rho$ is called unitary also
on the Lie algebra level. 

It is well known that
this representation is irreducible. By the classical Neumann-Stone
theorem the irreducible unitary representations of Heisenberg's
groups on a complex infinite dimensional Hilbert 
space are unique up-to multiplications with
unit complex numbers.

The {\it Zeeman-Hamilton operator} was the result of a longstanding
creative effort, searching for a Hamilton operator explaining the 
Zeeman effect. The
exact form of the Hamiltonians was established by considering
an electron revolving in a magnetic field 
$\mathbf K$. The classical Hamilton function
of this system can be determined by the Maxwell equations 
\begin{equation}
\mathbf K=\nabla\times\mathbf a\quad ,\quad
\mathbf E=-{1\over c}{\partial\mathbf a\over \partial t}-\nabla\phi ,
\end{equation}
describing the electromagnetic field $(\mathbf E,\mathbf K)$ in
terms of the vector-, $\mathbf a$, and scalar-potential, $\phi$.

According to these equations, the complete classical Hamilton function  
of a charged particle
in an electromagnetic field that is specified externally
(produced by charges and currents other then the one considered)
\begin{equation}
H={1\over 2M}|\mathbf p-{Q\over c}\mathbf a|^2+Q\phi + V +
{Q^2\over 16Mc^2\pi}\int(\mathbf{E}^2+\mathbf{K}^2)dq ,
\end{equation}
where $M$ is the mass, $\mathbf p$ is the kinetic moment, $Q$ is the
charge (for an electron $Q=-e$), $V$ is that part of the potential 
energy which is of non-electromagnetic origin, and the last
term is the electromagnetic energy of the field in which the 
electron is staying.
By (\ref{hrep}), the Hamilton operator corresponded to this
Hamilton function is:
\begin{eqnarray}
\label{magnH}
\rho (H)=-{\hbar^2\over 2M}\Delta -{\hbar e\over 2Mc}\mathbf i
(\mathbf a\cdot\nabla +\nabla\cdot\mathbf a)
+{e^2\over 2Mc^2}|\mathbf a|^2+\\
+e\phi + V +
{e^2\over 16Mc^2\pi}\int(\mathbf{E}^2+\mathbf{K}^2)dq.\nonumber
\end{eqnarray}

In case of the Zeeman-Hamilton operator one assumes that the electron
is orbiting in the $(x,y)$-plane in a
constant magnetic field $\mathbf K=\partial_z$ that is directed in the 
$z$-direc\-tion and derived from the vector potential
\begin{equation}
\mathbf a=-{1\over 2}(x\partial_y-y\partial_x).
\end{equation}
Then $\mathbf K=\nabla\times\mathbf a=\partial_z$ holds. 
There is also supposed that $\mathbf E=0$, from which $\phi=0$ follows.
The energy of the constant magnetic field provided by the last
integral term of (\ref{magnH}) is obviously infinity. This integral
is infinity also for Coulomb electric fields $\mathbf E$ even if
$\mathbf K=0$. These were the first infinities appearing in the
early days of Quantum Field Theory which were followed by many
others in the history of Quantum Theory. Without this confusing term
the Zeeman Hamiltonians, function and operator, are of the form
(cf. \cite{bo}, formula (15.26.49)):
\begin{eqnarray}
\label{Zee}
H_z={1\over 2M}|\mathbf p+{eB\over c}J_z(x,y)|^2 + V ,\\
H_Z=\rho (H_z)=-{\hbar^2\over 2M}\Delta 
-{\hbar eB\over 2Mc}\mathbf i D_z\bullet
+{e^2B^2\over 8Mc^2}(x^2+y^2) +V .
\end{eqnarray}

So far the constant magnetic field was externally specified and not
attributed to the charged particle itself. If one supposes that
the particle has ``built in'' angular momentum then the influence of 
the magnetic dipole moment
\begin{equation}
\mathcal S=-{e\over 2Mc}\mathcal L
\end{equation}
associated with the angular momentum
\begin{equation}
\mathcal L=h(J_x=yp_z-zp_y,J_y=xp_z-zp_x,J_z=xp_y-yp_x)
\end{equation}
will be felt on subjecting the particle to a magnetic field.
The first experiment to observe the magnetic moment of the
electron directly is due to Stern and Gerlach (1922), justifying
the hypothesis of a ``built in'' orbital angular momenta in these
particles. 
Note that the angular momentum operator in the Zeeman operator 
commutes with the rest part, $\mathbf O$, of $H_Z$. 
Thus the spectrum appears on common eigenfunctions, meaning
that the angular momentum operator splits the spectral lines
of operator $\mathbf O$. This so called Zeeman effect was
observed by the Stern-Gerlach experiment, which proved that
the angular momentum is quantized.
 
Still in this
chapter the $H_Z$ for a free particle, $V=0$, is 
established as the Laplacian on a Riemannian
manifold. An interesting feature of this interpretation is
that there will be a well defined finite constant term corresponding
to the controversial electromagnetic self energy in the above formulas.
In the Point-Spread Theory developed later in this paper 
this term, $H_f$, corresponds to the field
energy of the ``little magnet inside of the point-spread''.
This combined particle-field Hamiltonian will be denoted by 
$H_{Zf}=H_Z+H_f$.

\medskip

\noindent{\bf 2D-Zeeman operators and extended Fock representations.} 
The above Hamiltonians are considered  
in the literature mostly on $\mathbf R^3$,
meaning that operator $\Delta$
is the Laplacian on the 3-space. However, also
the restriction of this 3-dimensional operator to the $(x,y)$-plane 
is well known in the literature. Note that (\ref{Zee}) 
becomes a 2D-operator by considering the
2D-Laplacian $\Delta_{(x,y)}$. Then the 3- and 2-dimensional versions
differ from each other just by the operator
$
-{(\hbar^2/2\mu)}\partial^2_z
$,
thus the spectral computations with respect to these two cases
can be easily compared. The 2D-operators have intensely been 
investigated, since its 
first appearance in \cite{ac}, in 
connection with the {\it Aharonov-Bohm (AB) effect} \cite{ab}. This
AB-phenomena got a lot of attention in the near past. A brief account
on this problem can be found in \cite{sz5}.

The 2D-Zeeman operators can directly be established by complex
Heisenberg group representations.
Complex Heisenberg groups have a natural
representation 
on the Hilbert space $\mathcal H=L_{\mathbf C\eta}^2$ 
of complex valued functions depending on both of the 
canonical coordinates
$(q_i,p_i)$, where the density
$\eta$, defining the inner product
$\langle f,g\rangle =\int f\overline g\eta dX$,
is 
$
\eta =e^{-|X|^2}
$
or, in general,
$
\eta_\lambda =e^{-\lambda |X|^2},
$
where $\lambda >0$ is a real constant. Even more general density is
introduced later.
Comparing with the real case, one should point out numerous other 
differences. 

In the complex case the canonical coordinates
are established such that $\partial_{p_i}=J(\partial_{q_i})$ holds. 
Thus, one can
define the holomorphic coordinates $z_i=q_i+\mathbf ip_i$. 
For recreating the whole Hilbert space,
$\mathcal H=L_{\mathbf C\eta}^2$, 
one should introduce also the antiholomorphic coordinates
$\overline z_i$. Then the density is, in the simplest case, 
of the form   
$
\eta =e^{-\sum z_i\overline{z}_i}
$,
whose most general form is
$
\eta_{\lambda_i} =e^{-\sum\lambda_i z_i\overline{z}_i}.
$
The Hilbert space is spanned by the polynomials written in terms
of both holomorphic and antiholomorphic coordinates. This Hilbert space
is isomorphic to the standard Hilbert space 
$L^2_{\mathbf C}$
having the standard Euclidean density $\eta =1$ by the map
\begin{equation}
L^2_{\mathbf C\eta} \to L^2_{\mathbf C}\quad ,\quad
\psi\to\psi e^{-{1\over 2}\sum\lambda_i z_i\overline{z}_i}.
\end{equation}

Differentiations 
$\partial_{z_i}\, ,\,\partial_{\overline{z}_i}$ 
are defined by means of partial differentiations 
$\partial_{q_i}$ and $\partial_{p_i}$ by
\begin{equation}
\partial_{z_i}={1\over 2}\partial_{q_i-\mathbf ip_i}\quad ,\quad
\partial_{\overline{z}_i}={1\over 2}\partial_{q_i+\mathbf ip_i}.
\end{equation}
Then the complex Poisson bracket is introduced by
\begin{equation}
\{f,g\}_{\mathbf C}=\sum
\partial_{z_i}(f)\partial_{\overline{z}_i}(g) -
\partial_{\overline{z}_i}(f)\partial_{z_i}(g).
\end{equation}

The complex Heisenberg algebra is defined by restricting this bracket
onto the linear space spanned by the functions
$\{z_i,\overline{z}_i,1\}$ (now the set of holomorphic and 
antiholomorphic
coordinates are extended by the constants $c\in\mathbf C$). 
The real Heisenberg algebra
is hidden inside of this complex algebra which can be 
uncovered by the formulas 
\begin{equation}
{q_i}={1\over 2}(z_i+{\overline{z}_i})\quad ,\quad
{p_i}=-{1\over 2}\mathbf i(z_i-{\overline{z}_i}).
\end{equation}

The representation of this complex Heisenberg algebra is introduced by
\begin{equation}
\label{chrep}
\rho_{\mathbf c} (z_i)(\psi )=(-\partial_{\overline{z}_i}+
\lambda_i z_i\cdot )\psi\quad ,\quad
\rho_{\mathbf c} (\overline{z}_i)(\psi )=\partial_{z_i}\psi . 
\end{equation}
This representation also satisfies the corresponding Heisenberg
relations (\ref{comrel}).
It should be emphasized, however, 
that this representation is not unitary on the
whole complex algebra. It becomes unitary by restricting it
onto the real sub-algebra. This statement follows from formulas
\begin{equation}
(-\partial_{\overline{z}_i}
+\lambda_i z_i\cdot )^*=\partial_{z_i}\quad ,\quad
(\partial_{z_i})^*=
-\partial_{\overline{z}_i}+
\lambda_i {z_i}\cdot \, , 
\end{equation}
since one has:
\begin{equation}
\rho_{\mathbf c} (q_i)=
{1\over 2}(\partial_{{z}_i}^*
+\partial_{{z}_i})\quad ,\quad
\rho_{\mathbf c} (p_i)= 
-{1\over 2}\mathbf i(\partial_{{z}_i}^*-\partial_{{z}_i}).
\end{equation}
Thus $e^{\mathbf i\rho_{\mathbf C}}$ defines, indeed, a unitary 
representation of the real subalgebra. 
    
The most important difference 
between the real and complex representations  
is that the complex representation $\rho_{\mathbf C}$ is not an
irreducible representation on the whole Hilbert space
$\mathcal H=L_{\mathbf C\eta}^2$. Indeed, the holomorphic subspace 
$\mathcal H^{(0)}$
spanned by the holomorphic polynomials 
$z_1^{a_1}\dots z_k^{a_k}$ 
is obviously invariant and irreducible under the actions of operators 
determined by the representation.
In the literature this irreducible representation is called 
complex (Fock) representation and, 
up-to the knowledge of this author, no
thorough investigation of the whole reducible representation has been
implemented so far. For the sake of clarity, we call it 
{\it extended Fock representation}. 
The 2D-Zeeman operator $H_Z$, described in (\ref{Zee}),
can be established from the Hamilton function $H_z$ by the
correspondence principle.
\medskip

\noindent{\bf Zones established by Gram-Schmidt orthogonalization.}
Next Hilbert space $\mathcal H$ is decomposed into the direct sum
of pairwise orthogonal semi-irreducible invariant subspaces. 
The main tool used in this decomposition is the
standard Gram-Schmidt orthogonalization. 
It turns out, later, that this decomposition is 
exactly the one established
by explicit spectrum computation such that the eigenfunctions are sorted
into zones according to their magnetic quantum states. 
In order to make clear distinction between these two constructions,
the latter's are called {\it spectrally defined Zeeman zones}, 
while the ones introduced in this chapter are the so called
{\it Heisenberg-Zeeman zones}.    

The {\it Gram-Schmidt process} is applied to the series 
$G^{(a)},\, a=0,1,\dots$,
of subspaces, where $G^{(a)}$ is spanned by the
subspaces 
$\overline{z}_1^{a_1}\dots \overline{z}_k^{a_{k/2}}\mathcal H^{(0)}$ 
satisfying
$a_1+\dots +a_{k/2}=a$. 
Clearly, $\mathcal H^{(0)}=G^{(0)}$ holds. 
The next subspace,
$\mathcal H^{(1)}$, is defined as orthogonal complement of 
$\mathcal H^{(0)}$ in
$G^{(0)}\oplus G^{(1)}$. 
The higher order subspaces are defined inductively. Thus,
$\mathcal H^{(a)}$ is defined as orthogonal complement of
$G^{(0)}\oplus\dots\oplus G^{(a-1)}$ in 
$G^{(0)}\oplus\dots\oplus G^{(a)}$. 

Since the subspaces 
$G^{(0)}\oplus\dots\oplus G^{(a)}$ 
are obviously invariant under the actions of operators defined in 
(\ref{chrep}), by induction and by the above unitary property one gets
that the {\it subspaces $\mathcal H^{(a)}$ are 
invariant with respect to the complex Heisenberg algebra 
representation}. The irreducibility 
is well known on the holomorphic zone.
If $k=dim(\mathbf v)=2$, 
the irreducibility can be established also on the
zones of higher order by the same simple proof. 
Whereas, in case of $k>0$, the higher order zonal
functions defined by a term 
$\overline{z}_1^{a_1}\dots \overline{z}_k^{a_{k/2}}$
defined with fixed exponents 
$a_1,\dots ,a_k$
form an irreducible invariant 
subspace denoted by 
$\mathcal H^{(a_1\dots a_{k/2})}\subset\mathcal H^{(a)}$. 

The subspaces $\mathcal H^{(a)}$ are called {\it gross zones}
on which the formulas appear in a much simpler form than on
the irreducible zones. Each gross zone decomposes into
$
{a+(k/2)-1\choose a}= 
{a+(k/2)-1\choose (k/2)-1}
$ 
number of irreducible zones. The gross zones are independent of
the complex coordinate system $(z_1,\dots ,z_{k/2})$, while
the irreducible zones do depend on them.
 
\medskip

\noindent{\bf Mathematical modeling; Zeeman manifolds.}
One of the most surprising features of the Zeeman operator $H_Z$ 
is that it can be pinned down
as the Laplacian on a Riemannian manifold. As far as the
author knows, this interpretation has not been recognized 
in the literature so far. We use this Riemannian manifold as the
fundamental mathematical model describing the space-time on 
quantum level. Although the metric is positive definite, 
this model fulfills the relativistic criteria \cite{sz5}.

By the first version of these manifolds single Zeeman-, or, 
Pauli-particles are modeled.
This fundamental Zeeman manifold is a Riemannian 
circle bundle over $\R^2$, defined by factorizing the center
of the 3-dimensional Heisenberg group endowed with a left-invariant 
metric. 

The Lie algebra $\mathbf n=\R^2\times\R=\R^3$ 
(where $\R$ is the center) of the 3D-Heisenberg
group can be described in terms of the
natural complex structure $J$, acting on $\R^2$, and the natural
inner product $\langle ,\rangle$, defined on $\mathbf n=\R^3$, 
by the formula
$\langle [X,Y],Z\rangle =\langle tJ(X),Y\rangle$, where the 3-vectors 
$X,Y$ and $Z$ are
in $\R^2$ and $\R$ respectively, furthermore, $t$ is
the coordinate of $Z$ in $\R$.
The map $Z\to tJ=J_Z$ associates skew endomorphisms acting on $\R^2$
to the elements, $Z$, of the center. They satisfy
the relation 
$J^2_Z=-|Z|^2id$. 
Thus the metric Lie algebra is completely
determined by the system 
\begin{eqnarray}
\label{n}
\{\mathbf n=\mathbf v\oplus\mathbf z,\langle ,\rangle ,J_Z\},
\label{nlie}
\end{eqnarray}
where $\mathbf v=\R^2$ and $\mathbf z$ are called X- and
Z-space respectively. With higher dimensional X- and Z-spaces this
system defines the Heisenberg type Lie algebras introduced by
Kaplan ~\cite{ka}. If the Clifford condition    
$J^2_Z=-|Z|^2id$
is dropped for the skew endomorphisms, the above system defines a
most general 2-step nilpotent Lie algebra. The considerations
will be extended to these general cases, however, the discussion
proceeds with the fundamental 3-dimensional case. 

Note that there
are two options, 
$J$ or $J^\prime =-J$, 
for choosing a complex structure on $\R^2$.
The two Lie algebras,
$\mathbf n$ and $\mathbf n^\prime$ are isometrically isomorphic by the
map $(X,Z)\to (X,Z^\prime =-Z)$.
 
The Lie group defined by this Lie algebra is denoted by $N$, 
furthermore, 
$g$ is the left-invariant extension of the inner product
$\langle ,\rangle$ defined on the tangent space 
$T_{(0,0)}(N)=\mathbf n$
at the origin. Then the exponential map is a one-to-one map 
whose inverse identifies
the group $N$ with its Lie algebra $\mathbf n$. Thus also the group  
lives on the same linear space $(X,Z)$ and the
group multiplication is given by:
\begin{eqnarray}
(X,Z)(X^* ,Z^* )
=(X+X^* ,Z+Z^* +{1\over 2}
[X,X^* ]).
\label{nilprod}
\end{eqnarray}

On the linear coordinate systems 
$\big\{x^1,x^2,t\big\}$, defined by the natural basis
$\big\{E_1, E_2, e_t\big\}$, the left-invariant
extensions of the vectors $E_i;e_t$ are of the form
\begin{eqnarray}
\mathbf
X_i=\partial_i + \frac {1} {2} 
\langle [X,E_i],e_t \rangle  \partial_t
=\partial_i + \frac {1} {2}  \langle
J\big(X\big),E_i\rangle
\partial_t\quad ; \quad
\mathbf T=\partial_t,
\label{invect}
\end{eqnarray}
where $\partial_i=\partial /\partial x^i$, $\partial_t=
\partial/\partial t$.
Then for the Laplacian, $\Delta$, acting on functions we have: 
\begin{eqnarray}
\label{lapl1}
\Delta=\Delta_X+ (1+\frac 1 {4}|X|^2)
 \partial_{tt}^2+\partial_t D \bullet,
\label{debullet}
\end{eqnarray}
where $\Delta_X$ is the Euclidean Laplacian on the X-space and
$D\bullet$ means differentiation (directional derivative)
with respect to the vector field
\begin{eqnarray}
D : X \to J\big (X\big )
\label{D}
\end{eqnarray}
tangent to the X-space \cite{sz2,sz3}. 

The above Laplacian is not the desired Zeeman operator yet. 
This surprising interpretation
can be established on center-periodic Heisenberg 
groups defined by an L-periodic lattice 
$\Gamma_Z=\{Z_{\gamma L}=\gamma L|\gamma\in\Z\}$ 
on the center.  
Since the $\Gamma_Z$ is a discrete subgroup of isometries, 
one can consider the
factor manifold $\Gamma_Z \backslash N$ with the factor metric. The
factor manifold is a principal circle bundle over the base space 
$\mathbf v$ such that the circles $C_X=\pi^{-1}(X)$ over the points 
$X\in \mathbf v$
are of constant length $L$. Then the projection 
$\pi : \Gamma_Z \backslash N \to \mathbf v$
projects the inner product from the horizontal subspace
(defined by the orthogonal complement to the circles) 
to the Euclidean inner product $\langle ,\rangle$ on the X-space.

By using the Fourier-Weierstrass decomposition
\begin{eqnarray}
L^2(\Gamma\backslash N)=\oplus FW^{(\gamma)},
\label{L^2}
\end{eqnarray}
where $FW^{(\gamma)}$ consists of functions 
of the form
\begin{eqnarray}
\phi^{(\gamma)}(X,Z)=
\varphi (X)e^{\mathbf i\gamma 2\pi t/ L},
\label{phi}
\end{eqnarray}
the Laplacian can be established in the following particular form. 

By (\ref{lapl1}), the function spaces 
$FW^{\gamma}$ are invariant under 
the action of the Laplacian. More precisely we have:
\begin{eqnarray}
\Delta \phi^{(\gamma )}
=(\Box_{(\lambda )}\varphi )e^{ 
\mathbf i\gamma 2\pi t/ L},\quad\quad\mbox{where}\\
\Box_{(\lambda )}
=\Delta_X +{2 \mathbf i} D_\lambda\bullet -
{4\lambda^2}\big(1 + 
\frac 1 {4} |X|^2\big)\quad ,\quad \lambda ={\pi\gamma\over L},
\label{lapl}
\end{eqnarray}
and $D_\lambda\bullet =\lambda D\bullet$ means directional derivative
along the the vector field $X\to\lambda J(X)=J_\lambda (X)$. 
If $\lambda <\, 0$,
the $J$ and $\lambda$ are exchanged for $-J$ and $-\lambda$ 
respectively. Thus one can assume that $\lambda >0$.

Apart from the constant term $-4\lambda^2$, operator
$\Box_{(\lambda )}$ is nothing but the Zeeman operator $H_Z$ 
described in (\ref{Zee}). 
The surplus constant term will be identified later
with the field-energy of the 
constant magnetic field in the charge spread, thus the above
operator is, actually, $H_{Zf}$ described earlier. The precise  
description of identification of the Zeeman operator (\ref{Zee}) 
with the Laplacian $\Box_\lambda$ acting on the invariant subspace 
$FW^{\gamma}$ is as follows. The macroscopic Zeeman operator is defined
by $\hbar =\mu =1$. Then 
$H_{Zf}=-(1/2)\Box_{(\lambda )}$, 
where $\lambda =-eB/2c$. Note that particles with negative charge 
correspond to the cases $\gamma >0$, i. e., they are attached to $J$, 
while particles with positive charges are attached to $-J$.

On the quantum (microscopic) level, the periodicity $L$
and the parameter $\lambda$ are exchanged for $L_\hbar =\hbar L$
and $\lambda_\hbar =\lambda /\hbar$ respectively. This process means
nothing but scaling of the periodicity by $\hbar$. Then we have
$H_{Zf}=-(\hbar^2/2\mu )\Box_{(\lambda_\hbar )}$. By scaling also the
Euclidean metric on the X-space by $\hbar/\sqrt\mu$, one has 
$H_{Zf}=-(1/2)\Box_{(\lambda_\hbar )}$. In the following we proceed
with the macroscopic operator, however, the previous formulas allow
an easy transfer from the macroscopic level to the microscopic one.

By a straightforward generalization, described later, these operators 
can be introduced on higher dimensional Heisenberg groups defined
by a complex structure, $J$, on an even dimensional Euclidean space
$\R^k$.
Let $(z_1,\dots ,z_{k/2})$, where 
$z_i=q_i+\mathbf ip_i$ and
$\partial_{p_i}=J(\partial_{q_i})$, 
be a complex coordinate system regarding $J$. Then this system 
identifies $\R^k$ with $\mathbf C^{k/2}$. The circle bundle,
defined by factorizing the center, $\R$, determines 
quantum operators depending just on one parameter $\lambda$. 
One can easily introduce operators depending on 
different parameters $\lambda_i>0$ defined for 
each complex coordinate plane $z_i$. 
(Such operators are constructed in the next section.)
These operators correspond to the Hamiltonians of systems 
where $k/2$ number of charged particles are circulating in a constant
magnetic field. When there is only one $\lambda$ involved, the particles
are considered to be identical up-to the sign of the charge.

\medskip

\noindent{\bf Zeeman manifolds with higher dimensional centers.}
The mathematical model for interpreting the Zeeman operator
as the Laplacian on a Riemannian manifold has been, so-far, 
a Riemannian circle bundle, defined by factorizing the centers on 
Heisenberg groups endowed with left invariant metrics.
This idea works out also on metric two-step nilpotent Lie groups 
whose center, $\mathbf z$, is factorized by a lattice $\Gamma_Z$. 
This center is considered as
an abstract higher dimensional space such that an element 
$Z\in\mathbf z$
is identified with the endomorphism $J_Z:\mathbf v\to\mathbf v$ and 
its natural inner product is defined by
$\langle Z_1,Z_2\rangle =-Tr(J_{Z_1}\circ J_{Z_2})$.
Formulas (\ref{n})-(\ref{lapl1})
apply also to these general cases, just the Laplacian (\ref{lapl1})
appears in a slightly different form. 
Up-to isomorphism, the Lie algebra of such a group is uniquely 
determined by a linear space, $J_{\mathbf z}$, of skew endomorphisms 
acting on the Euclidean space $\mathbf v$. Two 2-step nilpotent groups
are isometrically isomorphic if and only if the corresponding 
endomorphism spaces are conjugate. 

The rather large 
class of Riemannian torus bundles introduced in this way are  
called also Zeeman manifold. Below also particular 
Zeeman manifolds are introduced. It is
remarkable that for the so called Clifford-Zeeman manifolds even 
classification can be implemented. 
This classification can be used for classifying
the charged particles investigated in this theory. 
 
The Laplacian on the Riemannian group $(N_{J_{\mathbf z}},g)$, defined
by the endomorphism space $J_{\mathbf z}$, 
has the explicit form:

\begin{eqnarray}
\label{lapl2}
\Delta=\Delta_X+\Delta_Z+\frac 1 {4} \sum_{\alpha,\beta =1}^r 
\langle J_\alpha
\big (X\big),J_\beta\big (X\big)\rangle
 \partial_{\alpha\beta}^2
+\sum_{\alpha =1}^r\partial_\alpha D_\alpha \bullet,
\label{delta}
\end{eqnarray}
which leaves the function spaces $FW^{(\gamma )}$ spanned by
the functions of the form 
$\Psi^{(\gamma )}(X,Z)
=\psi (X)e^{2\pi\mathbf i\langle \mathcal Z_\gamma ,Z\rangle}
=\psi (X)e^{2\mathbf i\langle Z_\gamma ,Z\rangle}
$,     
for all lattice points 
$\mathcal Z_\gamma\in\Gamma_Z$ 
(resp. $ Z_\gamma\in \pi\Gamma_Z$), 
invariant.
Its action on such a 
function space can be described in the form
$\Delta (\Psi^{(\gamma )})(X,Z)=\Box_{(\gamma)}(\psi )(X)
e^{2\pi\mathbf i\langle\mathcal Z_\gamma ,Z\rangle}$,  
where operator $\Box_{(\gamma )}$, acting on 
$L^2(\mathbf v)$, is of the form
\begin{eqnarray}
\label{zeem}
\Box_{(\gamma )}
=\Delta_X + 2\pi\mathbf i D_{(\gamma )}\bullet -4\pi^2
\big(|\mathcal Z_\gamma |^2 + 
\frac 1 {4} |J_{\mathcal Z_\gamma}(X)|^2\big)\label{lapla} \\
=\Delta_X + 2\mathbf i D_{ Z_\gamma }\bullet -
4\big(| Z_\gamma |^2 + \frac 1 {4} |J_{Z_\gamma}(X)|^2\big).
\nonumber
\end{eqnarray}
 
Thus the Zeeman operator appears on the invariant
subspaces defined by the Fourier-Weierstrass decomposition. The 
spectral investigations on these manifolds are reduced to investigate
this operator on each Fourier-Weierstrass subspace separately. 

The particles represented
by these Riemannian torus bundles are called {\it Zeeman 
molecules}. A Zeeman molecule consists of $k/2$ number of 
charged particles. The higher dimensional center represents the variety
of the individual constant magnetic fields defined for the particles
in the molecule. It is explained in \cite{sz5} that the constant 
magnetic field fixes an inertia system which defines the self-time 
for a particle. Thus this general model matches Dirac's famous
multi-time model, introduced for implementing relativistic 
criteria on the quantum level. Lattice $\Gamma$ in the Z-space
relates our model to the crystal models of quantum theory.

There are special Z-molecules, defined by particular endomorphism 
spa\-ces, which are particularly interesting.
The {\it Heisenberg-type} or {\it Cliffordian
endomorphism spaces} are attached to Clifford modules 
(representations of Clifford algebras).
They are characterized by the property $J^2_Z=-|Z|^2id$, for all
$Z\in \mathbf z$, \cite{ka}. The corresponding molecules are called
{\it Clifford-Zeeman molecules}.
The well known {\it classification} of Clifford modules
provides classification also for the Clifford endomorphism
spaces and molecules. A brief account  on this 
classification theorem is as follows.

{\it If $r=dim(J_{\mathbf z})\not =3(mod 4)$, then 
there exist (up to equivalence) exactly one
irreducible H-type endomorphism space acting on a $\R^{n_r}$,
where the dimension $n_r$, depending on $r$, 
is described below. This endomorphism space
is denoted by $J_r^{(1)}$. If $r=3(mod 4)$, then there exist 
(up to equivalence) exactly
two non-equivalent irreducible H-type endomorphism spaces acting on
$\R^{n_r}$ which are denoted by 
$J_r^{(1,0)}$ and
$J_r^{(0,1)}$ 
respectively. They are connected by the relation 
$J_r^{(1,0)}\simeq
-J_r^{(0,1)}$.
 
The values $n_r$ corresponding to
$
r=8p,8p+1,\dots ,8p+7
$
are

\begin{eqnarray}
n_r=2^{4p}\, ,\, 2^{4p+1}\, , \, 2^{4p+2}\, , \,
2^{4p+2}\, ,
\, 2^{4p+3}\, ,\, 2^{4p+3}\, , \, 2^{4p+3}\, , \,
2^{4p+3}.
\label{cliff}
\end{eqnarray}

The reducible Clifford endomorphism spaces can be built up by these
irreducible ones. They are denoted by 
$J_r^{(a)}$ resp. $J_r^{(a,b)}$.
The corresponding Lie algebras are denoted by 
$h^{(a)}_r$ resp. 
$h^{(a,b)}_r$. In the latter case the X-space 
is defined by the $(a+b)$-times product 
$\R^{n_r}\times\dots\times\R^{n_r}$
such that on the last $b$ component the action of a $J_Z$ is defined by 
$J^{(0,1)}_Z$ 
and on the first $a$ components this action is defined by
$J^{(1,0)}_Z$. In the first case this process should be applied only on 
the corresponding $a$-times product.} 

In a Clifford endomorphism space each endomorphism anticommutes with
all perpendicular endomorphisms. In other words, all endomorphisms
are anticommutators. A more general concept can be 
introduced by  the {\it anticommutative endomorphism spaces} 
where all endomorphisms are anticommutators.
They can be built up, in a non-trivial way, by Clifford endomorphism
spaces. Roughly speaking, a CZ-molecule is the compound of
irreducible molecules of the same type while an {\it anticommutative
Z-molecule} is an indecomposable compound of CZ-molecules of 
different types in general.

Originally, the metric groups   
$(N_J,g)$ were used, in many different ways,  for constructing 
isospectral Riemannian metrics with different local geometries. The
author's results regarding such constructions are published in
\cite{sz1,sz2,sz3,sz4} which contain also
detailed history about this topic. 
These examples 
include isospectral pairs of metrics on
ball$\times$torus-, sphere$\times$to\-rus-, ball-, and sphere-type
manifolds. Among these examples the most striking are
those where one of the metrics in the isospectral pair
is homogeneous while the other is not even locally homogeneous.
Such examples have been constructed so-far on 
sphere-, sphere$\times$sphere, and sphere$\times$torus-type
manifolds. 
These isospectrality constructions are implemented such that on some
of the irreducible subspaces $\R^{n_r}$ the endomorphism spaces
$J_r^{(1,0)}$ 
(resp. $J_r^{(0,1)}$) are switched to 
$J_r^{(0,1)}$ 
(resp. $J_r^{(1,0)}$). It turns out that the Riemannian space, 
resulted by this switching, has a completely different local geometry,
yet, the considered domains in the original and the
new Riemann spaces are isospectral.  
Endomorphism spaces $J_r^{(1,0)}$ and $J_r^{(0,1)}$ are considered 
to be representing irreducible CZ-particles having opposite charges.
They are called also antiparticles. 
Thus the isospectrality theorem can be physically interpreted as 
follows: 

{\it By exchanging some of the irreducible CZ-particles 
in a CZ-molecule with their antiparticles the spectra of the 
considered domains 
remain the same, however, the local
geometry is drastically changed in general.}

Most of these isospectrality statements are established by constructing
intertwining operators, while some are proved by explicit
computations of the spectrum. These computations are
different from the one developed for the Zeeman zones.

\section{\bf Zonal projections, zonal point-spreads}

{\bf Establishing the projection kernels.}
In this section the  operators $\mathbf P^{(a)}$ projecting 
onto the zones $\mathcal H^{(a)}$ will be explicitly described. 
If $\{\varphi_i^{(a)}\}_{i=1}^\infty$ is an 
orthonormal basis in $\mathcal H^{(a)}$, 
then the corresponding projection can be
formally defined as convolution with the kernel 
\begin{equation}
P^{(a)}(z,w)=\sum_i \varphi^{(a)}_i(z)\overline{\varphi^{(a)}_i(w)}.
\end{equation}
Interestingly enough, these self-adjoint operators are integral
operators having a smooth  
Hermitean integral kernel 
$P^{(a)}(z_j,w_j)$. The projection kernel regarding 
the holomorphic zone is
the well known Bergman kernel. However, up-to the knowledge 
of this author, the projections onto the other zones
have never been considered in the literature so far.  
First, we describe these kernels for the 2D-Zeeman operator defined 
on the complex plane
$\mathbf C$ (i. e., $k=2$). Note that in this case 
only a single parameter 
$\lambda$ is involved to the computations.

On the zone $\mathcal H^{(0)}$ the holomorphic polynomials $z^i$ form an
orthogonal basis, thus the holomorphic polynomials 
\begin{equation}
\label{phi^0}
\varphi_i^{(0)}(z)=\mathbf z^i=(1/\int (z\overline{z})^i
d\eta )^{1\over 2}z^i
\end{equation}
form an orthonormal basis in $\mathcal H^{(0)}$. 
The sought integral operator
(convolution operator) on this subspace is 
\begin{eqnarray}
\label{P^0}
P^{(0)}(z,w)*_wf(w)=\int P^{(0)}(z,w)f(w)d\eta (w)=\\
\int {1\over \pi}e^{z\overline{w}}f(w)e^{-w\overline{w}}dw
=\int {\lambda\over \pi}e^{\lambda z\overline{w}}f(w)
e^{-\lambda w\overline{w}}dw,\nonumber \\
\mathbf P^{(0)}_\lambda (f)(Z)=\int {\lambda\over \pi}
e^{\lambda (\langle Z,W\rangle +\mathbf i\langle Z,J(W)\rangle )}f(W)
e^{-\lambda |W|^2}dW,
\end{eqnarray}
where $Z,W$ denotes the the corresponding complex numbers as vectors
and $J$ is the complex structure on this vector space.
(The latter formula is used for an easy generalization to the higher
dimensional cases. If the real dimension is $k$ and there is a single
$\lambda$ involved, only the constant should be changed to 
$(\lambda /\pi )^{k/2}$.) Formula (\ref{P^0}) can be established
by showing that the image space of this Hermitean integral 
operator is
$\mathcal H^{(0)}$ and also $P^{(0)}(z,w)*_ww^i=w^i$ holds.

Now we pass to the next zone, $\mathcal H^{(1)}$. Also this case is
considered on the complex 
plane $\mathbf C$ first. Now the functions
\begin{equation}
\phi^{(1)}_i(z)=
\overline{z}z^i
-\mathbf P^{(0)}(\overline{z}z^i)=
\overline{z}z^i
-(1/\lambda)iz^{i-1}
\end{equation}
form an orthogonal basis on the zone and the normalization provides
an orthonormal basis. Also in this case the 
projection 
$P^{(1)}$ is an integral
operator whose kernel can be constructed as follows. 

First note that
this projection can be described by the generalized kernel 
\begin{equation}
\label{prefor1}
 P^{(1)}_{zw}=(\overline{z}P^{(0)}_z)*_w\partial_{\overline{w}} 
-P^{(0)}_z*_v(\overline{v}P^{(0)}_v)*_w\partial_{\overline w}
\end{equation}
since the range of this Hermitean operator
is $\mathcal H^{(1)}$ which 
acts on this subspace as the identity. Also this 
operator, which still contains differentiations, is an
integral operator with a smooth Hermitian
kernel. This statement follows by the following computations:
\begin{eqnarray}
 P^{(1)}_{zw}=
(\overline{z}(w-z)P^{(0)}_z)*_w 
-(\overline{w}P^{(0)}_z)*_w\partial_{\overline w}=\\
((\overline{z}(w-z)+1 +\overline{w}(z-w))P^{(0)}_z)*_w 
=((1 -|z-w|^2)P^{(0)}_z)*_w\, ,\nonumber \\ 
\mathbf P^{(1)} (f)_Z=\int {1\over \pi}(1-|Z-W|^2)
e^{(\langle Z,W\rangle +
\mathbf i\langle Z,J(W)\rangle )}f(W)e^{-|W|^2}dW.
\end{eqnarray}
The derivative $\partial_{\overline{w}}$, 
acting on functions $f$ first, is 
handled by integration by parts such that
it is substituted by its dual's action on the function 
$
(\overline{z}P^{(0)}_z) 
-P^{(0)}_z*_v(\overline{v}P^{(0)}_v).
$

Passing to the second zone, $\mathcal H^{(2)}$, the same 
computational technique, applied in a 
more complicated situation, yields the following formulas:
\begin{eqnarray}
\label{(2)}
2!\mathbf P_z^{(2)}(f)
=\big((\overline{z}^2P^{(0)}_z)*_w\partial^2_{\overline{w}} \\
-(P^{(0)}_z+P^{(1)}_z)*_v(\overline{v}^2P^{(0)}_v)*_w
\partial^2_{\overline w}\big)(f)
=\{\overline{z}^2(w-z)^2P_z^{(0)}*_w\nonumber\nonumber \\
-\overline{v}^2(w-v)^2(2-(z-v)(\overline{z}-\overline{v}))P_z*_vP_v*_w\}
(f)\nonumber \\
=\overline{z}^2(w-z)^2P_z^{(0)}*_w
-\partial_{v}^2\{(w-v)^2(1-(z-v)(\overline{z}-\overline{v}))\nonumber \\
-2(z-v)(w-v)\}P_z*_vP_v*_w=\{\overline{z}^2(w-z)^2+(w-z)^2\overline{w}^2
\nonumber \\
-2\overline{w}\overline{z}(w-z)^2-4(z-w)(\overline{z}-\overline{w})+2\}
P^{(0)}_z*_w(f)\nonumber \\
=\int {1\over \pi}(2-4|Z-W|^2+|Z-W|^4)
e^{Z\cdot\overline{W}
}f(W)e^{-|W|^2}dW.\nonumber
\end{eqnarray}
Term $Z\cdot\overline W$ often appears later in the form
$
Z\cdot\overline{W}=\langle Z,W\rangle +\mathbf i\langle Z,J(W)\rangle .
$

Regarding the polynomials arising in these formulas,
note that  
$L_0^{(0)}(t)=1, L^{(0)}_1(t)=1-t,L^{(0)}_2(t)=(2-4t+t^2)/2!$
are the first three Laguerre polynomials parametrized by $\alpha =0$ 
and one gets the 
above ones by the substitution $t=|Z-W|^2$.
These polynomials 
depend, in general, on a parameter $\alpha$
\cite{sze}, such that the $n^{th}$-order 
polynomials $L_{n}^{(\alpha )}$ are defined by 
the eigenfunctions 
of the operator
\begin{equation}
\Lambda_{\alpha} (u)(t)=tu^{\prime\prime}+(\alpha +1-t)u^\prime.
\end{equation}
The corresponding eigenvalue is $-n$. The above polynomials
belong to parameter $\alpha =0$. On the $k$-dimensional case this
parameter will be $\alpha =(k/2)-1$.

An analogous definition for these polynomials uses 
the following Rodrigues formula:
\begin{equation}
e^{-t}t^{\alpha} L_a^{(\alpha )}(t)={1\over a!}
\partial^a_t(e^{-t}t^{a+\alpha}).
\end{equation}
The well known explicit form of these polynomials is
\begin{equation}
\label{rec0}
L^{(\alpha )}_a(t)=\sum_{b=0}^a{a+\alpha\choose a-b}{(-t)^b\over b!}.
\end{equation}
The followings recurrence formulas will be often used in establishing 
the general projection kernels.
\begin{eqnarray}
\label{rec1}
(a+1)L^{(\alpha )}_{a+1}(t)=(\alpha +2a+1-t)L_a^{(\alpha )}(t)-
aL_{a-1}^{(\alpha )}(t)\, ,\\
\label{rec2}
\sum_{b=0}^aL_b^{(\alpha )}=L_a^{(\alpha +1)}(t)\, ,\,
\partial_tL_a^{(\alpha )}=
-L_{a-1}^{(\alpha +1)}\, , \\
\label{rec3}
tL_a^{(\alpha )\,\prime}(t)=aL_a^{(\alpha )}(t)
-(a+\alpha )L_{a-1}^{(\alpha )}(t)\, .
\end{eqnarray}
 
The pre-formula  (\ref{prefor1}) for defining the
projections can be recurrently introduced for the higher order cases.  
Indeed, if $P^{(0)},\dots ,P^{(a-1)}$ 
are already defined, the pre-formula for $P^{(a)}$ is:
\begin{eqnarray}
\label{rec4}
\mathbf P^{(a)}
(f)_z={1\over a!}\big( (\overline{z}^aP^{(0)})*_w 
- (P_z^{(0)}+\dots +P_z^{(a-1)})*_v(\overline{v}^aP^{(0)})*_w\big)
\partial^a_{\overline{w}}(f_w)
\end{eqnarray}
Strictly speaking, these formulas concern the case 
$k=2\, ,\,\lambda=1$.
If for a general $k$; but still assuming
$\lambda=1$; the corresponding projections with respect to a
complex coordinate $z_i$ are denoted by 
$
\mathbf P^{(a_i)},
$ 
then the sought pre-formula is:
\begin{equation}
\label{rec5}
\mathbf P^{(a)}=\sum_{a_1+\dots +a_{k/2}=a}
\mathbf P^{(a_1)}\dots
\mathbf P^{(a_{k/2})}.
\end{equation}
The long computation in (\ref{(2)}) shows that a 
direct establishing of the
Laguerre polynomial form of the projection operators would be an arduous
task. This is why mathematical induction has been chosen below 
for concluding the proof of the corresponding theorem.
\begin{theorem}{\bf (Zonal Theorem)}
Representation (\ref{chrep}) is a reducible 
unitary representation of the
complex Heisenberg algebra. The function space $L^2_{\mathbf C\eta}$
is an orthogonal direct sum of invariant subspaces $\mathcal H^{(a)}$, 
called gross Zeeman zones. In case $k=2$, this representation is 
irreducible on each zone. It is always irreducible also on the 
holomorphic zones, however, it  is reducible if $k>2$ and $a>0$. 
(The irreducible zones are considered in the end of this theorem.)

One can construct the gross zone decomposition by 
the Gram-Schmidt orthogonalization of subspaces 
$G^{(a)}=\overline{z}_1^{a_1}\dots 
\overline{z}_k^{a_k}\mathcal H^{(0)}$, 
where $a_1+\dots +a_k=a$. The first zone 
$\mathcal H^{(0)}=G^0$ is the holomorphic 
zone spanned by the holomorphic polynomials 
$z_1^{a_1}\dots z_k^{a_k}$.

The projections onto the gross zones are integral operators
having analytic integral kernels. In case of $k=2\, ,\,\lambda =1$
these integral operators are
\begin{equation}
\label{P^a}
\mathbf P^{(a)}(f)_{Z}
=\int {1\over \pi}L_a^{(0)}(|Z-W|^2)
e^{\langle Z,W\rangle +\mathbf i\langle Z,J(W)\rangle}f(W)e^{-|W|^2}dW,
\end{equation}
where $L_a^{(0)}(t)$ is the corresponding Laguerre polynomial.
In the $k$-dimension\-al cases, where there is only a single parameter
$\lambda >0$ involved, the corresponding integral operators are
\begin{equation}
\label{P^a_l}
\mathbf P_\lambda^{(a)}(f)_{Z}
=\int {\lambda^{k/2}\over \pi^{k/2}}L_a^{((k/2)-1)}(\lambda |Z-W|^2)
e^{\lambda Z\cdot\overline{W}}f(W)e^{-\lambda |W|^2}dW,
\end{equation}
where
$Z\cdot\overline{W}=\langle Z,W\rangle +
\mathbf i\langle Z,J(W)\rangle$. 
If more $\lambda_i$'s are involved which are 
defined for $k_i$-dimensional
subspaces such that $k=\sum_{i=1}^p k_i$ holds then the projection is: 
\begin{eqnarray}
\label{P^a_l.l}
\mathbf P_{\lambda_1\dots \lambda_p}^{(a)}(f)_{Z}=\\
\int {\prod\lambda_i^{k_i/2}\over \pi^{k/2}}L_a^{((k/2)-1)}
(\sum\lambda_i |Z_i-W_i|^2)
e^{\sum\lambda_i Z_i\cdot\overline{W}_i}f(W)
e^{-\sum\lambda_i |W_i|^2}dW.\nonumber
\end{eqnarray}

Each gross zone decomposes into
$
{a+(k/2)-1\choose a}= 
{a+(k/2)-1\choose (k/2)-1}
$ 
number of irreducible zones. The gross zones are independent of
the complex coordinate system $(z_1,\dots ,z_{k/2})$, while
the irreducible zones do depend on them.
The irreducible zones, denoted by 
$\mathcal H^{(a_1\dots a_{k/2})}\subset\mathcal H^{(a)}$, 
are spanned by the higher order zonal
functions defined by the antiholomorphic terms 
$\overline{z}_1^{a_1}\dots \overline{z}_k^{a_{k/2}}$
for fixed exponents 
$a_1,\dots ,a_k$. The projection,
$
\mathbf P_{\lambda_1\dots \lambda_p}^{(a_1\dots a_{k/2})}(f)_{Z},
$
onto an irreducible zone differs from the previous operator just in the 
Laguerre polynomial term what should be exchanged for
$
\prod_{i=1}^{k/2}
L^{(0)}_{a_1}(\lambda_i |Z_i-W_i|^2).
$
 
\end{theorem}
\heading{Proof.}
The proof goes through the cases
indicated in the theorem step by step. 
If $k=2,\lambda =1$, then (\ref{P^a}) is proved
for $a\leq 2$ by (\ref{phi^0})-(\ref{(2)}). 
Suppose that it is valid for $a\geq 2$. 
Then for its validity regarding $a+1$ one should prove that 
the integral operator defined by the kernel
\begin{equation}
\mathbf L^{(0)}_{a+1}(z,w)=
{1\over \pi}L^{(0)}_{a+1}(|z-w|^2)
e^{z\overline{w}}
\end{equation}
is nothing but the projection onto $\mathcal H^{(a+1)}$.
 
The range of this Hermitean operator is in 
$\mathcal H^{(\leq a+1)}
=G^{(0)}\oplus \dots\oplus G^{(a+1)}$, therefore, it maps the 
perpendicular complement
$\mathcal H^{(\leq a+1)\perp}$
to zero. Thus, one should consider this operator only on
$\mathcal H^{(\leq a+1)}$.
Therefore, it is enough to prove that it maps also the functions 
$
\overline{w}^bw^i
\in \mathcal H^{(\leq a)}$, where $b\leq a$, to zero.
By (\ref{rec1}), (\ref{rec3}), induction, and integration by parts 
one can carry out the following computations.
\begin{eqnarray}
\int (a+1)\mathbf L^{(0)}_{a+1}(z,w)
(\overline{w}^bw^i)d\eta\\
=\int\big((2a+1-|z-w|^2)\mathbf L_a^{(0)}(z,w)-
a\mathbf L_{a-1}^{(0)}(z,w)\big)
(\overline{w}^bw^i)d\eta
\nonumber\\
=\int\big((2a+1)
\mathbf L_a^{(0)}(z,w)
-a\mathbf L_{a-1}^{(0)}(z,w)\big)
(\overline{w}^bw^i)d\eta\nonumber \\
-\int\partial_w(|z-w|^2
\mathbf L_a^{(0)}(z,w)
(\overline{w}^{b-1}w^i))d\eta
\nonumber \\
=\int\big(a
\mathbf L_a^{(0)}(z,w)
(\overline{w}^bw^i)
-a\overline{z}\mathbf L_{a-1}^{(0)}(z,w)
(\overline{w}^{b-1}w^i)\big)
d\eta\nonumber \\
-\int (-a\mathbf L^{(0)}_{a-1}(z,w)\overline w^{b-1}
\partial_ww^i+|z-w|^2
\mathbf L_a^{(0)}(z,w)
\partial^b_w(w^i))d\eta =\phi (z)
\, .\nonumber
\end{eqnarray}

Since the formula is established for $a=0,1,2$, we can suppose 
$a\geq 2$.
Then the integral of the very last function in the last line is zero, 
since it refers to the $\mathbf P^{(a)}$-projection of a function from 
$\mathcal H^{\leq 1}$ onto $\mathcal H^{(a)}$. Thus the
the $\phi$ decomposes into three integral terms,
$\phi =\phi_1 +\phi_2 +\phi_3$,
where $\phi_3$ is derived from the last line and the other two from
the line before the last one.

If $b<a$, then $\phi =0$, since the integral vanishes by terms,
due to the induction. This means that the 
$P^{(a+1)}=\mathbf L^{(0)}_{a+1}(z,w)*$ 
operates trivially on 
$\mathcal H^{(\leq a-1)}$.
 
If $a=b$, then $\phi$ is in
$\mathcal H^{(\leq a-1)}$. In fact, the $\phi_3$ is obviously in that
subspace, furthermore, the first two terms corresponding to
$\phi_1 +\phi_2$ is a
difference of functions from $\mathcal H^{(\leq a)}$ 
such that they have the same 
leading term, $a\overline{z}^az^i$. 
Thus also this function must be in
$\mathcal H^{(\leq a-1)}$. 
Therefore, by applying 
$P^{(a+1)}=\mathbf L^{(0)}_{a+1}(v,z)*$ 
on $\phi$, one gets the zero function:
$\mathbf (P^{(a+1)})^2(\overline w^bw^i)=0$.
Since the operator considered is Hermitean, also the $\phi$
must be zero then. This proves the desired identity.

Proof of (\ref{P^a}) will be completed by 
showing that the integral operator
with kernel $\mathbf L^{(0)}_{a+1}(z,w)$ acts on $G^{(a+1)}$ 
as the projection $\mathbf P^{(a+1)}$. By the induction hypothesis,
we have. 
\begin{eqnarray}
\int\mathbf L_{a+1}^{(0)}(z,w)
\overline{w}^{a+1}w^id\eta
=
\overline{z}^{a+1}z^i 
-\sum_{b=0}^a\mathbf P^{(b)}
(\overline{w}^{a+1}w^i)
\\
=
\overline{z}^{a+1}z^i
-\int\sum_{b=0}^a
\mathbf L_b^{(0)}(z,w)
\overline{w}^{a+1}w^id\eta
=\nonumber \\
\overline{z}^{a+1}z^i
-\int
\mathbf L_a^{(1)}(z,w)
\overline{w}^{a+1}w^id\eta .\nonumber
\end{eqnarray}
The last equation is due to (\ref{rec2}). By (\ref{rec1}), 
(\ref{rec3}), and the first part of this proof we have
\begin{eqnarray}
\int\mathbf L_{a+1}^{(0)}(z,w)
\overline{w}^{a+1}w^id\eta
=
\int\partial_w(\mathbf L_{a+1}^{(0)}(z,w)
\overline{w}^{a}w^i)d\eta
=\\
\int\mathbf L_{a}^{(1)}(z,w)(\overline{z}-\overline{w})
\overline{w}^{a}w^id\eta
=
\overline{z}^{a+1}z^i
-\int
\mathbf L_a^{(1)}(z,w)
\overline{w}^{a+1}w^id\eta .\nonumber
\end{eqnarray}
This concludes the proof of (\ref{P^a}).

Formula (\ref{P^a_l}), concerning the general 
$k$-dimensional cases, will be
pro\-ved, first, for $\lambda =1$. In light of 
(\ref{rec4}) and (\ref{rec5}) the
statement immediately follows from
\begin{equation}
\label{comp}
L_n^{(\alpha )}(t_1+\dots +t_{\alpha +1})=
\sum_{a_1+\dots +a_{\alpha +1}=n}
L^{(0)}_{a_1}(t_1)\dots L^{(0)}_{a_{\alpha +1}}(t_{\alpha +1}),
\end{equation}
which can be readily established by the following simpler formula:
\begin{equation}
\label{simp}
L_n^{(\alpha +1)}(t_1+t_2)=
\sum_{a=0}^{n}L^{(0)}_{a}(t_1) 
L^{(\alpha )}_{n-a}(t_2).
\end{equation}
The latter formula can be established by (\ref{rec0})
or by mathematical induction. Below the latter method is chosen,
since this was the natural way in which 
(\ref{P^a})-(\ref{P^a_l.l}) were found.
In establishing (\ref{simp}) the induction is used twice. 
First, in proving
of the starting step
\begin{equation}
\label{L^1}
L^{(1)}_n(t_1+t_2)=\sum L^{(0)}_a(t_1)L^{(0)}_{n-a}(t_2),
\end{equation}
where the induction
concerns $n$. Then the induction concerns $\alpha $ in proving of 
(\ref{simp}).

Since 
$L_1^{(0)}(t_1)+L^{(0)}_1(t_2)=2-(t_1+t_2)$, 
formula (\ref{L^1}) is valid for $n=1$. Assume that it is valid for
$degrees\leq (n-1)$. 
Then, by applying (\ref{rec1})
both for $\alpha =1$ and $\alpha =0$, one has
the following sequence of equations
{\small
\begin{eqnarray}
nL_n^{(1)}(t_1+t_2)=
(2n-(t_1+t_2))L_{n-1}^{(1)}(t_1+t_2)
-nL^{(1)}_{n-2}(t_1+t_2)=
\\
\sum_{a=0}^{n-1} ((2a+1-t_1)+(2(n-1-a)+1-t_2))
L_{a}^{(0)}(t_1)L_{n-1-a}^{(0)}(t_2)-
\nonumber \\
-nL^{(1)}_{n-2}(t_1+t_2)
=\sum_{a=0}^{n-1} \big(
(a+1)L_{a+1}^{(0)}(t_1)L_{n-1-a}^{(0)}(t_2)+
(n-a)L_{a}^{(0)}(t_1)L_{n-a}^{(0)}(t_2)\nonumber \\
+aL_{a-1}^{(0)}(t_1)L_{n-1-a}^{(0)}(t_2)+
(n-1-a)L_{a}^{(0)}(t_1)L_{n-a-2}^{(0)}(t_2)\big)
-nL^{(1)}_{n-2}(t_1+t_2)=\nonumber \\
=n\sum_{b=0}^{n}
L_{b}^{(0)}(t_1)L_{n-b}^{(0)}(t_2)+
n\sum_{c=0}^{n-2}
L_{c}^{(0)}(t_1)L_{n-c-2}^{(0)}(t_2)
-nL^{(1)}_{n-2}(t_1+t_2)=\nonumber \\
=n\sum_{b=0}^{n}
L_{b}^{(0)}(t_1)L_{n-b}^{(0)}(t_2).\nonumber
\end{eqnarray}
}
The terms in the first and second summation formulas are equal 
by the recurrence equation
\begin{equation}
(c+1)L_{c+1}^{(0)}(t)=
(2c+1-t)L_{c}^{(0)}(t)
-cL_{c-1}^{(0)}(t).
\end{equation}

To complete the proof of (\ref{simp}) 
suppose that it is valid for $\alpha $.
Then its validity for $\alpha +1$ follows by using  
(\ref{rec2}) twice:
\begin{eqnarray}
L_n^{(\alpha +1)}(t_1+t_2)=
\sum_{a=0}^nL_a^{(\alpha )}(t_1+t_2)
=\sum_{a=0}^n\sum_{p=0}^aL_p^{(\alpha -1)}(t_1)L^{(0)}_{a-p}(t_2)\\
=\sum_{b=0}^n\big(\sum_{a=b}^nL_{a-b}^{(\alpha -1)}(t_1)\big)
L^{(0)}_{b}(t_2)
=\sum_{b=0}^nL_{n-b}^{(\alpha )}(t_1)
L^{(0)}_{b}(t_2).\nonumber
\end{eqnarray}
Formulas (\ref{P^0}) and (\ref{comp}) 
establish the most general case described in (\ref{P^a_l.l})
by the substitutions $t_i=\lambda_i|Z_i-W_i|^2$.
Finally, the statement concerning the projection kernels regarding
the irreducible zones can be established by (\ref{comp}). These 
arguments complete the proof of the theorem.

$\Box$
\medskip

\noindent{\bf Point-spreads.}
Since the zonal projection kernels are the restrictions of the global 
Dirac delta distribution 
$\delta_Z(W)=\sum \varphi_i(Z)\overline{\varphi}_i(W)$,
they are denoted also in the following form 
\begin{eqnarray}
\label{pointspread}
\delta_{\lambda_1\dots \lambda_{k/2}}^{(a)}(Z,W)=
P_{\lambda_1\dots \lambda_{k/2}}^{(a)}(Z,W)=
\\ 
{\prod\lambda_i^{k_i/2}\over \pi^{k/2}}L_{(k/2)-1}^{(a)}
\big(\sum\lambda_i |Z_i-W_i|^2\big)
e^{\sum\lambda_i (Z_i\cdot\overline{W}_i-{1\over 2}(|Z_i|^2+ |W_i|^2)}.
\nonumber
\end{eqnarray}
The kernels
$
\delta_{\lambda_1\dots \lambda_{k/2}}^{(a_1\dots a_{k/2})}(Z,W)
$
on the irreducible zones 
$\mathcal H^{(a_1\dots a_{k/2})}\subset\mathcal H^{(a)}$
differ from these ones just in the 
Laguerre polynomial term what should be exchanged for
$
\prod_{i=1}^{k/2}
L^{(0)}_{a_1}(\lambda_i |Z_i-W_i|^2).
$

As opposed to the global Dirac-delta distribution kernel, these ones
are smooth functions. They represent one of most important 
fundamental concepts
in this article.
Also note the simple rule by which the zonal kernels are derived
from the one defined for the holomorphic zone. This holomorphic 
point-spread is just
multiplied by the corresponding radial Laguerre polynomials.  
These zonal kernels can be interpreted such that 
a point particle concentrated
at a point $Z$ appears on the zone as an object which spreads around 
$Z$ as a wave packet with wave-function 
described by the above kernel explicitly. 

There is a complete matching
between the point-spreads and de Broglie's wave packets (cf.
\cite{bo}, pages 61), which concept became one of the corner stones
of Quantum Theory.
De Broglie's theory  was finalized by the Schr\"odinger
equation. The mathematical formalism did not follow this development,
however, and the Schr\"odinger theory is built up on such mathematical
background which does not exclude the existence of the controversial
point objects. On the contrary, 
an electron must be considered as a point-object in
the Schr\"odinger theory as well 
(cf. Weisskopf's argument on this problem in \cite{schw, sz5}). 
An other demonstration for 
the presence of
point particles in classical theory is the duality principle, 
stating that objects manifest themselves sometime as waves and sometime
as point particles. The bridge between the two visualizations is
built up in Born's probabilistic theory,
where the probability for that that a particle, 
attached to a wave $\xi$, can be found on a domain $D$ is measured by
$\int_D\xi\overline\xi$.

These controversial point-objects, by having infinite self-mass or
self-charge attributed to them by the Schr\"odinger, launched one
of the deepest crisis-es in the history of physics.
In the zonal theory de Broglie's idea is established on a mathematical
level. Although the points are ostracized from this theory, 
the point-spreads 
still bear some reminiscence of the point-particles. For instance, they
are the most compressed wave-packets and all the 
other wave-functions in the 
zone can be expressed as a unique superposition of the point-spreads.
If $\xi$ is a zone-function, the above integral 
measures the probability for that that
the center of a point-spread is on the domain $D$. This 
interpretation restores, in some extend, the duality principle in the
zonal theory.

Function 
$
\delta^{(a)}_{\lambda Z}
\overline{\delta}^{(a)}_{\lambda Z}
$
is called the density of the spread around $Z$. By this reason, 
function 
$
\delta^{(a)}_{\lambda Z}
$
is called spread-amplitude. Both the spread-amplitude and 
spread-density generate well defined measures on the path-space
consisting of continuous curves connecting two arbitrary points.
Both measures can be constructed by the method applied in constructing
the Wiener measure.
 
The point-spread 
concept bears some remote reminiscence of Heisenberg's 
suggestion (1938) for the existence of 
a fundamental length $L$, analogously
to $h$, such that field theory was valid only for distances larger than
$L$ and so divergent integrals would be cut off at that distance.
This idea has never became an effective theory, however. 
Other distant relatives of the point-spread concept are the 
smeared operators, i. e. those suitably averaged over small 
regions of space-time, considered by Bohr and Rosenfeld in quantum 
field theory. There are also other theories where an electron is 
considered to be extended. Most of them fail on lacking the
explanation for the question: 
Why does an extended electron not blow up? The zonal theory is checked
against this problem in \cite{sz5}, section 
(F) ``Linking to the blackbody radiation; Solid zonal particles''.

\section{\bf Spectral definition of zones}

{\bf Technicalities.} In this chapter the zones 
are established by the explicit 
spectra of the Zeeman-Laplacians 
defined on center-periodic 2-step nilpotent manifolds 
$\Gamma\backslash N$.
On these manifolds the Laplacian leaves the subspaces defined by 
the Fourier-Weierstrass decomposition
\begin{equation}
L^2(\Gamma\backslash N)=\oplus FW^{(\gamma)}
\end{equation}
invariant. The $FW^{(\gamma)}$ is defined 
for any fixed lattice point 
$\mathcal Z_\gamma\in \Gamma_Z$ (resp. $Z_\gamma\in \pi\Gamma_Z$) 
such that it consists the functions of the form
\begin{equation}
\label{fi_inFW}
\phi^{(\gamma)}(X,Z)=
\varphi (X)e^{2\pi\mathbf i\langle\mathcal Z_\gamma ,Z\rangle}=
\varphi (X)e^{2\mathbf i\langle Z_\gamma ,Z\rangle}.
\end{equation}
(Lattice $\pi\Gamma_Z$ is 
introduced, mostly, for simplifying the formulas.) Note 
that (\ref{fi_inFW}) defines $\Gamma_Z$-periodic functions.
Thus, this problem is ultimately reduced to computing the spectrum 
of the Zeeman Laplacian $\Box_{(\gamma )}$ 
introduced in (\ref{zeem}). The computations are 
carried out by the second
form of this operator, where the $\pi$ is involved in the norm of 
$Z_\gamma$ 
($|Z_\gamma |=\pi |\mathcal Z_\gamma|$). Since the Hamiltonian 
$H_{Z_\gamma f}=-{1\over 2}\Box_\gamma$ differs from the classical
Zeeman-Hamilton operator $H_{Z_\gamma}$ 
just by an additive constant, this
spectrum computation readily provides the spectrum also of a free 
electron orbiting around the origin in a constant magnetic field.

The subspace corresponding to $\gamma =0$ is nothing but the space of 
functions depending only on the $X$-variable. Thus the corresponding
Laplacian is the Euclidean Laplacian $\Delta_X$ and  
the spectrum is the well known continuous spectrum there.
On the other subspaces, where $\gamma\not =0$ and the endomorphism 
$J_{\gamma}$ is non-degenerated,
the spectrum is discrete which
can explicitly be computed. If the $J_\gamma$ is degenerated, then the 
Laplacian is the direct sum of the Euclidean Laplacian, 
defined on the maximal subspace where $J_\gamma$ is degenerated, and of 
the non-degenerated magnetic Laplacian defined on the complement 
subspace. Thus the spectrum is the product of the continuous 
Euclidean spectrum by the discrete spectrum of the non-degenerated
operator. By writing the corresponding decomposition of the 
Fourier-Weierstrass function spaces in the form
\begin{equation}
FW^{\gamma}=
FW_0^{\gamma}\odot
FW_{ML}^{\gamma},
\end{equation} 
one has the decomposition
\begin{equation}
FW=
FW_0\odot
FW_{ML}=
(\oplus FW_0^{\gamma})\odot
(\oplus FW_{ML}^{\gamma})
\end{equation}
of the $\Gamma_Z$-periodic functions such that the spectrum is 
continuous on $FW_0$ (it is determined by the corresponding
Euclidean spectra) and it is discrete on $FW_{ML}$. Note that on
H-type groups the subspace $FW_0$ is nothing but $FW^0$ and $FW_{ML}$
is spanned by all $FW^\gamma$ where $\gamma\not =0$. In what follows
the spectrum is computed on $FW_{ML}$ what is called {\it magnetic
spectrum}, or, {\it Zeeman spectrum}. It simply means that
$J_{\gamma}$ is non-degenerated in the Laplacian (\ref{zeem}).

Despite the discreteness, this spectrum
is very different from those defined on compact manifolds. For instance,
each eigenvalue has infinite multiplicity. Roughly speaking, this is 
one of the reasons for infinities appear in free electron theory. This
paper approaches to this problem by decomposing the corresponding 
$L^2$-space into {\it Zeeman-zones}
on which these multiplicities together with other important quantities
become finite ones.

The differential operator (\ref{zeem}) is built up by the following 
two commuting operators
\begin{equation}
\mathbf O=
\Delta_X -
4\big(|Z_\gamma |^2 + \frac 1 {4} |J_{Z_\gamma}(X)|^2\big)\, ,\,
\mathbf D=2\mathbf i D_{(\gamma )}\bullet.
\end{equation}
Apart from the constant term $4|Z_\gamma |^2$, the first operator  
can be interpreted as the Schr\"odinger wave operator 
of a harmonic oscillator. The textbook solution of the 
eigenvalue problem of this operator uses
Hermitian polynomials. We adopt this technique to the spectral
computations developed in this paper. It should 
be mentioned, however, that 
these functions provide eigenfunctions only for the first operator
above. In order to get eigenfunctions also with respect to the
second operator (which commutes with the first one!), a
certain decomposition technique on 
Hermitian polynomials shall be applied. 

This computation technique
can be described, in short, as {\it splitting the spectral lines of the 
harmonic oscillator operator $\mathbf O$ with the angular 
momentum operator $\mathbf D$ (Zeeman effect).} 
There is explained later that
this is indeed different from the standard method, which neglects
the oscillator potential, retains the Coulomb potential and the 
spectrum is computed by splitting the spectrum of the 
Coulomb-Schr\"odinger operator by $\mathbf D$. It turns out that
the latter standard method is ``blindfolded'' regarding the 
Zeeman zones, meaning that the zones can not be constructed by it.
\medskip

\noindent{\bf Eigenfunctions of the harmonic oscillator operator} 
$\mathbf O$.
The standard form of the Hamilton operator for a 
harmonic oscillator is
$\Delta_X-|X|^2$, which acts on $L^2(\mathbf R^k)$ by the Friedrichs
extension. The 1-dimensional 
version of this operator is $S=d^2/dx^2 -x^2$, whose eigenvalues
can be computed as follows \cite{ta}. 
(We need this standard computations later, for
establishing the spectrum of the complete operator $\Box_\gamma$.)

By introducing the first order operators
\begin{equation}
s^+={d\over dx}+x\quad ;\quad s^-={d\over dx}-x\, ,
\end{equation}
we have
\begin{equation}
S=s^+s^-+1=s^-s^+-1\,\, ,\,\, [s^+,s^-]=s^+s^--s^-s^+=-2\, .
\end{equation}

Since $s^+(e^{-x^2/2})=0$, 
the $L^2$-eigenfunctions of $S$ are the functions
\begin{equation}
\big({d\over dx}-x\big)^le^{-x^2/2}=H_l(x)e^{-x^2/2}\, ,
\end{equation}
where $H_l(x)$ are the so called Hermite polynomials, given by
\begin{equation}
H_l(x)=(-1)^re^{x^2}\big({d\over dx}\big)^le^{-x^2}
=\sum_{j=0}^{[l/2]}(-1)^j{l!\over j!(l-2j)!}(2x)^{l-2j}\, .
\end{equation}
The computations below, implemented in the general cases, show that
the eigenvalue of $S$, corresponding to the function (\ref{K}), is
$-(2l+1)$. Since the functions
\begin{equation}
c_lH_l(x)e^{-x^2/2}\,\, ;\,\, c_l=\big(\pi^{1/2}2^l(l!)\big)^{-1/2}
\,\, ;\,\, l=0,1,\dots ,
\end{equation}
form an orthonormal basis on $L^2(\R)$, the eigenvalue problem
of $S$ is solved in the 1-dimensional case.

On $\R^k$ 
the eigeinfunctions of $S=\Delta_X-|X|^2$ are  
\begin{equation}
\label{h_l..l}
h_{l_1\dots l_k}(X)=
H_{l_1}(X^1)\dots H_{l_k}(X^k)e^{-|X|^2/2}\,\, ;\,\, l=l_1+\dots +l_k
\end{equation}
with eigenvalues $-(2l+k)$. In this case the functions
$c_{l_1}\dots c_{l_k}h_{l_1\dots l_k}(X)$ form an orthonormal basis
in $L^2(\R^k)$.

Since the function 
\begin{equation}
\label{K}
\langle -J^2_{Z_\gamma}(X),X\rangle=\langle K(X),X\rangle
\end{equation}
differs from $|X|^2$ in general, one should modify the above formulas,
since the corresponding spectrum, 
defined in (\ref{h_l..l}), 
depends also on the eigenvalues of
the positive self-adjoint operator $K$. 
In fact, let $\sqrt K$ be the uniquely 
determined positive square root of
the operator $K$. If $\sqrt K$ has a single eigenvalue $\lambda$
then the corresponding eigenfunctions are
\begin{equation}
\label{eigfh}
h_{\lambda ,l_1,\dots ,l_k}(X)=
H_{l_1}(\lambda X^1)\dots H_{l_k}(\lambda X^k)e^{-\lambda |X|^2/2}
\end{equation}
with eigenvalues $-(2l+k)\lambda$. The multiplicity of this eigenvalue
is ${{l+k-1}\choose{k-1}}$, since this is the number of options for 
$l$ can be expressed
in the form $l=l_1+l_2+\dots +l_k$, where $l_i$'s
are non-negative integers. (A reformulation of this 
combinatorial question is this:
In how many ways can be $l$ identically looking balls painted by
$k$ different colors?)

If $\sqrt K$ has the distinct
eigenvalues $0\leq \lambda_1<\dots <\lambda_s$ with multiplicities
$k_1,\dots ,k_s$ on the corresponding $k_i$-dimensional subspaces
$\mathcal X_i$ then the eigenfunctions have the following product
form
\begin{equation}
\label{eigfH}
{\small\prod_{i=1}^s}H_{r_{i1}}(\lambda_i X_i^1)\dots H_{r_{ik_i}}
(\lambda_i X_i^{k_i})
e^{-\lambda_i |X_i|^2/2},
\end{equation}
where   
$X_i=(X_i^1,\dots ,X_i^{k_i})$ 
is a vector from $\mathcal X_i$ whose components are defined 
by an orthonormal basis on $\mathcal X_i$.
The eigenvalues corresponding to these functions are
\begin{equation}
{\small \sum_{i=1}^s}-\lambda_i(k_i+2{\small\sum_{j=1}^{k_i}}l_{ij}).
\end{equation}
The multiplicity of this eigenvalue is the product,
{\small
\begin{equation}
\prod_{i=1}^m{{k_i-1+ \sum_{j=1}^{k_i}l_{ij}}\choose{k_i-1}},
\end{equation}
}
of the corresponding multiplicities determined on the eigensubspaces
of $K$.

\medskip
\noindent{\bf Eigenfunctions of the complete Laplacian 
$\Box_\gamma$.}
Functions (\ref{eigfh}) and (\ref{eigfH}) 
are not yet the eigenfunctions of the
operator $\mathbf iD_{\gamma}\bullet$, which is also included into
the complete magnetic Laplacian (\ref{zeem}). 
Since this operator commutes with $\mathbf O$, 
one can generate a subspace which is invariant
with respect to the complete operator by letting the
operator $D_{\gamma}\bullet$ successively acting on functions 
(\ref{eigfH}). 
Next we describe this generated subspace which shall be
decomposed into the direct sum of subspaces which are invariant with
respect to the action of the magnetic Laplacian. 
Exactly this step in the computations is called
{\it splitting the spectrum of $\mathbf O$ with the angular momentum
operator.} 

The main idea in these
computations is that the linear functions 
$\langle Q+\mathbf iJ(Q),X\rangle$
resp. 
$\langle Q-\mathbf iJ(Q),X\rangle$,
where the unit skew endomorphism $J$ satisfying
$J^2=-id$ is defined by normalizing $J_\gamma$, 
are eigenfunctions of the corresponding operator
$D\bullet$, defined by $J$, with eigenvalue $\mathbf i$ resp.
$-\mathbf i$. Thus the eigenvalue of $D\bullet$ 
on the $l^{th}$ order polynomial $P_+^{(p)}(X)P_-^{(l-p)}(X)$, where
\begin{eqnarray}
\label{P+-}
P_+^{(p)}(X)=\langle Q^+_{1}+\mathbf iJ(Q^+_{1}),X\rangle\dots
\langle Q^+_{p}+\mathbf iJ(Q_p^+),X\rangle,\\
P_-^{(l-p)}(X)=\langle Q^-_{1}-\mathbf iJ(Q^-_{1}),X\rangle\dots
\langle Q^-_{l-p}-\mathbf iJ(Q^-_{l-p}),X\rangle \nonumber
\end{eqnarray}
is $(p-(l-p))\mathbf i=(2p-l)\mathbf i$. 
This formula, where $p=0,\dots ,l$,
describes the complete set of eigenvalues of $D\bullet$
on the $l^{th}$ order polynomials.

First we deal with the case when $\sqrt K$ has only a single
eigenvalue and, thus, the eigenfunction of $\mathbf O$ 
has the form (\ref{eigfh}). Let
$(E_1,\dots ,E_k)$ be the orthonormal basis on the X-space
defining the coordinates $(X^1,\dots ,X^k)$ and write
each coordinate in the form
\begin{equation}
X^i={1\over 2}\langle E_i+\mathbf iJ(E_i),X\rangle
+{1\over 2}\langle E_i-\mathbf iJ(E_i),X\rangle.
\end{equation}
Plug in these formulas in (\ref{eigfh}) and perform all of the 
powering and multiplications indicated
in the Hermite polynomials. Then the polynomial
$H_{l_1}(\lambda X^1)\dots H_{l_k}(\lambda X^k)$
becomes a linear combination of polynomials of the form (\ref{P+-}).

For a fixed $m$, the function 
$H_{\lambda ,l_1,\dots ,l_k}^{(m)}(X)$ is 
defined by the sum of all those
terms for which $2p-l=m$ holds. 
Note that these functions always
contain an $l^{th}$-order term and, except for the functions 
$H_{\lambda ,l_1,\dots ,l_k}^{(l)}$ and
$H_{\lambda ,l_1,\dots l_k}^{(-l)}$, 
they always include also lower order terms. In other words, these
polynomials are inhomogeneous and, therefore,
they are different from the spherical harmonics used for similar
spectral computations in general. 
(See more about this question in the last section of this chapter.) 
The complete range of possible $m$'s is described by $2p-l$, where
$p=0,\dots ,l$. Such a function is an eigenfunction of the operator
$\mathbf D=\mathbf i D\bullet$ 
with eigenvalue $-(2p-l)\lambda$. Therefore,
by letting $\mathbf D$ successively acting on
\begin{equation}
\label{H}
H_{\lambda ,l_1,\dots ,l_k}(X)=
\sum_{p=0}^lH_{\lambda,l_1,\dots ,l_k}^{(2p-l)}(X),
\end{equation}
we have
\begin{equation}
\label{DH}
\mathbf D^i(
H_{\lambda ,l_1,\dots ,l_k})(X)=\sum_{p=0}^l
(-(2p-l)\lambda )^iH_{\lambda,l_1,\dots ,l_k}^{(2p-l)}(X).
\end{equation}

Since 
$\mathbf De^{-\lambda |X|^2/2}=0$, thus 
(\ref{H}) and (\ref{DH}) hold also
for the functions
\begin{equation}
\label{He}
h_{\lambda ,l_1,\dots ,l_k}(X)
\,\, ,\,\, 
h_{\lambda,l_1,\dots ,l_k}^{(2p-l)}(X)=
H^{(2p-l)}_{\lambda ,l_1,\dots ,l_k}(X) 
e^{-\lambda |X|^2/2} .
\end{equation}

Note that  
equations (\ref{DH}), considered for the values
$i=0,\dots ,p$, allow to express the functions
$H_{\lambda,l_1,\dots ,l_k}^{(2p-l)}$
via the functions
$\mathbf D^i
H_{\lambda ,l_1,\dots ,l_k}$. 
In fact, these equations establish
a linear connection between the two fields of $(l+1)$-vectors whose
coordinate functions 
are the above functions. The
matrix of this linear transform is non-degenerated. In fact, if $l$
is an odd number then it is a Vandermonde matrix,
$V$, with entries $V_{ip}=
(-(2p-l)\lambda )^i\not =0$. If $l$ is even, the first row consists of
$1$'s and except the first one all the elements are $0$'s 
in the middle column. The rest part is again a Vandermonde matrix.
In either cases the $V$ is invertible and
\begin{equation}
H^{(2p-l)}_{\lambda ,l_1,\dots ,l_k}=\sum_{i=0}^rV^{-1}_{pi} 
\mathbf D^i
H_{\lambda ,l_1,\dots ,l_k}.
\end{equation}
The same equations hold corresponding to the functions (\ref{He}). 
By the commutativity of 
operators
$
\mathbf O=\Delta_X   
-4\big(|Z_\gamma |^2 + \frac 1 {4} |J_{Z_\gamma}(X)|^2\big)
$
and
$
2\mathbf D
$ 
we get that the functions
\begin{equation}
h_{\lambda,l_1,\dots ,l_k}^{(2p-l)}(X)=
H^{(2p-l)}_{\lambda ,l_1,\dots ,l_k}(X) 
e^{-\lambda |X|^2/2}
\end{equation}
are eigenfunctions of the complete Zeeman-Laplace operator 
with eigenvalue
\begin{equation}
\label{eigv1}
-((4p+k)\lambda +4k\lambda^2).  
\end{equation}
The first term was computed as the sum of 
$-(2l+k)\lambda$ and $-2(2p-l)\lambda$.

Next the multiplicities of eigenvalues are determined for a fixed  
value $l$. One should bear in mind this restriction in
these computations. The  multiplicities
on the whole $L^2(\mathbf v)$ will be discussed later. 
 
For a fixed $l$ the computation of multiplicities can be established by
noticing that the spectrum-elements (\ref{eigv1}) depend only on the 
holomorphic parameter $p$ and they do not depend on the 
anti-holomorphic parameter $\upsilon =l-p$. Therefore, the complex
valued $l^{th}$-order polynomial space splits into 
$(l+1)$ number of
invariant subspaces determined by the values 
$p=0,1,\dots ,l$. For a fixed
$p$ the complex dimension of the corresponding subspace is
\begin{equation}
\label{mult1}
{{p+k/2-1}\choose{k/2-1}}{{\upsilon +k/2-1}\choose{k/2-1}},
\end{equation}
which number is the multiplicity of the eigenvalue corresponding to 
$p$. (This multiplicity is computed on the complex space 
$\mathbf C^{k/2}$. The first term in (\ref{mult1}) is the number of 
decompositions $p=l_1+\dots +l_{k/2}$. The second factor has the 
same meaning corresponding to $\upsilon$.)
This shows that the anti-holomorphic parameter, 
$\upsilon =l-p$, plays role only in the
multiplicity of a spectrum-element.

When $\sqrt K$ has distinct eigenvalues, $\lambda_i\,;\,i=1,\dots ,s$,
then the eigenfunctions are the products of functions
\begin{equation}
h_{\lambda_i,l_{i1},\dots ,l_{ik_i}}^{(2p_i-\sum_jl_{ij})}(X)=
H^{(2p_i-\sum_jl_{ij})}_{\lambda_i ,l_{i1},\dots ,l_{ik_i}}(X) 
e^{-\lambda_i|X_i|^2/2}, 
\end{equation}
where $X_i$ represents a generic vector from the eigensubspace 
$\mathcal X_i$ belonging to the
eigenvalue $\lambda_i$. Thus we have
\begin{theorem}{\bf (Explicit Spectral Theorem)}
Laplacian $\Delta_N$ leaves the function spa\-ces $FW^{(\gamma)}$,
defined for the lattice points $\mathcal Z_{\gamma}\in \Gamma$ by
functions (\ref{fi_inFW}), invariant and it acts on such a
subspace as the magnetic Laplacian (\ref{zeem}). 
This latter operator is the sum of the harmonic oscillator
operator $\mathbf O$ and of the angular momentum operator 
$\mathbf D_{(\gamma)}$. These two operators commute.

For $\gamma =0$, the Laplacian $\Box_{(0)}$ is nothing but
the Euclidean Laplacian $\Delta_X$
which has the well known continuous spectrum on $FW^{(0)}$. For a
$\gamma \not =0$, where the $J_{\gamma}$ is non-degenerated, the
spectrum is discrete
on $FW^{(\gamma)}$. More precisely,
the eigenfunctions 
can be expressed in the product form
\begin{equation}
\prod_{i=1}^mh_{\lambda_{\gamma i},l_{i1},\dots ,l_{ik_i}}^
{(2p_i-\sum_jl_{ij})}(X_i)
\end{equation}
which have the corresponding eigenvalues
\begin{equation}
\label{eigv2}
\sum_{i=1}^s-(
\lambda_{\gamma i}
(4p_i +k_i)+4\lambda_{\gamma i}^2k_i).
\end{equation}
For fixed values
$p_i\, ,\,\upsilon_i =\sum_j l_{ij}-p_i$, 
where $i=1,2,\dots ,s$,
the multiplicity (the complex dimension of the corresponding 
eigensubspace) of such an eigenvalue is
\begin{equation}
\label{mult2}
\prod_{i=1}^s{{p_i+k_i/2-1}\choose{k_i/2-1}}
{{\upsilon_i+k_i/2-1}\choose{k_i/2-1}}.
\end{equation}
 
If the $J_\gamma$ is degenerated, the spectrum is continuous which is,
formally, the product of the continuous spectrum corresponding to the
subspace where $J_\gamma$ is degenerated and of the discrete spectrum
corresponding to the subspace where $J_\gamma$ is non-degenerated
and where the spectrum can be computed by the above method. On 
Heisenberg-type groups the only
subspace on which the spectrum is continuous is $FW^{(0)}$ and the
spectrum is discrete on the other subspaces. 
\end{theorem}
\medskip

\noindent{\bf Zones defined spectrally.}
The gross Zeeman zones were defined, in Section 2, as the
semi-irreducible invariant subspaces of the complex Heisenberg
group representation. There is also described the number of irreducible
zones a gross zone decomposes into. 
In this section the zones are defined by means of the
explicit spectrum established above.
If there is only a single parameter 
$\lambda\not =0$ involved, this decomposition is defined, in 
this new way, by  the eigensubspaces spanned by functions  
$H^{(m)}_{\lambda,\mathbf l}$, defined in (\ref{H}). In the following
these subspaces are denoted also by $\mathbf H^{(p,l-p)}$.
These polynomials are constructed by the $l^{th}$-order polynomials 
$P^{(p,l-p)}$, where $p$ resp.
$(l-p)$ are the holomorphic resp. antiholomorphic degrees, by a
splitting technique using Vandermonde matrices. The orders $p$
and $\upsilon =l-p$ are called {\it holomorphic} and 
{\it antiholomorphic indexes} also with respect to the eigenfunctions.
The connection between the pairs $(m,l)$ and
$(p,l-p)$ can be established by the formula $m=2p-l$.
The elements of the spectrum are 
$(4p+k)\lambda +4k\lambda^2$. 

When more $\lambda$'s,
$(\lambda_1 ,\dots ,\lambda_s)$,
are involved the corresponding subspaces are denoted by
\begin{equation}
\mathbf H^{(p_1,\dots ,p_s,\upsilon_1,\dots ,\upsilon_s)}\quad ,\quad
\upsilon_i=l_i-p_i,
\end{equation}
which can be considered as the tensor product of the function spaces
$
\mathbf H_{i}^{(p_i,\upsilon_i)}
$.

By formulas (\ref{eigv1}) and (\ref{eigv2}), 
the actual eigenvalues do not depend
on the antiholomorphic indexes, therefore, each eigenvalue is constant
on the infinite dimensional function space
\begin{equation}
IM^{(p_1,\dots ,p_s)}=\sum_{(\upsilon_1,\dots ,\upsilon_s)}
\mathbf H^{(p_1,\dots ,p_s,
\upsilon_1,\dots ,\upsilon_s)}e^{-{1\over 2}\sum\lambda_i|X_i|^2}.
\end{equation}
In other words, each eigenvalue has infinite multiplicity. 

For fixed antiholomorphic degrees 
${(\upsilon_1,\dots ,\upsilon_s)}$
consider the direct sum: 
\begin{equation}
\label{FSGZ}
FSGZ^
{(\upsilon_1,\dots ,\upsilon_s)}
=\sum_
{(p_1,\dots ,p_s)}
\mathbf H^{(p_1,\dots ,p_s,
\upsilon_1,\dots ,\upsilon_s)}e^{-{1\over 2}\sum\lambda_i|X_i|^2},
\end{equation}
which defines a so called Fine Spectral Gross Zone
corresponding to the Fine Zone Indexes
$FZI={(\upsilon_1,\dots ,\upsilon_s)}$. 
The Spectral Gross Zone, 
$SGZ^\upsilon$, 
corresponding to the index
$GZI=\upsilon =\upsilon_1+\dots +\upsilon_s$ 
is the direct sum
\begin{equation}
SGZ^\upsilon =\sum_
{\upsilon =\upsilon_1+\dots +\upsilon_s} 
FSGZ^
{(\upsilon_1,\dots ,\upsilon_s)}.
\end{equation}

In order to define the Irreducible Spectral Zones, consider a complex
coordinate neighborhood $z_1\dots ,z_{k/2}$ of the
complex structure $J$ defined by normalizing $J_\gamma$ such that
the $J_\gamma$ has the eigenvalue $\lambda_i$
on a complex plane $z_i$. 
If $p_i,l_i,u_i=l_i-p_i$ are given for all complex plane $z_i$, then
the Irreducible Spectral Zone 
$
ISZ^
{(u_1,\dots ,u_{k/2})}
$ 
regarding the index 
$
{(u_1,\dots ,u_{k/2})}
$
is defined in the same way as the fine zone (\ref{FSGZ}).

For an eigenfunction the {\it gross holomorphic, antiholomorphic,
zonal azimuthal, and magnetic indexes (quantum numbers)} are defined by 
$p=\sum p_i,\\ u=\sum u_i,l=\sum l_i$ and $m=2p-l$ respectively.
For a fixed complex plane $z_i$ the corresponding indexes are called
individual indexes. These names indicate that each coordinate represent
an individual charged particle.
Then we have:
\begin{theorem}{\bf (Splitting Theorem)}
The gross spectral zones are invariant under the action 
of the magnetic Laplacian (\ref{zeem}). 
If $\lambda_i\not =0$, for all $i$,
the spectrum restricted to a gross zone is discrete such that 
each element has finite multiplicity. 

The spectrum-elements on two different gross zones
are the same. The multiplicities, 
described in (\ref{mult2}) for fixed fine zone index
$FZI={(\upsilon_1,\dots ,\upsilon_s)}$, 
are different, however. Two spectra yielding this property are called
isochromatic. If $k=2$, any two gross zones are isospectral. If
$k>2$, any two distinct gross zones are properly isochromatic having
distinct multiplicities for the same eigenvalues.  

By the formula ${1\over 2}(l-m)=u$,
the gross zonal azimuthal number and the magnetic number of the
eigenfunctions uniquely
determine the gross zone they belong to.

Each spectral gross zone decomposes into
$
{\upsilon +(k/2)-1\choose \upsilon}= 
{\upsilon+(k/2)-1\choose (k/2)-1}
$ 
number of irreducible spectral zones, 
$
ISZ^
{(u_1,\dots ,u_{k/2})}
$. The irreducible zone to which an eigenfunction belongs to 
is uniquely determined by the individual magnetic and zonal 
azimuthal numbers of the eigenfunction.
The multiplicity of each eigenvalue on an irreducible zone is $1$ 
if and only if the eigenvalues $\lambda_i$ are distinct.  
Thus higher zonal multiplicities indicate the presence of
particles which are identical up to the sign of the charge.
Any two irreducible spectral zones are isospectral.
\end{theorem}

Finally, the identity of the earlier introduced and the 
spectral zones is established.
\begin{theorem}
The spectral gross zones are 
semi-irreduc\-ible invariant subspaces 
under the action of the complex Heisenberg algebra representation
\begin{equation}
\rho_{\mathbf c} (z^i)(\psi )=
(-\partial_{\overline{z}^i}+
\lambda_i z^i\cdot )
\psi\quad ,\quad
\rho_{\mathbf c} (\overline{z}^i)(\psi )=\partial_{z^i}\psi 
\end{equation}
defined in (\ref{chrep}) on the complex Hilbert space 
$L^2_{\mathbf C\eta_{\lambda_i}}$, where
$\eta_{\lambda_i}=e^{-\sum\lambda_iz^i\overline{z}^i}$. 

Actually, the spectral zones (gross and irreducible) are nothing 
but the ones defined earlier by the extended Fock representation.
\end{theorem}

\heading{Proof.} 
For a fixed zonal azimuthal quantum number $l$, 
the spectral Zeeman zones are distinguished by
the distinct eigenvalues $m=2p-l$ of the skew operator $\mathbf D$.
Therefore two distinct zones are perpendicular to each other.
Operators $\partial_{z_i}$ act, by definition, on 
function spaces
$\mathbf H_{\gamma}^{(p_1,\dots ,p_s,
\upsilon_1,\dots ,\upsilon_s)}$ 
and thus the spectral Zeeman zones are
invariant under these actions. Since the operators
$-\partial_{\overline{z}_i}+
\lambda_i z_i\cdot$
and $\partial_{z_i}$ are dual, 
the invariance with respect to both operators
follows by the above orthogonal property.

The first spectral zone is obviously $\mathcal H^{(0)}$. 
The second one is
in $G^{(0)}\oplus G^{(1)}$ such that it is perpendicular to $G^{(0)}$. 
Therefore, it must be $\mathcal H^{(1)}$. 
Mathematical induction proves that
the zones introduced by the extended Fock representation 
are the same as the spectral zones.
$\Box$
\medskip

\noindent{\bf Standard computation of the Zeeman spectrum.}
In order to establish a
precise connection between the above and classical
quantum numbers of particle theory, 
in this section the spectral computations are accomplished in
the standard way they are carried out
in classical quantum mechanic for computing the spectrum
of an electron. These classical computations always involve the Coulomb
potential of the nucleus.
It should be noted that the zones are
not invariant with respect to multiplications with radial functions, 
thus the technique establishing the Zeeman zones can not be applied 
to the Coulomb-Schr\"odinger operator. The only
method for computing the spectrum in this case is the one
which traces back the eigenvalue
problem to the radial functions 
$f(\langle X,X\rangle )$ by seeking the eigenfunctions in the form
$fH$, where the $H$ is a homogeneous harmonic polynomial. 

This method is really different from the technique of 
splitting the spectrum of
$\mathbf O$ with $\mathbf D_\gamma$,
by which the Zeeman zones are constructed.
On the complex plane, for instance, the homogeneous 
harmonic polynomials can be divided into two classes. One of them
contains the holomorphic polynomials, the other the antiholomorphic
ones. By this classification the zones can not be discovered. 
Probably this is why the Zeeman
zones were not recognized in the literature earlier. 

This standard method is a variant of the technique 
of {\it splitting the spectrum of
the Coulomb-Schr\"odinger operator with the angular momentum operator,}
by which the Zeeman effect (quantization of magnetic dipole moment)
was originally explained.
The word ``variant'' is justified, since 
the standard computations
usually neglect the quadratic oscillator potential and retain 
the Coulomb 
potential $V$, while in our case $V=0$, however, both the
harmonic oscillator potential and the field energy is retained.
 
The computations will be carried out for such a
Zeeman operator $H_{Zf}$ which depends only on a single parameter 
$\lambda$. The eigenfunctions of this operator are sought
now in the form
$F(X)=f(\lambda \langle X,X\rangle )H^{(\tilde l,m )}(X)$, 
where $f$ is a real valued smooth even function 
defined on $\R$ and
$H^{(\tilde l,m )}(X)$ is a complex valued
homogeneous harmonic polynomial of order
$\tilde l$ on the X-space such that it is eigenfunction 
also of the operator
$D\bullet$, introduced by means of the unit endomorphism $J$ in \S 2,
with eigenvalue $m =2\tilde p-\tilde l$. When more parameters, 
$\{\lambda_i\}$, are involved, the eigenfunctions are represented
as product of functions of the form
$F_{(i)}(X)=f_{(i)}(\lambda_i \langle X,X\rangle)
H_{(i)}^{(\tilde l_i,m_i )}(X)$,
where the functions in the formula are defined on the maximal
eigensubspace corresponding to the parameter $\lambda_i$.

One can trace back 
the eigenvalue problem of $\Box_{\lambda}$ to the eigenvalue problem
of an ordinary differential operator acting on the radial functions
$f (\langle X,X\rangle)$ as follows.  
First the simplest case is described
where $\lambda =1$.
Since $D_{\lambda }\bullet f=0$, $|Z_\lambda |^2=k$, and 
$|J_\lambda (X)|^2=\langle X,X\rangle$, we get
\begin{eqnarray}
(\Box_{(\lambda )}F)(X)=\big(4\langle X,X\rangle f^{\prime\prime}
(\langle X,X\rangle )
+(2k+4\tilde l)f^\prime (\langle X,X\rangle )\\
-(2m +4(k +{1\over 4}\langle X,X\rangle )
f(\langle X,X\rangle ))\big)H^{(\tilde l,m)} (X).
\end{eqnarray}
The eigenvalue problem is reduced, therefore, to the ordinary
differential operator
\begin{equation}
\label{Lf}
(L_{(\lambda =1,l,m )}f)(t)=
4tf^{\prime\prime}(t)
+(2k+4\tilde l)f^\prime (t)
-(2m +4(k +
{1\over 4}t))f(t).
\end{equation}
 
The function 
$e^{-{1\over 2}t}$ 
is an eigenfunction of this
operator with eigenvalue $-(4\tilde p+5k)$. The general
eigenfunctions are sought in the form
\begin{equation}
f(t)=u(t)e^{-{1\over 2}t}.
\end{equation}
Such a function is an eigenfunction of 
$L_{\tilde l,m}$
if and only if the $u(t)$ is an eigenfunction of the differential 
operator
\begin{equation}
\label{Pop}
(P_{(\lambda =1,\tilde l,m )}u)(t)=
4tu^{\prime\prime}(t)
+(2k+4\tilde l-4t)u^\prime (t)
-(4\tilde p+5k)u(t).
\end{equation}

The most remarkable property of this operator is that it has a
uniquely determined polynomial eigenfunction
\begin{equation}
\label{lag}
u_{(\lambda =1,n,\tilde l,m )}(t)
=t^n+a_1t^{n-1}+a_2t^{n-2}+\dots +a_{n-1}t+a_n
\end{equation}
with coefficients satisfying the recursion formulas
\begin{equation}
a_0=1\quad ,\quad a_i=-a_{i-1}(n-i)(n+\tilde l+{1\over 2}k+1-i)n^{-1}.
\end{equation}
One can easily establish explicit combinatorial formula for $a_i$ by
this recursion, yet we do not deal with this problem here. 
The eigenvalue corresponding to this polynomial is
\begin{equation}
\label{leigv}
\mu_{(\lambda =1,n,\tilde l,\nu )}=-(4n+4p+5k)\quad ,
\quad p={1\over 2}(m +\tilde l).
\end{equation}

The polynomials (\ref{lag}) are nothing but 
the Laguerre polynomials. It was mentioned earlier that
they are defined as the $n^{th}$-order polynomial
eigenfunctions
of the differential operator
\begin{equation}
\label{Lop}
\Lambda_{\alpha} (u)(t)=tu^{\prime\prime}+(\alpha +1-t)u^\prime ,
\end{equation}
with eigenvalues $-n$.
Therefore
\begin{equation}
P_{(\lambda =1,\tilde l,m )}=4\Lambda_{({1\over 2}k+
\tilde l-1)}-(4p+5k).
\end{equation}
Thus the eigenfunctions of operators (\ref{Pop}) and (\ref{Lop}) 
are the same.
The relations between the corresponding eigenvalues are properly
described in (\ref{leigv}). 
Particularly we get that, for fixed values for $k,\tilde l,m$ and thus
also for $p$, the functions 
$
u_{(\lambda =1,n,\tilde l,m )}\, ,\, n=0,1,\dots \infty
$ form a basis in $L^2([0,\infty ))$.

In case of a single $\lambda$, differential operator (\ref{Lf})
is of the form 
\begin{equation}
(L_{(\lambda ,\tilde l,m )}f)(t)=
4tf^{\prime\prime}(t)
+(2k+4\tilde l)f^\prime (t)
-(2m\lambda +4\lambda^2(k +
{1\over 4}t))f(t).
\end{equation}
The eigenfunctions are sought in the form
\begin{equation}
u(\langle X,X\rangle )e^{-{1\over 2}\lambda 
\langle X,X\rangle}\Phi_{(\tilde l,m)}(X),
\end{equation}
which leads to the radial operator 
\begin{equation}
\label{rad}
(P_{(\lambda ,\tilde l,m )}u)(t)=
4tu^{\prime\prime}(t)
+(2k+4\tilde l-4\lambda t)u^\prime (t)
-(2m\lambda +4k\lambda^2)u(t),
\end{equation}
corresponding to (\ref{Pop}). 
The polynomial eigenfunctions relate to the above Laguerre polynomials
by the simple formula
\begin{equation}
u_{(\lambda ,n,\tilde l,m)}(t)=
u_{(\lambda =1 ,n,\tilde l,m)}(\lambda t).
\end{equation}
Then the
eigenvalue regarding the corresponding eigenfunction of 
(\ref{rad}) is
\begin{equation}
\mu_{(\lambda,n,\tilde l,m)}=-((4n+4\tilde p+k)\lambda+4k\lambda^2).
\end{equation}

These actual eigenvalues depend on $p=n+\tilde p$. 
If this parameter is fixed, 
the multiplicity of the corresponding eigenvalue is
\begin{equation}
\sum_{n=0}^p{{p-n+k-1}\choose{k-1}}{{\tilde\upsilon +k-1}\choose{k-1}},
\end{equation}
Finally, if there are distinct eigenvalues 
$\{\lambda_1,\dots ,\lambda_s\}$ involved, the eigenfunctions 
are the products of functions of the form
\begin{equation}
u_{(\lambda =1,\tilde l,m ,n)}(\lambda_i \langle X,X\rangle_i)
e^{-{1\over 2}\lambda_i \langle X,X\rangle_i}
\Phi_{(\tilde{l}_i,m_i)}(X_i),
\end{equation}
where $X_i$ is the orthogonal projection of a vector $X$ onto the
maximal eigensubspace belonging to $\lambda_i$ and 
$\langle X,X\rangle_i=\langle X_i,X_i\rangle$.
Then the eigenvalue corresponding to this function is
\begin{equation}
\mu_{(\lambda,n,\tilde l,m)}=-\sum_{i=1}^s((4n_i+4\tilde p_i+k_i)
\lambda_i+4k_i\lambda_i^2),
\end{equation}
where $n=\sum n_i\, ,\, l=\sum \tilde{l}_i\, ,\, m =\sum m_i$,
and $\tilde{\upsilon}_i=\tilde{l}_i-\tilde{p}_i$.

The multiplicity of the eigenvalue corresponding to fixed values 
of parameters $p_i=n_i+\tilde{p}_i$ are
\begin{equation}
\prod_{i=1}^s\sum_{n_i=0}^{p_i}{{p_i-n_i+k_i-1}\choose{k_i-1}}
{{\tilde{\upsilon}_i+k_i-1}\choose{k_i-1}}.
\end{equation}
\medskip

\noindent {\bf Zeeman effect.}
According to the hypothesis, 
an electron of charge $-e$
and mass $M$ revolving about the
origin (nucleus) has {\it magnetic dipole moment}
$\mathcal S=-(eh/2Mc)\mathcal L$ associated with the {\it angular
moment}
$J_x=hL_x=yp_z-zp_y,\dots$ etc.. The influence
of $\mathcal L$ must be clearly felt under magnetic influence. This
influence was measured by the Stern-Gerlach experiment such
that a stream of one electron atoms was moved in the x-direction
through a constant magnetic field 
$\mathbf K=B\partial_z$ pointing into the
z-direction. In performing the experiment on silver atoms in normal
state the beam was broken up into 2 beams (Zeeman effect) proving
the reality of the magnetic dipole moment.

In the mathematical model of this paper the angular momentum, 
$(J_x,J_y,\\ J_z)$, 
corresponds to the endomorphisms $J_\gamma$ and the quantum 
angular momentum operator is
$2\mathbf D_\gamma\bullet$. 
The Zeeman effect can be interpreted as the influence
this operator exerts on the spectrum of the associated operators such
as the harmonic oscillator operator $\mathbf O$ or the Coulomb 
operator. Since it commutes with both of them, 
the angular momentum operator breaks up
(splits) the spectral lines of both associated operators, defining 
finer spectra. This effect is exhibited
in the classical splitting
\begin{equation}
\mathbf H^{(\tilde l)}
e^{-{1\over 2}\lambda|X|^2}
=\oplus_p\mathbf H^{(p,\tilde l-p)}
e^{-{1\over 2}\lambda|X|^2}
\end{equation}
as well as in the Zeeman zone decomposition.

In classical quantum theory the quantum numbers are established
by means of this standard spectral computation.
The quantum
angular momentum operator has the eigenvalues $m=2\tilde p-\tilde l$ on
such a common eigensubspace 
$\mathbf H^{(\tilde l)}$ 
and the range of these eigenvalues 
is $m=-\tilde l,\dots \tilde l$. Values $\tilde l$ and $m$ are called
{\it azimuthal quantum numbers} and {\it magnetic quantum numbers} 
respectively. The quantum physical interpretation of these numbers
is the following: The angular momentum is quantized and in a 
classical azimuthal quantum
state, $\tilde l$, it can take $\tilde l+1$ distinct values 
determined by the
formula $m=2\tilde p-\tilde l$, where $p$ runs through the values
$\tilde p=0,\dots ,\tilde l$. Thus a beam of electrons
being on the same azimuthal quantum state must split up into
$\tilde l+1$ smaller beams.
Comparing the zonal and classical quantum numbers note that the
corresponding azimuthal numbers are different while the magnetic 
numbers are the same.

This was the very first explanation given for the Zeeman effect.
In this interpretation only orbital angular momentum is attributed 
to the charged particle. Because
of the discrepancies between the computed and measured quantities, 
Pauli added a non-relativistic spin angular momenta to the orbital one.
Dirac developed the relativistic version of this operator. 
These operators are considered in separate papers \cite{sz5,sz6}.
The corresponding zones are called
{\it anomalous zones}. 

Due to the irreducibility, in quantum physics mostly the 
holomorphic zone, described by 
$\upsilon =0$, is involved to theoretical explanations.
One of the new features in this
paper is that the quantization
developed here applies reducible Heisenberg group 
representation (extended Fock representation). 
This reducibility is inevitable since 
the quantization of the angular momenta
can not be solved on a single zone.  
Indeed, for a fixed zonal azimuthal quantum number $l$, 
every zone represents only a single quantum number $m=2p-l$.
In other words, without using the reducible representation the
quantization of the angular momentum could have not been established.
 
In the original approach, where the Zeeman operator is defined on
$\R^3$, the 3-space is endowed by the polar (spherical) coordinate
system $(x,y,r)$ and the z-coordinate is eliminated by the
equation $x^2+y^2+z^2=r^2$. Then the eigenfunctions of the 
Hamilton operator
are sought again in the form $F(r)h^{(\tilde l)}(x,y)$, where the 
$h^{(\tilde l)}$ is an
$\tilde l^{th}$-order spherical harmonics whose degree defines the
azimuthal quantum number. According to H. Weyl \cite{we}, pages 66,  
the silver atom can be in three normal states corresponding to
$m=-1,0,+1$ associated with the functions 
$\overline z,z\overline z,z$.
Then Weyl continues:`` On performing the experiment on silver atoms
in the normal state two beams, corresponding to $m=+1,-1$, were 
observed. Why the unperturbed beam corresponding to $m=0$ did not appear
remained unexplained."

In establishing the Zeeman zones the eigenfunctions are sought not in
the above form, due to that that the zones are not invariant 
with respect to multiplications by radial functions.  
A possible explanation to the above problem may be that that 
the functions
$\overline z$ and $z\overline z$ are in the same zone, 
$\mathcal H^{(1)}$, thus one can talk about just two normal state 
in zonal theory.
\medskip

\noindent{\bf Standard method used for isospectrality constructions.}
So far, the explicit spectrum has been computed on the non-compact
manifolds $\Gamma_Z\backslash N$. If one considers compact 
sub-torus-bundles with boundaries (for instance, ball$\times$torus-type
domains which have sphere$\times$torus-type boundaries), the technique
described in the previous section can not be straightforwardly adopted
for computing the spectrum. One of the advantages of
determining the eigenfunctions in the standard way 
is that this technique can be applied also for these compact domains.
This method was used in \cite{sz3,sz4} for computing the
spectrum both on ball$\times$torus- and sphere$\times$torus-type domains
of Heisenberg-type Lie groups, in which case only a single
eigenvalue $\lambda$ is involved to the computations. Then, 
on ball$\times$torus-type domains,
one has to determine both the Dirichlet and Neumann spectra 
of operator
$
\Lambda_{({1\over 2}k+\tilde l-1)}
$
defined by formula (\ref{Lop}) on the interval $[0,R]$, 
where $R$ is the radius of the ball. 
Comparing with the non-compact case, there arise
numerous differences. First of all, the eigenfunctions are not 
polynomials and they are not known explicitly.
One can establish only
implicit formulas, and as a result, one can not entirely 
explicitly compute the 
eigenvalues $\rho_{(k,\tilde l,n)}\to -\infty$ of the Laguerre operator
$
\Lambda_{({1\over 2}k+\tilde l-1)}
$
with the given boundary conditions on a compact interval. 
The eigenvalues of the ball$\times$torus-type domain depend on this
implicitly determined spectrum by the semi-explicit formula
\begin{equation}
\mu_{(\lambda,\tilde l,m ,n)}=-((-4\rho_{(k,\tilde l,n)}+4p+k)\lambda+4k
\lambda^2).
\end{equation}
For fixed $n$ and $\tilde l$, 
the multiplicity of such an eigenvalue is the 
dimension of the space 
of the corresponding spherical harmonics $H^{(\tilde l,m)}$.
These computations are not used in this paper.

\begin{center} 
\large{\bf PART TWO}
\end{center} 
\begin{center} 
 \large{\bf NORMAL DE BROGLIE GEOMETRY}
\end{center} 
\medskip

\noindent{\bf Technical introduction.}
In the next chapters the zonal versions of the global heat,
$e^{-tH}$, 
and Feynman-Dirac kernel, 
$e^{-\mathbf itH}$, are established,
where $H$ stands for $H_Z$ or   
$H_{Zf}=-{1\over 2}\Box_\gamma$. 
Remember that these two operators differ from each other 
just by the constant term $2|Z_\gamma|^2$. 
The global Wiener-Kac (WK) and Feynman-Dirac (FD) kernels defined by 
$H_Z$ are denoted in the unified form
\begin{equation}
d_{\sigma}(t,X,Y)=e^{-\sigma tH_Z}(X,Y),
\end{equation}
where $\sigma =1$ resp. $\sigma =\mathbf i$ indicate the WK- resp.
FD-kernel. For $H_{Zf}$ these kernels will be denoted differently.
The heat kernel is the 
fundamental solution
of the heat equation by which the heat distribution $u_t(X)$
for a given initial function $u_0(X)$ can be constructed by the 
convolution
$u_t(X)=\int d_1(t,X,Y)u_0(Y)dY$. Using FD-kernels, $d_{\mathbf i}$, 
in this formula, one has the solutions of the Schr\"odinger 
equation. (The computations are carried
out only for these two values of $\sigma$, however, the formulas
established later easily extend to kernels defined by arbitrary 
unit complex number $\sigma$.)   

For a non-degenerated $J_\gamma$ the spectrum, 
$\{\mu_i\geq 0\}$, is discrete, thus 
these kernels can be introduced by the
infinite series
\begin{equation}
\label{d_sigm}
d_{\sigma} (t,X,Y)=
\sum e^{-t\sigma\mu_i}\psi_i(X)\overline{\psi}_i(Y),
\end{equation}
where the elements of the orthonormal basis $\{\psi_i\}$ 
on the $L^2$-space are eigenfunctions of operator 
$H_Z$. These kernels satisfy the limit-property
\begin{equation}
\label{lim}
\lim_{t\to +0}d_\sigma (t,X,Y)=\delta_X(Y).
\end{equation}
The kernels defined for $H_{Zf}$
are denoted by
$p_\sigma (t,X,Y)$. Because of the simple connection between 
the two kind of kernels,
it is enough to consider only the classical kernels $d_\sigma$.

The existence of the global kernels seems to be jeopardized by
the infinite multiplicities of the eigenvalues. Yet, surprisingly
enough, (\ref{d_sigm}) defines smooth kernels in either cases. They are
explicitly computed in the next section. The existence of the
heat kernel $p_1(t,X,Y)$ follows from that fact that operator $\Box$
is derived from the Laplacian of a Riemannian manifold. By this reason
there exist also well defined Wiener-Kac measures, $w_{xy}$, 
on the path-spaces 
$\mathcal P_{xy}$, consisting
of continuous curves connecting two arbitrary points $x$ and $y$.
The heat kernel is not of the trace class, however. 
Let it be mentioned that by the technique of regularization one
can define relative heat kernels of the trace class which define
relative zeta and eta functions.

More serious problems arise with respect to the 
Dirac-Feynman kernel
$d_{\mathbf i}(t,X,Y)$. By (\ref{d_i}) this kernel is of the form 
$A(t)e^{B(t,X,Y)\mathbf i}$. 
Thus, for fixed $t$ and $X$, the integral   
$\int d_{\mathbf i}(t,X,Y)dY$ 
does not exist (the integrand is not absolute integrable, anyway). 
Also the approximating measures $w^{(n)}_{\mathbf ixy}$,
defined analogously to the approximating measures of the WK-measure, 
are divergent. It is well known from the history of this problem that
Kac was able to define this measure only by using the heat kernel.
(Regarding the Euclidean Laplacian $\Delta_X$, 
this measure had been constructed 
by Wiener, earlier.) Later, Feynman and Kac established also the 
Radon-Nikodym derivative of the Wiener-Kac measure with respect to the
Wiener measure. The corresponding formula became known as the 
Feynman-Kac formula. 

Thus it is surprising 
that both zonal flows, introduced by restricting 
the global flows onto the zones, define trace class zonal kernels.
Additionally, both define the corresponding zonal measures, 
$w^{(a)}_\sigma$, on the path-spaces $\mathcal P_{xy}$ rigorously. 

\section{Global WK- resp. DF-kernels}

The explicit computation of the global heat kernel is greatly 
facilitated
by the explicit kernel 
$D=e^{-{1\over 2}t(-\mathbf i\nabla-\mathbf a)^2}$
established, in \cite{ahs,si}, 
on $\R^3$,
where the vector potential $a_1=-(1/2)Bx_2, a_2={1\over 2}Bx_1, a_3=0$ 
corresponds to the
constant magnetic field $(0,0,B)$. The other physical quantities such
as $\hbar, e,c, M$ are chosen to be unit numbers. 
According to this result, the 
kernel in the question is as follows: 
\begin{eqnarray}
D(t,x_1,x_2,x_3)=
{B\over{(2\pi t)^{1/2}4\pi sinh({1\over 2}Bt)}}exp\big(-{1\over{2t}}
(x_3-y_3)^2\\
-{B\over 4}coth({1\over 2}Bt)
[(x_2-y_2)^2+(x_1-y_1)^2]-{1\over 2}\mathbf iB
(x_1y_2-x_2y_1)\big).
\nonumber
\end{eqnarray}
The computations are carried out by the Feynman-Kac-Ito formula
regarding the well defined Wiener measures (Brownian motion technique).
Since the Hamiltonian in this case 
is the sum of $H_Z$, defined on the 
$(x_1,x_2)$-plane,
and $-{1\over 2}\partial^2_{x_3}$,
the above heat flow is the product of the 
1-dimensional Euclidean heat flow and
\begin{eqnarray}
d_{1\gamma} (t,X,Y)=e^{-tH_0}(t,X,Y)=\\
{B/2\over 2\pi sinh({1\over 2}Bt)}e^{
-{B\over 4}coth({1\over 2}Bt)
|X-Y|^2-{1\over 2}\mathbf iB
\langle X,J(Y)\rangle}\nonumber \\
={\lambda\over 2\pi sinh(\lambda t)}e^{
-{\lambda\over 2}coth(\lambda t)
|X-Y|^2-\mathbf i
\langle X,J_\gamma (Y)\rangle},
\nonumber
\end{eqnarray}
where
$\lambda =(1/2)B$. (Note that $J_\gamma$ involves $\lambda$, since
the eigenvalues of $-J^2_\gamma$ are $\lambda^2$.) 

One has a self-contained proof for this formula, however, by 
showing that this kernel satisfies the heat equation
$
(\partial_t+H_{ZX})d_{1\gamma}(t,X,Y)=0
$ 
as well as the limit
property (\ref{lim}). The heat equation can be readily established 
by the following computations, carried out for $k=2$.
\begin{eqnarray}
\label{heateq}
\partial_t(d_{1\gamma})(t,X,Y)
=\big(-\lambda coth(\lambda t)+
{\lambda^2|X-Y|^2\over 2 sinh^2(\lambda t)}\big)
d_{1\gamma}(t,X,Y),\\
H_{ZX}(d_{1\gamma})(t,X,Y)=
\big(\lambda coth(\lambda t)+{\lambda^2\over 2}\{-coth^2(\lambda t)
|X-Y|^2 \\
+|Y|^2-2\mathbf icoth(\lambda t)\langle X,J(Y)\rangle
-2\mathbf icoth(\lambda t)\langle J(X),Y\rangle\nonumber \\
-2\langle X,Y\rangle +|X|^2\}\big)d_{1\gamma}(t,X,Y)
=-\partial_t(d_{1\gamma})(t,X,Y).\nonumber
\end{eqnarray}
For a general $k$ the kernel considered is the product of kernels 
determined for the coordinate planes $z_i$. Thus the 
heat equation
is established in generality. Limit property (\ref{lim}) will be shown
after introducing the zonal flows 
$
d_{1}^{(a)}(t,X,Y)
$ 
and, then, by proving that
\begin{equation}
\label{lim^a}
\lim_{t\to +0}
d_{1\gamma}^{(a)}(t,X,Y)
=\delta_\gamma^{(a)}(X,Y),
\end{equation}
valid for general $k$-dimensions.
 
The eigenvalues of the classical Zeeman-Hamilton operator and the one
introduced in this paper differ from each other by
the extra $2|Z_\gamma|^2=2\lambda^2$ constant. 
Thus, for $k=2$, the sought heat kernel is
\begin{equation}
\label{p}
p_{1\gamma}(t,X,Y)=e^{-2t|Z_\gamma |^2}d_{1\gamma}(t,X,Y).
\end{equation}
The heat kernel on the torus bundle on $\Gamma\backslash N$
is the direct sum of kernels 
$
b_{1\gamma}(t,X,Y,Z_x,Z_y)
$ 
determined on the invariant subspaces $FW^\gamma$. One obviously has
\begin{equation}
\label{b}
b_{1\gamma}(t,X,Y,Z_x,Z_y)=
p_{1\gamma}(t,X,Y)
e^{2\mathbf i\langle Z_\gamma ,Z_x-Z_y\rangle}.
\end{equation}

Because of these straightforward connections, the computations are 
carried out just for the kernel $d_1(t,X,Y)$.
For determining $d_{1\gamma}(t,X,Y)$,
defined by a general endomorphism $J_\gamma$, one should 
proceed as follows. 

On the maximal subspace where the $J_\gamma $ is 
degenerated consider the standard Euclidean heat kernel. Then consider
the maximal invariant subspace on which the 
$J_\gamma$ is non-degenerated. 
Decompose it into orthogonal direct sum of 2-dimensional 
subspaces which are invariant under the action of $J_\gamma$. 
Such a decomposition can
be established by considering a complex coordinate neighborhood
$(z_1,\dots ,z_p)$ of the complex structure defined by normalizing
$J_\gamma$ on the non-degenerated subspace. If $\lambda_i\mathbf i$ is
the eigenvalue on the complex plane corresponding to $z_i$ then consider
the kernel $d_{1\lambda_i}(t,X,Y)$ considered in (\ref{heateq}).
Then the desired kernel, $d_{1\gamma}(t,X,Y)$, is the product of these 
kernels. Thus the global Wiener-Kac kernel is established in generality.
   
For two fixed points $x$ and $y$, this kernel defines 
the well known complex WK-measure, $w_{xy}$, on the 
complete separable metrizable space 
$\mathcal P_{ab}(M)$ of continuous maps
$I\to M$ carrying $(0,1)\to (x,y)$. This space is topologized 
with the topology of uniform convergence \cite{ee}. 
It is elementary that the Borel $\sigma$-algebra
(generated by the open sets) is also generated by the fibred sets
$\rho_{\mathbf t}^{-1}(B)\subset\mathcal P_{xy}(M)$, where 
$\mathbf t=(0<t_1<\dots <t_n<1$ and $\rho_{\mathbf t}:
\mathcal P_{xy}(M)\to
M^n=M\times\dots\times M$ is the evaluation map 
$\rho_{\mathbf t}(x)=(x(t_1),\dots ,x(t_n))$ and $B$ is a Borel subset
of $M^n$. Measure $w_{xy}$ on a fibred set 
$\rho_{\mathbf t}^{-1}(B)$ is defined by
\begin{eqnarray}
\label{w_1xy}
w_{\mathbf 1xy}\big(\rho_{\mathbf t}^{-1}(B)\big)=
\int_Bd_{1\gamma}(t_1,x,m_1)\cdot \\
\cdot d_{1\gamma}(t_2-t_1,
m_1,m_2)\dots d_{1\gamma}(1-t_n,m_n,y)dm_1\dots dm_n.\nonumber
\end{eqnarray}
Above the manifold is either compact, or, one considers the problem
on the manifold defined by the one-point compactification 
$M\cup \infty$.
Thus the $\mathcal P_{xy}$ is a compact topological space.
By classical results; such as Riesz' theorem concerning the measure
representation of bounded linear functionals on the Banach space of
continuous functions defined on a compact metrizable space (the Banach 
norm is defined by $sup |f|$) and the Stone-Weierstrass theorem 
asserting that the curves 
$\rho_{\mathbf t_n}^{-1}(x_1,\dots ,x_n)$ corresponding
to ${1\over n}<\dots <{n-1\over n}$ are dense in $\mathcal P_{xy}$; 
this construction determines a complex countably
additive regular Borel measure 
$w_{xy}$ on $\mathcal P_{xy}(M)$ satisfying
\begin{equation}
w_{\mathbf 1xy}\big(\mathcal P_{xy}(M)\big)=d_{1\gamma}(1,x,y).
\end{equation}

Next we proceed with the
FD-kernel, defined by $\sigma =\mathbf i$ 
in (\ref{d_sigm}), in a similar fashion. 
In-spite of the infinite
multiplicities of the eigenvalues, also this kernel exists which can 
be established by means of the FD-kernel 
\begin{eqnarray}
{-B\mathbf i/2\over 
(2\pi\mathbf it)^{1/2}
2\pi sin({B\over 2}t)}exp\big(
{\mathbf i\over{2}}\{
{(x_3-y_3)^2\over t}\\
+{B\over 2}(cot{Bt\over 2})[(x_2-y_2)^2+(x_1-y_1)^2]+B
(x_1y_2-x_2y_1)\}\big),\nonumber
\end{eqnarray}
established in \cite{fh}, formula (3-64), for a 
rotating charged particle in a constant magnetic field.
It is the product of the one-dimensional Euclidean Feynman-Dirac kernel
\begin{equation}
\label{K_0}
K_0(t,x,y)=
(2\pi\mathbf it)^{-1/2}
exp\big(\mathbf i{|x-y|^2\over{2t}}\big),
\end{equation}
and of the sought Zeeman-Feynman-Dirac kernel 
\begin{eqnarray}
\label{d_i}
d_{\mathbf i}(t,X,Y)=e^{-t\mathbf iH_Z}(t,X,Y)=\\
\big({-B\mathbf i/2\over 2\pi sin({1\over 2}Bt)}\big)^{k/2}e^
{{1\over{2}}\mathbf i\{
{B\over 2}(cot{Bt\over 2})|X-Y|^2-B\langle X,J(Y)\rangle\}} \nonumber \\
=\big({\lambda\over 2\pi\mathbf i sin(\lambda t)}\big)^{k/2}e^
{\mathbf i\{
{\lambda\over 2} cot(\lambda t)|X-Y|^2-\langle X,J_\gamma (Y)\rangle\}}.
\nonumber
\end{eqnarray}
Here we suppose that
$-J^2_\gamma$ has only one eigenvalue, $\lambda^2$. 
In the general cases this kernel is an appropriate product of 
kernels of the form (\ref{K_0}) and (\ref{d_i}).

Also this formula can be established in a self-contained manner by 
showing that both (\ref{lim}) and the Schr\"odinger equation
$
(\partial_t+\mathbf iH_{ZX})d_{\mathbf i\gamma}(t,X,Y)=0
$ 
hold. The limit property (\ref{lim}) 
will be proved later by pointing out that the zonal kernels 
approach, by the limit $\lim_{t\to +0}$, to the zonal 
Dirac $\delta$-spreads. The Schr\"odinger equation can 
be established
by the following computations,
what should be carried out also in this case only for $k=2$.
\begin{eqnarray}
\partial_t(d_{\mathbf i\gamma})(t,X,Y)
=\big(-\lambda cot(\lambda t)-\mathbf i
{\lambda^2|X-Y|^2\over 2 sin^2(\lambda t)}\big)
d_{\mathbf i\gamma}(t,X,Y),\\
\mathbf i(H_{Z})_{X}(d_{\mathbf i\gamma})(t,X,Y)=
\big(\lambda cot(\lambda t)+
{\lambda^2\mathbf i\over 2}\{cot^2(\lambda t)
|X-Y|^2 \\
+|Y|^2-2cot(\lambda t)\langle X,J(Y)\rangle
-2cot(\lambda t)\langle J(X),Y\rangle\nonumber \\
-2\langle X,Y\rangle +|X|^2\}\big)d_{\mathbf i\gamma}(t,X,Y)=-
\partial_t(d_{\mathbf i\gamma})(t,X,Y).\nonumber
\end{eqnarray}

There arise serious difficulties when one tries to 
construct the complex Feynman measure on the path-space by 
the same steps applied for introducing the 
Wiener-Kac measure. These problems due to the fact that,
for fixed $t$ and $x$, function $d_{\mathbf i}(t,x,y)$ 
(depending on $y$) is neither of the class $L^1$ nor $L^2$ and 
integral (\ref{w_1xy}) does
not exist in general. (Yet, the formally defined Feynman integral is
a very useful heuristic tool in perturbation theory.
In particular, perturbation expansions can be computed explicitly.) 
There have been also various attempts made for defining an appropriate 
Feynman measure. These problems are out of the scope of this paper.
Below a non-perturbative theory is offered which is the most important
new feature in this article.

For a fixed $t$ the function $d_\sigma(t,X,X)$ is constant, thus
the corresponding convolution operators are not of the trace class.
Therefore, important functions such as the partition function, 
zeta function, eta function, and the determinant of the Hamilton
operator can not be defined in the usual way. Actually, these
non-existing integrals cause the infinities appearing in
quantum field theory. In perturbation theory one gets rid off these
infinities by ``subtracting two infinities" in order to get the
desired finite quantity. A typical example is that one adds a
suitable potential function $V$ to a non-trace class Hamiltonian $H_Z$
such that the kernel $e^{-tH_Z}-e^{-t(H_Z+V)}$ is of the trace class.
Then, one defines the above functions (called relative functions)
with respect to this relative kernel. 
In physics the usual designation for this process is  
{\it renormalization}. The most important papers concerning this topic
are collected in \cite{schw}. 

To implement spectral 
investigations on non-compact Riemannian manifolds, this tool gained
ground also in mathematics. It is
called {\it regularization} which tool includes, besides perturbations,
also depicting and removing the divergent terms from functions. 
This type of investigations started 
out with \cite{ops}, which article
was apparently motivated by papers \cite{a,po}
written in physics. One can consult with \cite{mu} for 
finding more informations about the recent developments in this field. 
So-far this is the major tool for controlling
the infinities on both areas. In this article a new non-perturbative
approach is offered to these infinities. Namely, the kernels
$d_\sigma (t,X,Y)$ 
will be restricted to the zones and then shown that
both zonal kernels
$d_\sigma^{(a)}(t,X,Y)$
are trace class kernels. A much more striking discovery is that
also the zonal Feynman-Dirac kernels, 
$d_{\mathbf i}^{(a)} (t,X,Y)$, determine a well defined 
zonal Feynman measure, $w^{(a)}_{\mathbf ixy}$,
on the path-space.

The most general form of the statements discussed in this section
is as follows.
\begin{theorem}{\bf (Global Kernel Theorem)}
Infinite function series (\ref{d_sigm}) define smooth kernels both for
the Wiener-Kac setting, $\sigma =1$, and the Dirac-Feynman setting,
$\sigma =\mathbf i$, where $H$ is the classical Zeeman-Hamilton
operator $H_Z$, or, it is the Zeeman-Hamilton Laplacian 
$H_{Zf}$. The kernels corresponding to these operators 
are distinguished by $d_\sigma$ and $p_\sigma$ respectively. These
kernels are defined on the X-space $\R^k$. If the
kernel is defined on the torus bundle $\Gamma\backslash N$, it is 
denoted by $b_\sigma$. The simple connections among these kernels
are described in (\ref{p}) and (\ref{b}), thus, it is enough to describe
only the classical case.

Suppose that the Hamilton operator is non-degenerated such that the
distinct non-zero parameters $\{\lambda_i\}$,  $i=1,\dots ,r$,
are defined on $k_i$-dimensional subspaces. Then for the Wiener-Kac 
kernel we have
\begin{eqnarray}
\label{d_1gamm}
d_{1\gamma} (t,X,Y)=e^{-tH_Z}(t,X,Y)=\\
=\prod \big({\lambda_i\over 2\pi sinh(\lambda_i t)}\big)^{k_i/2}e^{
-\sum{\lambda_i}({1\over 2}coth(\lambda_i t)
|X_i-Y_i|^2+\mathbf i
\langle X_i,J (Y_i)\rangle}.\nonumber
\end{eqnarray}
This kernel satisfies the Chapman-Kolmogorov identity as well as the
limit property (\ref{lim}), however, it is not of the trace class.

The explicit form of the Feynman-Dirac kernel is
\begin{eqnarray}
d_{\mathbf i}(t,X,Y)=e^{-t\mathbf iH_Z}(t,X,Y)=\\
=\prod\big({\lambda_i\over 2\pi
\mathbf i sin(\lambda_i t)}\big)^{k_i/2}e^
{\mathbf i\sum{\lambda_i}\{
{1\over 2} cot(\lambda_i t)|X_i-Y_i|^2-\langle X_i,J(Y_i)\rangle \}}.
\nonumber
\end{eqnarray}
Since for fixed $t$ and $X$ function $d_{\mathbf i}(t,X,Y)$ is neither
an $L^1$- nor an $L^2$-function of the variable $Y$, 
the integral required
for the Chapman-Kolmogo\-rov identity is not defined. It is not of the 
trace class either. Nevertheless, it satisfies the 
limit property (\ref{lim}).
\end{theorem}

\section{Zonal WK-kernels}

{\bf Establishing the zonal WK-kernels.}
Since the zones are spanned by eigenfunctions, the corresponding 
{\it zonal Wiener-Kac- and zonal Feynman-Dirac-kernels} can be defined
by using only the eigenfunctions belonging to the zone in  
(\ref{d_sigm}).
One can use also the projections $P^{(a)}$
to define these kernels by the integral formula
\begin{eqnarray}
d_\sigma^{(a)}(t,X,Y)=
\int\int P^{(a)}(X,U)P^{(a)}(V,Y)
d_\sigma (t,U,V)dUdV\\
=\int P^{(a)}(X,U)
d_\sigma (t,U,Y)dU=
\int P^{(a)}(V,Y)
d_\sigma (t,X,V)dV.\nonumber
\end{eqnarray}
The last two equations follow from the invariance of
the zones under the actions of both flows.

The basic tool used in computations below is the well
known integral formula:
\begin{equation}
\label{intform}
\int_{\R^N}exp(-{1\over 2}Z\cdot A\cdot Z +C\cdot Z)dZ=
\big({(2\pi)^N\over det[A]}\big)^{1/2}
exp({1\over 2}C\cdot A^{-1}\cdot C),
\end{equation}
where $A$ is a complex diagonal matrix having the same entries $a$,
satisfying $Re(a)>0$, in the diagonal and $C$ is a complex
$N$-vector.  

First the zonal kernel
$p_1^{(0)}(t,X,Y)$ is established.
For the sake of simplicity we assume that 
$\lambda =1$. (In the end these formulas will be stated also in 
the most general form.) By (\ref{d_1gamm}) and (\ref{intform}), we have:
\begin{eqnarray}
\label{INT}
\int P^{(0)}(X,Z)d_{1}(t,Z,Y)dZ=\\
({1\over 2\pi^2 sinh(t)})^{k/2}
e^{-{1\over 2}(|X|^2+coth(t)|Y|^2)}INT,\nonumber \\
INT=\int e^{-{1\over 2}Z\cdot A\cdot Z +C\cdot Z}dZ=
\big({(2\pi)^{k}\over det[A]}\big)^{1/2}
e^{({1\over 2}C\cdot A^{-1}\cdot C)},
\end{eqnarray}
where the complex diagonal matrix $A$ has the constant entries
\begin{equation}
a=1+coth(t)
\end{equation}
on the main diagonal and the complex vector $C$ is:
\begin{equation}
C=X-\mathbf iJ(X)
+coth(t)Y-\mathbf iJ(Y),
\end{equation}
where $J$ is the complex structure defined by normalizing $J_\gamma$.
Therefore, by the identities
\begin{eqnarray}
a^{-1}=
e^{-t}sinh(t)
={1\over 2}(1-e^{-2t}),\\
{1\over 2}C\cdot A^{-1}\cdot C ={1\over 2}
e^{-t}sinh(t)
\big((coth^2(t)-1)|Y|^2+ \\
2(coth(t)-1)(\langle X,Y\rangle 
+\mathbf i\langle X,J(Y)\rangle\big)\nonumber
\end{eqnarray}
we have
\begin{eqnarray}
\label{d_1^0}
d_{1}^{(0)}(t,X,Y)=
\int P^{(0)}(X,Z)d_{1}(t,Z,Y)dZ=\\
={e^{-{kt\over 2}}\over \pi^{k\over 2}}
e^{-{1\over 2}(|X|^2+|Y|^2)+e^{-2t}
\langle X,Y+\mathbf iJ(Y)\rangle}.
\end{eqnarray}

Note that this kernel satisfies the limit property (\ref{lim^a}).
Unlike the global one, this kernel is of the trace class 
with trace:  
\begin{eqnarray}
\label{d^1_0}
Trd_{1}^{(0)}(t)
=\int d_{1}^{(0)}(t,X,X)dX\\
=\int
{e^{-{k\over 2}t}
\over \pi^{k\over 2}}
e^{-(1-e^{-2t})|X|^2}dX
={e^{-{k\over 2}t}\over
(1-e^{-2t})^{k\over 2}}\nonumber
\end{eqnarray}

The computations regarding the next kernel,
\begin{eqnarray}
d_1^{(1)}(t,X,Y)
=\int P^{(1)}(X,Z)d_{1}(t,Z,Y)dZ\\
=\int L_1^{((k/2)-1)}(|X-Z|^2)P^{(0)}(X,Z)d_{1}(t,Z,Y)dZ,\nonumber
\end{eqnarray}
can be traced back to the above ones by integrating by parts
in
\[
INT^{(1)}=\int L_1^{((k/2)-1)}(|X-Z|^2)
e^{-{1\over 2}Z\cdot A\cdot Z +C\cdot Z}dZ.
\]
According to this computations we have
\begin{eqnarray}
d_1^{(1)}(t,X,Y)=
\Lambda^{(1)}(t,X,Y)
d_1^{(0)}(t,X,Y)\\
=(L_1^{((k/2)-1)}(|X-Y|^2)+
LT_1^{(1)}(t,X,Y))
d_1^{(0)}(t,X,Y),\nonumber\\
\label{LT^1}
\mbox{where}\quad\quad\quad\quad\quad\quad
LT_1^{(1)}(t,X,Y))=
(1-e^{-2t})
lt_1^{(1)}(t,X,Y))=
\\
(1-e^{-2t})
(|X|^2+|Y|^2-1-
(1+e^{-2t})\langle X,Y+\mathbf iJ(Y)\rangle ).\nonumber
\end{eqnarray}

On a general zone the heat kernel is
\begin{eqnarray}
d_1^{(a)}(t,X,Y)=
(L^{((k/2)-1)}_a(|X-Y|^2)+
LT_1^{(a)}(t,X,Y))
d_1^{(0)}(t,X,Y),
\end{eqnarray} 
thus it is the sum of the {\it dominant term}
\begin{eqnarray}
D_{1}^{(a)}
={e^{-{kt\over 2}}\over \pi^{k\over 2}}
L^{((k/2)-1)}_a(|X-Y|^2)
e^{-{1\over 2}(|X|^2+|Y|^2)+e^{-2t}
\langle X,Y+\mathbf iJ(Y)\rangle}
\end{eqnarray}
and the {\it long term zonal heat kernel}
\begin{eqnarray}
\tau_{1}^{(a)}=
(1-e^{-2t})T^{(a)}_1(e^{-2t},|X|^2,|Y|^2
,\langle X,Y+\mathbf iJ(Y)\rangle )
d_1^{(0)}(t,X,Y),
\end{eqnarray}
where the $T_1^{(a)}$ is a $2a^{th}$-order polynomial regarding 
the $X,Y$ variables and a $(2a-1)^{th}$-order polynomial of the 
$e^{-2t}$ variable. This kernel vanishes for $t=0$, 
explaining its name. The method of integrating by parts
and recursions (\ref{rec1})-(\ref{rec3}) yield
a recursion formula for the long term kernel $T_1^{(a)}$.

Since a gross zone decomposes into
$
{a+(k/2)-1\choose a}
={a+(k/2)-1\choose (k/2)-1}
$
number of irreducible zones, each of them is isospectral to the
irreducible holomorphic zone, the partition function on a gross zone
is
\begin{eqnarray}
\mathcal Z^{(a)}(t)=
Trd_{1}^{(a)}(t)
={a+(k/2)-1\choose a}
{e^{-{k\over 2}t}\over
(1-e^{-2t})^{k\over 2}}.
\end{eqnarray}
On the other hand, by (\ref{rec0}),
$
L^{((k/2)-1)}_a(0)
={a+(k/2)-1\choose a}
={a+(k/2)-1\choose (k/2)-1}
$.
Thus the partition function is determined by the dominant heat kernel
and the long term heat kernel is in the 0 trace class.
\medskip

\noindent{\bf Zonal partition and zeta functions.}
These formulas allow
to introduce the zonal zeta functions by the following 
standard way:
\begin{equation}
\label{zeta^a}
\zeta^{(a)} (s)={1\over \Gamma (s)}\int^\infty_0t^{s-1}
\mathcal Z^{(a)}(t)dt=\sum {1\over\mu_i^s},
\end{equation}
where $\mu_i>0$ are the eigenvalues of the Hamilton operator which
can be either $H_Z$ or 
$-{1\over 2}\Box_\gamma =H_Z+2|Z_\gamma|^2:=H_{Zf}$. The infinite
series (\ref{zeta^a}) is absolute convergent on the half-plane
$Re(s)>1$ such that the 
function defined has meromorphic extension onto the whole complex plane.

Interestingly enough, these zeta functions strongly relate
to the Riemann,
$\zeta_{R}(s)=\sum {1\over n^s}=\Gamma^{-1}(s)\int_0^\infty t^{s-1}
(e^t-1)^{-1}dt$, 
and the Hurwitz zeta function,  
$\zeta_{Hu}(s,x)=\sum {1\over (n+x)^s}$, respectively. 
(Though Hurwitz defined
these functions only for $0\leq x\leq 1$, we allow arbitrary positive
numbers. If $x$ is a natural number, one gets the Hurwitz function 
by removing the first $x$ terms from the Riemann zeta function.)
 
In fact, in the original electron case, defined by 
$k=2,\lambda =1$, the multiplicity of eigenvalues on each zone is one
and, by (\ref{eigv2}), one has
\begin{eqnarray}
\label{Rzeta}
\zeta^{(a)}_{H_Z}(s)+2^{-s}\zeta_R(s)=\zeta_R(s)\, ,\,
\zeta^{(a)}_{H_{Zf}}(s)+2^{-s}\zeta_{Hu}(s,4)=\zeta_{Hu}(s,4),\\
\zeta_{R}(s)={2^{s}\over 2^s-1}\zeta_{H_Z}(s)\, ,\,
\zeta_{Hu}(s,4)={2^{s}\over 2^s-1}\zeta_{H_{Zf}}(s).
\end{eqnarray}
These equations can be easily established 
also by the first equation of (\ref{zeta^a}). They
provide a quantum physical interpretation for these classical zeta
functions. 

Also other zonal functions such as zonal determinants, zonal eta- resp. 
shift-functions can
be introduced by the standard formulas originally established on compact
Riemannian manifolds. We do not go into these details.
\medskip

\noindent{\bf Comparing with the oscillator kernel} 
$e^{-{1\over 2}t \mathbf O_\gamma}$.
Finally, the zonal Wiener-Kac kernels are  
compared with the global kernel $k_{\gamma}(t,X,Y)$ defined for the
harmonic oscillator $\mathbf O_\gamma$ by (\ref{d_sigm}). 
By (\ref{eigv1})
this kernel is:
\begin{equation}
\label{k_gamma1}
k_{\gamma} (t,X,Y)=\sum e^{ 
\sum_{a=1}^s-
(2l_a +k_a)
\lambda_{\gamma a}
t}\psi_i(X)\overline{\psi}_i(Y),
\end{equation}
 
This infinite function 
series converges to a kernel which is
explicitly described by the following Mehler formula: 
\begin{equation}
k_{\gamma}(t,X,Y)={exp\{ {B\over sinh(2Bt)}[-{1\over 2}(cosh(2Bt))
(|X|^2+|Y|^2)+\langle X,Y\rangle ]\over (2\pi sinh(2Bt))^{k/2}}
\end{equation}
To be more precise, this formula describes the flow with respect to the
operator $(1/2)(-\Delta +B|X|^2)$. If there are more constants 
($B_1,\dots ,B_s$) involved, the kernel is the product of kernels
belonging to these constants. By this formula one immediately gets
that the global Mehler kernel is of the trace class, for any fixed 
time $t$.

The infinite function series defining the zonal heat flow is a partial
sum of the following infinite function series
\begin{equation}
\label{k_gamm}
k_{\gamma} (2t,X,Y)e^{t\sum (4\upsilon_a +k_a)\lambda_a},
\end{equation}
therefore, there is reestablished that the zonal heat kernels 
are of the trace class on any zone.

Formula (\ref{k_gamm}) provides an opportunity to compare the
zonal partition function $\mathcal Z^{(a)}(t)$
of the heat flow both with the global and zonal partition
function, $\mathcal Z_k^{(a)}(t)$, of the oscillator operator 
$\mathbf O_{\gamma}$. Indeed, one has:
\begin{eqnarray}
0< \mathcal Z^{(a)}(t)=\mathcal Z^{(a)}_k(t)
exp\big(\sum ((4\upsilon_a +k_a)\lambda_a)t\big)\\
< exp\big(\sum ((4\upsilon_a +k_a)\lambda_a)t\big)
/(B(cosh(2Bt)-1))^{k/2},\nonumber
\end{eqnarray}
where the last term is computed by the well known partition function 
of the harmonic oscillator operator. 

By summing up we have
\begin{theorem}{\bf (Zonal WK-Flow Theorem)} 
Suppose that the classical
Zee\-man-Hamilton operator is non-degenerated, 
having the non-zero parameters
$\{\lambda_i\}$. Then the zonal 
Wiener-Kac kernels are of the trace class, which can be
described, along with their partition functions, by the following
explicit formulas. 
\begin{eqnarray}
\label{d_1^a}
d_{1}^{(0)}(t,X,Y)=
\prod\big({\lambda_i e^{-\lambda_i t}\over \pi}\big)^{k_i\over 2}
e^{\sum \lambda_i(-{1\over 2} (|X_i|^2+|Y_i|^2)+ e^{-2\lambda_i t}
\langle X_i,Y_i+\mathbf iJ(Y_i)\rangle )},
\\
d_{1}^{(a)}(t,X,Y)=
(L^{({k\over 2}-1)}_a(\sum\lambda_i |X_i-Y_i|^2)+LT_1^{(a)}
(t,X,Y))d_1^{(0)}(t,X,Y),
\end{eqnarray}
where 
$LT_1^{(1)}$ 
is defined in (\ref{LT^1}) and for the general terms,
$LT_1^{(a)}$,
recursion formula can be established. Furthermore,
\begin{eqnarray}
\mathcal Z_1^{(a)}(t)=
Trd_{1}^{(a)}(t)
={a+(k/2)-1\choose a}\prod {e^{-{k_i\lambda_i t\over 2}}\over
(1-e^{-2\lambda_i t})^{k_i\over 2}}=
TrD_{1}^{(a)}(t),
\end{eqnarray}
where 
$
D_{1}^{(a)}(t,X,Y)=
L^{({k\over 2}-1)}_a(\sum\lambda_i|X_i-Y_i|^2)
d_{1}^{(0)}(t,X,Y)
$
is the dominant zonal kernel. The remaining long term kernel
in the WK-kernel vanishes for $\lim_{t\to 0_+}$ and is of the 0 trace
class.
The zonal WK-kernels satisfy the Chapman-Kolmogorov identity 
along with the limit property
(\ref{lim^a}). 
(In the text the zonal partition functions and the corresponding zonal 
zeta functions are described in (\ref{d_1^0})-(\ref{zeta^a}).)

In (\ref{Rzeta}) these zeta functions are expressed by the Riemann-
resp. Hurwitz zeta functions and vice versa. These formulas provide
a quantum physical interpretation for these classical zeta functions.
 
In (\ref{k_gamma1})-(\ref{d_1^a}) 
the Wiener-Kac heat kernels are compared with the 
heat kernel belonging to the harmonic oscillator operator 
$\mathbf O_\gamma$, which is of the trace class.

According to (\ref{p}) and (\ref{b}), 
both zonal kernels $p_1^{(a)}$ and $b_1^{(a)}$ can be expressed by
$d_1^{(a)}$. 
The traces of these kernels
are the multiple of $Trd_1^{(a)}$ by the function 
$e^{-2t|Z_\gamma |^2}$. 

The kernels
$
d_{1\lambda_1\dots \lambda_{k/2}}^{(a_1\dots a_{k/2})}
$
on the irreducible zones 
$\mathcal H^{(a_1\dots a_{k/2})}\subset\mathcal H^{(a)}$
differ from the above ones just in the 
Laguerre polynomial term what should be exchanged by
$
\prod_{i=1}^{k/2}
L^{(0)}_{a_1}.
$
The partition function on each irreducible zone is the same as on
the holomorphic zone.
\end{theorem}

\section{Zonal DF-kernels}

Integral formula (\ref{intform}) is used also for establishing the 
zonal Feynman-Dirac kernels. On the holomorphic zone, by  assuming 
$\lambda =1$, the computations are as follows.
\begin{eqnarray}
\int P^{(0)}(X,Z)d_{\mathbf i}(t,Z,Y)dZ=\\
({-\mathbf i\over 2\pi^2\sin t})^{k/2}
e^{-{1\over 2}(|X|^2-cot(t)|Y|^2\mathbf i)}INT,\nonumber \\
INT=\int e^{-{1\over 2}Z\cdot A\cdot Z +C\cdot Z}dZ=
\big({(2\pi)^{k}\over det[A]}\big)^{1/2}
e^{({1\over 2}C\cdot A^{-1}\cdot C)},
\end{eqnarray}
where the complex diagonal matrix $A$ has the constant entries
\begin{equation}
a=1-\cot t\mathbf i=-\mathbf i\sin^{-1}te^{t\mathbf i}
\end{equation}
on the main diagonal and the complex vector $C$ is:
\begin{equation}
C=X-\mathbf iJ(X)-\mathbf i
(\cot tY+J(Y)),
\end{equation}
where $J$ is the complex structure defined by normalizing $J_\gamma$.
Therefore, by the identities
\begin{eqnarray}
a^{-1}=\sin^2t
(1+\cot t\mathbf i)=\mathbf i\sin te^{-t\mathbf i},\\
{1\over 2}C\cdot A^{-1}\cdot C ={1\over 2}\big(-|Y|^2+
2\langle X,\cos 2tY+\sin 2tJ(Y)\rangle + \\
+(2\langle X,-\sin{2t}Y+\cos{2t})J(Y)\rangle-\cot t|Y|^2)\mathbf i\big),
\nonumber
\end{eqnarray}
we have
\begin{eqnarray}
d_{\mathbf i}^{(0)}(t,X,Y)=
\int P^{(0)}(X,Z)d_{\mathbf i}(t,Z,Y)dZ=\\
{(\cos t-\mathbf i\sin t)^{k\over 2}\over \pi^{k\over 2}}
e^{-{1\over 2}\big(|X|^2+|Y|^2\big)+(\cos{2t}-\mathbf i\sin{2t})
\langle X,Y+\mathbf iJ(Y)\rangle}\nonumber \\
={e^{{-kt\over 2}\mathbf i}\over \pi^{k\over 2}}
e^{-{1\over 2}(|X|^2+|Y|^2)+e^{-2t\mathbf i}
\langle X,Y+\mathbf iJ(Y)\rangle}.\nonumber
\end{eqnarray}

Also this kernel is of the trace class with trace:
\begin{eqnarray}
Trd_{\mathbf i}^{(0)}(t)=\int d_{\mathbf i}^{(0)}(t,X,X)dX\\
=\int
{(\cos t-sin t\mathbf i)^{k\over 2}\over\pi^{k\over 2}}
e^{(\cos{2t}-1-\mathbf i\sin{2t})|X|^2}dX\nonumber \\
={(\cos t-\sin t\mathbf i)^{k\over 2}\over
(1-\cos{2t}+\mathbf i\sin{2t})^{k\over 2}}.\nonumber
\end{eqnarray}

The method of integrating by parts applies also to the computations of
the higher order zonal DF-kernels. Then one has
\begin{eqnarray}
d^{(a)}_{\mathbf i}(t,X,Y)=
(L_a^{((k/2)-1)}(|X-Y|^2)+LT^{(a)}_{\mathbf i}(t,X,Y))
d^{(0)}_{\mathbf i}(t,X,Y),\\
\label{LT^1_i}
\mbox{where}\quad\quad\quad\quad\quad\quad
LT_{\mathbf i}^{(1)}(t,X,Y))=
(1-e^{-2t\mathbf i})
lt_{\mathbf i}^{(1)}(t,X,Y))=
\\
(1-e^{-2t\mathbf i})
(|X|^2+|Y|^2-1-
(1+e^{-2t\mathbf i})\langle X,Y+\mathbf iJ(Y)\rangle ).\nonumber
\end{eqnarray}
The long term kernel $LT^{(a)}_{\mathbf i}$ on general zones can be 
described recursively.
Also this kernel decomposes into the {\it dominant term}
\begin{eqnarray}
D_{\mathbf i}^{(a)}
={e^{-{kt\mathbf i\over 2}}\over \pi^{k\over 2}}
L^{((k/2)-1)}_a(|X-Y|^2)
e^{-{1\over 2}(|X|^2+|Y|^2)+e^{-2t\mathbf i}
\langle X,Y+\mathbf iJ(Y)\rangle}
\end{eqnarray}
and the {\it long term}
\begin{eqnarray}
\tau_{\mathbf i}^{(a)}=
(1-e^{-2t\mathbf i})T^{(a)}_1(e^{-2t\mathbf i},|X|^2,|Y|^2
,\langle X,Y+\mathbf iJ(Y)\rangle )
d_{\mathbf i}^{(0)}(t,X,Y).
\end{eqnarray}
The same arguments yield that the 
$T_{\mathbf i}^{(a)}$ is a $2a^{th}$-order polynomial regarding 
the $X,Y$ variables and a $(2a-1)^{th}$-order polynomial of the 
$e^{-2t}$ variable such that it vanishes for $t=0$. Furthermore,
the long term kernel is of the $0$ trace class and the partition
function is determined by the dominant term.

By summing up we have
\begin{theorem}{\bf (Zonal Dirac-Feynman Flow Theorem)} 
Suppose that the classical
Zee\-man-Hamilton operator is non-degenerated, 
having the non-zero parameters
$\{\lambda_i\}$. Then the gross zonal 
Dirac-Feynman kernels are of the trace class which, together with their
partition functions, can be described 
by the following explicit formulas.

\begin{eqnarray}
\label{d_i^a}
d_{\mathbf i}^{(0)}(t,X,Y)=
\prod\big({\lambda_i e^{-\lambda_i t\mathbf i}\over \pi}
\big)^{k_i\over 2}
e^{\sum \lambda_i(-{1\over 2} (|X_i|^2+|Y_i|^2)+ 
e^{-2\lambda_i t\mathbf i}
\langle X_i,Y_i+\mathbf iJ(Y_i)\rangle )},
\\
d_{1}^{(a)}(t,X,Y)=
(L^{({k\over 2}-1)}_a(\sum\lambda_i |X_i-Y_i|^2)+LT_{\mathbf i}^{(1)}
(t,X,Y))d_1^{(0)}(t,X,Y),
\end{eqnarray}
where 
$LT_{\mathbf i}^{(1)}$ 
is described in (\ref{LT^1_i}) and a
general long term, 
$LT_{\mathbf i}^{(a)}$,
can be defined recursively. Furthermore, 
\begin{eqnarray}
\mathcal Z_{\mathbf i}^{(a)}(t)=
Trd_{\mathbf i}^{(a)}(t)
={a+(k/2)-1\choose a}\prod {e^{-{k_i\lambda_i t\mathbf i\over 2}}\over
(1-e^{-2\lambda_i t\mathbf i})^{k_i\over 2}}=
TrD_{\mathbf i}^{(a)}(t),
\end{eqnarray}
where 
$
D_{\mathbf i}^{(a)}(t,X,Y)=
L^{({k\over 2}-1)}_a(\sum\lambda_i|X_i-Y_i|^2)
d_{\mathbf i}^{(0)}(t,X,Y)
$
is the dominant kernel. The remaining long term in the zonal
DF-kernel is of the 0 trace class.

The zonal DF-kernels are zonal fundamental solutions 
of the Schr\"o\-dinger equation: 
$(\partial_t+\mathbf i(H_{Z})_{X})d^{(a)}_{\mathbf i\gamma}(t,X,Y)=0$,
satisfying the Chapman-Kolmogorov identity as well as 
the limit property (\ref{lim^a}).
 
According to (\ref{p}) and (\ref{b}), 
the zonal kernels $p_{\mathbf i}^{(a)}$ and $b_{\mathbf i}^{(a)}$ 
can be expressed by $d_{\mathbf i}^{(a)}$. 
The trace of both kernels
is the multiple of $Trd_{\mathbf i}^{(a)}$ by the function
$e^{-2t|Z_\gamma |^2\mathbf i}$. 

The kernels
$
d_{\mathbf i\lambda_1\dots \lambda_{k/2}}^{(a_1\dots a_{k/2})}
$
on the irreducible zones 
$\mathcal H^{(a_1\dots a_{k/2})}\subset\mathcal H^{(a)}$
differ from the above ones just in the 
Laguerre polynomial term what should be exchanged by
$
\prod_{i=1}^{k/2}
L^{(0)}_{a_1}.
$
The partition function on each irreducible zone is the same as on
the holomorphic zone.
\end{theorem}

\section{Zonal path-integrals}

\noindent{\bf Technicalities and constructions on the 0-zone.}
The existence of zonal Wiener-Kac measure on the path-space follows
from the existence of the corresponding global measure.
Since such completely additive complex measure does not exist 
with respect to the global Dirac-Feynman flow, $d_{\mathbf i}(t,X,Y)$, 
it is surprising that the zonal Dirac-Feyn\-man kernels rigorously 
define complex measures on the path-space.
These measures can be established by the standard method 
used for constructing the Wiener measure. 
Though it is common, we describe
this construction in details.

In the first step consider the one-point compactification
$M=\mathbf v\cup\{\infty\}$ of the X-space. The {\it path-space}, on 
which the measure is to be constructed, consists of curves started 
out from an arbitrarily chosen point $x\in\mathbf v$. Then, for any 
fixed value $0<T$, the set $\mathcal P_x^T$ of continuous curves 
$c: [0,T]\to M $ satisfying $c(0)=x$
is a compact metrizable space where the topology is defined by 
uniform convergence. A continuous curve $c(t)$ is uniquely
determined by the rational points corresponding to the parameters
$t=rT\in [0,T]$, where
$r\in \mathbf Q_{[0,1]}$ is an arbitrary rational number
on the interval $[0,1]$. Therefore,
this path-space can be identified with the infinite product space
\begin{equation}
\label{P_x^T}
\mathcal P_x^T=\prod_{t\in\mathbf Q_{[0,1]}T}M_t.
\end{equation}
The Banach space of continuous complex valued functions defined on 
$\mathcal P_x^T$, where the norm is defined by $||F||=sup|F|$, 
is denoted by 
$C(\mathcal P_x^T)$. 
According to the Stone-Weierstrass theorem, this
space is generated by continuous functions depending only on
finite many of the factors in (\ref{P_x^T}). 
For a fixed $n$, the corresponding
function space, determined by the rational values
$0<{1\over n}<\dots <{n-1\over n}<1$, 
is denoted by $C^{\sharp}_n(\mathcal P_x^T)$.
Then subset $\cup_{n\in \mathbf N}C^{\sharp}_n$ 
is dense in $\mathcal P_x^T$.

The constructions are implemented, first, on the holomorphic 
zone. By formula (\ref{w_1xy}), the well defined global 
WK-measure can directly be 
construc\-ted. An other slightly different construction technique
establishes an appropriate {\it linear functional on}
$C(\mathcal P_x^T)$, which defines 
the desired finite complex measure,  
$w_{\mathbf 1x}^{T(0)}$, on $\mathcal P_x^T$. 
We follow this idea for constructing the zonal Feynman measures. 
This functional is defined, first, on subspaces $C^{\sharp}_n$ by 
\begin{eqnarray}
\label{W_i^T}
\mathcal W^{T(0)}_{\mathbf in}(F)=\\
\int d^{(0)}_{\mathbf i}({T\over n},
a,m_1)\dots d^{(0)}_{\mathbf i}({T\over n},m_{n-1},m_n)
F(m_1,\dots ,m_n)dm_1\dots dm_n\nonumber \\
=\int e^{-{kT\over 2}\mathbf i}
e^{{e^{-2(T/n)\mathbf i}-1\over T/n}(T/n)(
\langle a,m_1+\mathbf iJ(m_1)\rangle +\dots +
\langle m_{n-1},m_n+\mathbf iJ(m_n)\rangle )}\nonumber \\
\cdot\delta^{(0)}(
a,m_1)\dots \delta^{(0)}(m_{n-1},m_n)
F(m_1,\dots ,m_n)dm_1\dots dm_n ,\nonumber
\end{eqnarray}
where $dm$ is normalized such that
$
\int |d^{(0)}_{\mathbf i}(t,p,m)|^2dm=
||d^{(0)}_{\mathbf i(t,p)}(m)||^2=1
$ 
holds for any fixed $p$ and $t$. By (\ref{intform}) one gets that 
$dm=\pi^{k/2}dm_E$,
where $dm_E$ is the Euclidean volume measure.

At this point of the proof only the first equation is needed. 
The second one will be used, later, 
for establishing the Feyn\-man-Kac formula. Let it be also mentioned
that functional
$\mathcal W^{T(0)}_{\mathbf 1n}(F)$ 
is similarly defined by means of the zonal
Wiener-Kac kernel. In this case all the $\mathbf i$'s, but the last ones
in formulas $\mathbf iJ(m_i)$,
should be exchanged for $1$ in the above formula. 
The Wiener-Kac functional is
uniformly bounded and absolutely continuous which, by the
limit $n\to\infty$, extends to the desired absolutely continuous
functional  
$\mathcal W^{T(0)}_{\mathbf 1}(F)$ which then
defines the desired measure 
$w^{T(0)}_{\mathbf 1x}$.

One can straightforwardly apply these ideas for constructing 
the zonal Feynman measure. The uniform 
boundedness of functionals (\ref{W_i^T}) 
can be established by proving the inequality 
$|\mathcal W_{\mathbf in}^{T(0)}(F)|\leq (2\pi )^{k\over 2}||F||$, 
for all $n$. The proof of this inequality can be carried out by 
the following formulas calculated by means of 
(\ref{intform}).
\begin{eqnarray}  
\label{intd_id_i}
\int |d^{(0)}_{\mathbf i}({T\over n},m,p)
d^{(0)}_{\mathbf i}({T\over n},p,q)|dp
\leq ||d^{(0)}_{\mathbf i({T/ n},m)}(p)||
||d^{(0)}_{\mathbf i({T/ n},q)}(p)||=1,\\
\label{int|d_i|}
\int |d^{(0)}_{\mathbf i}({T\over n},m,p)|||F||dp=\\
={||F||\over \pi^{k\over 2}}\int
e^{-{1\over 2}(|m|^2+|p|^2)+
cos({2T\over n})\langle m,p\rangle -sin({2T\over n})
\langle J(m),p\rangle }dp=(2\pi )^{k\over 2}||F||.\nonumber
\end{eqnarray}
If $n$ is even, by applying the first inequality in pairs in 
(\ref{W_i^T}), 
one proves that the norms of the functionals 
$\mathcal W_{\mathbf in}^{T(0)}$ 
are uniformly bounded by $1$. If $n$ is odd of the form
$n=2m+1$, then applying the first inequality for the first $m$ pairs
of $d_{\mathbf i}$'s in (\ref{W_i^T}) and the second one for the last 
$d_{\mathbf i}$, one has the desired inequality.
Therefore, by 
Riesz's and Stone-Weierstrass's theorems,
the approximating functionals have a unique continuous 
extension onto $C(\mathcal P_x^T)$.

Note that functional 
$\mathcal W^{T(0)}_{\mathbf in}$ 
defines approximatig measure
$w^{T(0)}_{\mathbf inxy}$ on the approximating path spaces 
$\mathcal P_{nxy}=M\times\dots\times M$
(n-times product) where the curves have fixed starting
point, $x$, and end point, $y$.
Then we have
\begin{equation}
w^{T(0)}_{\mathbf inxy}(\mathcal P_{nxy})
=d^{(0)}_{\mathbf i}(T,x,y),
\end{equation}
what immediately follows from the first equation of 
(\ref{W_i^T}). Therefore,
\begin{equation}
w^{T(0)}_{\mathbf ixy}(\mathcal P_{xy})=
\lim_{n\to\infty}w^{T(0)}_{\mathbf inxy}(\mathcal P_{nxy})
=d^{(0)}_{\mathbf i}(T,x,y).
\end{equation}
\medskip

\noindent{\bf Path-space measure induced by the 
holomorphic point-spread.}
Constructions of the zo\-nal measures $\nu_x^{T(0)}$ and 
$(\nu\overline\nu )_x^{T(0)}$  
on the path-space  
$\mathcal P_x^T$
by means of the point spread kernel  
\begin{equation}
\delta^{(0)}(X,Y)=\sum_i\varphi_i(X)\overline{\varphi}_i(Y)=
{1\over \pi^{k\over 2}}e^{-{1\over 2}(|X|^2+|Y|^2)+
\langle X,Y+\mathbf iJ(Y)\rangle}
\end{equation}
and the corresponding density kernel 
$\delta^{(0)}\overline{\delta}^{(0)}$ 
on the holomorphic zone can be established by the very same steps 
described above. These kernels
can be used for defining the bounded functionals 
$\mathcal N_{nxy}^{T(0)}$ and
$(\mathcal N\overline{\mathcal N})_{nxy}^{T(0)}$ 
as well as the corresponding finite measures $\nu_{nxy}^{T(0)}$
and $(\nu\overline\nu )_{nxy}^{T(0)}$. The functionals turn out to be 
uniformly bounded and, therefore, defining the sought measures by limit.
For the total measure we have
\begin{equation}
\nu^{T(0)}_{xy}(\mathcal P_{xy})=\nu^{T(0)}_{nxy}(\mathcal P_{nxy})
=\delta^{(0)}(x,y).
\end{equation}
The corresponding formula holds for the density measure.

To describe these measures more accurately, 
we compare them with the well
defined Wiener-Kac measure. First note the connection 
\begin{equation}
p^{(0)}_{1}(t,X,Z)=
e^{-{kt\over 2}}
e^{(e^{-2t}-1)
\langle X,Z+\mathbf iJ(Z)\rangle}\delta^{(0)}(X,Z)
\end{equation}
between $\delta^{(0)}$ and $p^{(0)}_{1}$.  
Thus, by taking the limit
$n\to \infty$
in the second equation of (\ref{W_i^T}), we have
\begin{equation}
\label{intdw1}
\int_{\mathcal P^T_x}f(\omega)
dw_{\mathbf 1x}^{T(0)}(\omega )=
\int_{\mathcal D^T_x} f(\omega) 
e^{-{kT\over 2}}
e^{-2\int_0^T|\omega (\tau )|^2d\tau}d\nu_x^{T(0)}(\omega ).
\end{equation}
This proves the Radon-Nikodym formula
\begin{equation}
dw_{\mathbf 1x}^{T(0)}(\omega )=
e^{-{kT\over 2}}
e^{-2\int_0^T|\omega (\tau )|^2d\tau}d\nu_x^{T(0)}(\omega ).
\end{equation}
The reversed Radon-Nikodym derivative is
\begin{equation}
\label{RN1}
d\nu_x^{T(0)}(\omega)=e^{{kT\over 2}
+2\int_0^T|\omega (\tau )|^2d\tau}dw_{\mathbf 1xy}^{T(0)}(\omega ),
\end{equation}
describing the considered measure in terms of the Wiener-Kac measure.
\medskip

\noindent
{\bf The Feynman measure on the 0-zone} 
is constructed by the same steps. In this
case the corresponding formulas are
\begin{equation}
\label{p_i^0}
p^{(0)}_{\mathbf i}(t,X,Z)=
e^{-{kt\over 2}\mathbf i}
e^{(e^{-2t\mathbf i}-1)
\langle X,Z+\mathbf iJ(Z)\rangle}\delta^{(0)}(X,Z).
\end{equation}
Thus, for any function $f\in C\mathcal P^T_x$, we have
\begin{eqnarray}
\label{int_P}
\int_{\mathcal P^T_x}f(\omega)dw_{\mathbf ix}^{T(0)}(\omega )
=\int_{\mathcal D^T_x} f(\omega) e^{-{kT\over 2}\mathbf i}
e^{-2\mathbf i\int_0^T|\omega (\tau )|^2d\tau}d\nu_x^{T(0)}(\omega ),\\
p_{\mathbf i}^{(0)}(T,x,y)=\int_{\mathcal D^T_{xy}} e^{(-{kT\over 2}
-2\int_0^T|\omega (\tau )|^2d\tau )
\,\mathbf i}d\nu_{xy}^{T(0)}(\omega ),
\nonumber \\
dw_{\mathbf ix}^{T(0)}(\omega)=e^{(-{kT\over 2}
-2\int_0^T|\omega (\tau )|^2d\tau )\mathbf i}
d\nu_{xy}^{T(0)}(\omega ).\nonumber
\end{eqnarray}

By the last equations of (\ref{W_i^T}) resp. (\ref{int_P}) we get
\begin{equation}
dw_{\mathbf ix}^{T(0)}(\omega)=e^{({kT\over 2}
+2\int_0^T|\omega (\tau )|^2d\tau )(1-\mathbf i)}
dw_{\mathbf 1xy}^{T(0)}(\omega ),
\end{equation}
which is the {\it Radon-Nikodym derivative of the zonal Feynman measure
with respect to the zonal Wiener-Kac measure}. This formula describes
the most direct connection between the two zonal measures. 

Formulas (\ref{intdw1}) and (\ref{int_P}) can be interpreted as 
{\it Feynman-Kac type 
formulas}. Originally, they stand for the Radon-Nikodym
derivative of the well defined Wiener-Kac measure with respect
to the Wiener measure concerning the Euclidean Laplacian $\Delta_X$.
On the zonal setting this idea can not be carried out
because the zones are not invariant with respect
to the Euclidean Laplacian and the zonal decomposition
with respect to this Laplacian is not defined. In the above version,
measure $d\nu_x^{T(0)}$ substitutes the Wiener measure.

The zonal Feynman measure of the set of curves connecting $x$ 
and $y$ is
\begin{equation}
w^{T(0)}_{\mathbf ixy}\big(\mathcal P_{xy}^T(M)\big)=
d_{\mathbf i}^{(0)}(T,x,y).
\end{equation}
This function is called zonal (in this case, ``holomorphic") 
{\it probability amplitude} and the function defined by
\begin{equation}
\rho^{T(0)}_{xy}=
d_{\mathbf i}^{(0)}(T,x,y)
\overline{d}_{\mathbf i}^{(0)}(T,x,y)
\end{equation}
is called zonal (holomorphic) {\it probability density}.
For a Borel set
$B$ on the X-space, integral $\int_B\rho^{(0)}_{x,y}dy$ (where dy is 
the normalized density described in (\ref{W_i^T})) 
measures the probability of that that
the point spread about $x$ can be caught, at the time $T$, 
among the point-spreads spread-ed about the points of $B$, meaning that 
$d_{\mathbf i}^{(0)}(T,x,y)$ 
is in the function space spanned by functions
$d_{\mathbf i}^{(0)}(T,b,y), b\in B$. For $B=\mathbf R^k$, 
this probability is 1.
 
Probability density can be defined also for any zonal function 
$\psi^{(0)}(t,X)$ satisfying 
the Schr\"odinger equation. 
By the convolution formula
\begin{equation}
\psi^{(0)}(t,X)=d^{(0)}_{\mathbf i}(t,X,Z)*_Z\psi (0,Z),
\end{equation}
such a function
is uniquely determined by the function defined for $t=0$. 
If the initial function $\psi^{(0)}(0,X)$ 
is normalized, so is $\psi^{(0)}(t,X)$, for any $t$. 
This statement readily follows from 
\begin{equation}
{1\over\mathbf i}
{\partial d^{(0)}_{\mathbf i}\over\partial t}
=-(H_{Z})_{X} d^{(0)}_{\mathbf i}\quad ,\quad
{1\over\mathbf i}{\partial \overline{d}^{(0)}_{\mathbf i}
\over\partial t}
=\overline{(H_{Z})_{X} d^{(0)}_{\mathbf i}}
\end{equation}
and from
$\int\overline{\Box f}g=\int\overline{f}{\Box g}$, for all
$f,g\in\mathcal H^{(0)}$.
Indeed, by them we have 
$\partial_t\int \psi\overline{\psi}=0$. This statement is known as
{\it conservation of probability under the action of the DF-flow.} 

The probability density is defined by 
$\rho_\psi^T=\psi\overline{\psi}$.
In (\ref{p_i^0}) the Dirac delta spread 
about an arbitrary point $x$ is the
initial function by which the probability amplitude is defined.
In general the initial function defines a starting continuum-spread 
(zonal object) and  
$\int_B\psi (T,X)\overline{\psi}(T,X)dX$
measures the probability that the zonal object can be caught on the 
Borel set $B$, at the time $T$.

The actual mathematical theorem exhibited in the above theorem 
can be formulated such that that {\it the holomorphic
Feyn\-man-Dirac flow defines a unitary semi-group}, 
$U_t^{(0)}$, on the holomorphic zone.
\medskip

\noindent{\bf Constructions on general zones.} Since the conservation
of probability holds on each zone,
formula (\ref{intd_id_i}) extends from the 0-zone to any zone. Also  
(\ref{int|d_i|}) extends to each zone. In fact, functions 
$d_{\mathbf i}^{(a)}({T\over n},m,p)$ 
depending on $p$ are 
absolute integrable for any fixed $n$ and $m$, furthermore,
\[
\lim_{n\to infty}
\int |d_{\mathbf i}^{(a)}({T\over n},m,p)|dp=
\int |\delta^{(a)}(m,p)|dp.
\]
From this convergence an appropriate upper bound for the generalized
integral (\ref{int|d_i|}) follows, proving that the zonal
Feynman measures are well defined on each zone. Thus we have
\begin{theorem}{\bf (Zonal Path Integral Theorem)}
Let $H_Z$ be a non-degene\-rated Zee\-man-Hamilton operator
with parameters $\lambda_1,\dots ,\lambda_r$. Then, for each zone,
both zonal kernels $p_\sigma^{(a)}, \sigma =1,\mathbf i$, generate
continuous measures $dw_{\sigma xy}^{T(a)}$ on the
path-space 
$\mathcal P_{xy}^{T}$ 
such that the corresponding measures 
of the whole
path-space are 
$p_\sigma^{(a)}(T,x,y)$. 
They are called zonal 
Wiener-Kac and Feynman measures respectively. 

Functions
$\pi^{k/2}d_{\mathbf i}^{(a)}(T,x,y)
\overline{d}_{\mathbf i}^{(a)}(T,x,y)$ resp. 
$\psi^{(a)}\overline{\psi}^{(a)}$,
where $\psi^{(a)}(t,X)$ is a zonal function satisfying the 
Schr\"odinger equation and $||\psi^{(a)}(0,X)||=1$, are called zonal 
probabilistic density functions. For any $T$, the norm of these
functions is $1$. This phenomena is called 
{\it conservation of probability}.
Since $\psi^{(a)}(t,X)=\int d_{\mathbf i}(t,Z)\psi^{(a)}(0,Z)$,
the mathematical meaning of this statement is that the Feynman-Dirac
flow generates, on each zone, a semi-group $U^{(a)}_t$ 
of unitary transformations.

A zonal measure, 
$d\nu_{xy}^{T(a)}(\omega)$, 
is determined on 
$\mathcal P_{xy}^{T}$
also by the zonal Dirac spread 
$\delta^{(a)}_{\lambda_1\dots \lambda_r}$. On the holomorphic zone
the Radon-Nikodym derivative of the measure $dw_{\mathbf ix}^{T(0)}$ 
with respect to the $d\nu_{x}^{T(0)}$
is described in the following Feynman-Kac type formulas.
\begin{eqnarray}
\label{int_Pfdw_i}
\int_{\mathcal P^T_x}f(\omega)dw_{\mathbf ix}^{T(0)}
=\int_{\mathcal P^T_x} f(\omega) 
e^{\sum\lambda_i(-{k_iT\over 2}\mathbf i
-2\mathbf i\int_0^T|\omega_i (\tau )|^2d\tau )}d\nu_x^{T(0)},\\
d_{\mathbf i}^{(0)}(T,x,y)=\int_{\mathcal P^T_{xy}} e^{\sum \lambda_i
(-{k_iT\over 2}
-2\int_0^T|\omega_i (\tau )|^2d\tau )\,\mathbf i}
d\nu_{xy}^{T(0)}(\omega ).\nonumber
\end{eqnarray}
The Radon-Nikodym derivative of
$\nu_{x}^{T(0)}$
with respect to
$w_{\mathbf 1x}^{T(0)}$ is established in (\ref{RN1}).

On the holomorphic zone, the Feynman-Kac type formula regarding 
the WK-measure is:  
\begin{equation}
d_{1}^{(0)}(T,x,y)=\int_{\mathcal P^T_{xy}} e^{\sum \lambda_i
(-{k_iT\over 2}
-2\int_0^T|\omega_i (\tau )|^2d\tau )}d\nu_{xy}^{T(0)}(\omega ).
\end{equation}

These two Feynman-Kac type formulas establish the Radon-Nikodym
derivative
\begin{equation}
dw_{\mathbf ixy}^{T(0)}(\omega )=e^{\sum\lambda_i({k_iT\over 2}
+2\int_0^T|\omega_i (\tau )|^2d\tau )
(1-\mathbf i)}dw_{\mathbf 1xy}^{T(0)}
(\omega )
\end{equation}
of the holomorphic Feynman measure with respect to the holomorphic 
Wiener-Kac measure, providing the most direct connection between 
these two holomorphic measures.  

These zonal measures can be constructed also by means of
the zonal kernels $p_{\sigma}^{(a)}$ and $b_{\sigma}^{(a)}$. 
The two modifications should be implemented regarding the
Feynman-Kac type formulas are the following ones: (1) 
The above functions are
multiplied by $e^{-2t|Z_\gamma |^2}$. (2) The Dirac $\delta$-spread 
on the torus bundle regarding the function space $FW^\gamma$
is defined such that $\delta^{(a)}_\gamma$ is multiplied by
$e^{2\mathbf i\langle Z_\gamma,Z_x-Z_y\rangle}$ and measure $d\nu$
regards this $\delta$-spread. Then the Feynman-Kac formulas for
these two zonal kernels have the same form.  
\end{theorem}


\begin{thebibliography}{GSWH}


\bibitem[AB]{ab}
Y. Aharanov and D. Bohm:
\newblock Significance of electromagnetic potentials in the 
quantum theory.
\newblock {\em Phys. Rev.}, 115:485--491, 1959.

\bibitem[AC]{ac}
Y. Aharanov and A. Casher:
\newblock Ground state of spin-$1/2$ charged particle in a 
two-dimensional magnetic field.
\newblock{\em Phys. Rev. A},
19:2461-2462, 1979.


\bibitem[A]{a}
O. Alvarez:
\newblock Theory of strings with boundaries: fluctuations, topology
and quantum geometry.
\newblock {\em Nucl. Phys.}, 216:125--184, 1983.

\bibitem[AHS]{ahs}
J. Avron, I. Herbst, B. Simon:
\newblock Schr\"odinger operators with magnetic fields, I.,
General interaction.
\newblock{\em Duke Math. J.},
45:847--883, 1978.


\bibitem[Bo]{bo}
D. Bohm:
\newblock Quantum Theory,
\newblock Dover, 1979.

\bibitem[BSB]{bsb}
S. Duplij, W. Siegel, J. Bagger (editors):
\newblock Concise Encyclopedia of Supersymmetry and Noncommutative
Structures in Mathematics and Physics (SUSY Encyclopedia).
\newblock Kluwer, 2003.

\bibitem[Ee]{ee}
J. Eells:
\newblock Random walk on the fundamental group.
\newblock{\em Proc. Symp. Pure Math.},  27:211--217, 1975.

\bibitem[F]{f}
R. P. Feynman:
\newblock QED.
\newblock Princeton Univ. Press, 1988.


\bibitem[FH]{fh}
R. P. Feynman, A. R. Hibbs:
\newblock Quantum Mechanics and Path Integrals.
\newblock McGraw-Hill, 1965.

\bibitem[GSW]{gsw}
M. B. Green, J. H. Schwarz, E. Witten:
\newblock Superstring Theory. 
\newblock Cambridge Univ. Press, 1999.

\bibitem[H]{h}
S. W. Hawking:
\newblock Zeta function regularization of path integrals in 
curved space time.
\newblock{\em Comm. Math. Phys.}, 55:133--148, 1977.

\bibitem[Ka]{ka}
A. Kaplan:
\newblock Riemannian nilmanifolds attached to Clifford modules.
\newblock{\em Geom. Dedicata}, 11:127--136, 1981.

\bibitem[LL]{ll}
L. D. Landau, E. M. Lifshitz:
\newblock Quantum Mechanics.
\newblock Pergamon Press LTD, 1958.

\bibitem[Me]{me}
A. Messiah:
\newblock Quantum Mechanics.
\newblock Dover Publ. INC, 1999.

\bibitem[M\"u]{mu}
W. M\"uller:
\newblock Relative zeta functions, relative determinants and
scattering theory.
\newblock{\em Comm. Math. Phys.}, 192:309--347, 1998.

\bibitem[N]{n}
E. Nelson:
\newblock Feynman integrals and the Schr\"odinger equation.
\newblock{\em J. Math. Phys.}, 5(3):332--345, 1963.

\bibitem[OPS]{ops}
B. Osgood, R. Philips, P. Sarnak:
\newblock Extremals of determinants of Laplacians.
\newblock{\em J. Func. Anal.}, 80:148--211, 1988.

\bibitem[Po]{po}
A. M. Polyakov:
\newblock Quantum geometry of bosonic and Fermionic strings.
\newblock{\em Phys. Lett. B}, 103:207--213, 1981.

\bibitem[RS]{rs}
M. Reed, B. Simon:
\newblock Methods of Modern Mathematical Physics, II.
\newblock Academic Press, London, 1978.

\bibitem[RSz]{rsz}
F. Riesz, B. Sz\"okefalvi-Nagy:
\newblock Functional Analysis.
\newblock Dover, 1990.

\bibitem[Sch]{sch}
S. S. Schweber:
\newblock QED and the Man Who Made It: Dyson, Feynman, Schwinger, and
Tomonaga.
\newblock Princeton Univ. Press, 1994.

\bibitem[Schw]{schw}
J. Schwinger (Editor):
\newblock Selected Papers on Quantum Electrodynamics.
\newblock Dover Publ. INC, 1958.

\bibitem[Si]{si}
B. Simon:
\newblock Functional Integration and Quantum Physics.
\newblock Academic Press, London, 1979.

\bibitem[Sz1]{sz1}
Z. I. Szab\'o:
\newblock Locally non-isometric yet super isospectral spaces.
\newblock{\em Geom. funct. anal. (GAFA)}, 9:185--214, 1999.

\bibitem[Sz2]{sz2}
Z. I. Szab\'o:
\newblock Isospectral pairs of metrics on 
balls, spheres, and other manifolds with different local geometries.
\newblock{\em Ann. of Math.}, 154:437--475, 2001.

\bibitem[Sz3]{sz3}
Z. I. Szab\'o:
\newblock A cornucopia of isospectral pairs of metrics on 
spheres with different local geometries.
\newblock{\em Ann. of Math.}, 161:343--395, 2005.

\bibitem[Sz4]{sz4}
Z. I. Szab\'o:
\newblock Correction and improvement added to ``Isospectral pairs 
of metrics on balls and 
spheres with different local geometries".
\newblock{\em DG/0510202 (submitted)}

\bibitem[Sz5]{sz5}
Z. I. Szab\'o:
\newblock Theory of zones on Zeeman manifolds: 
A new approach to the infinities
of QED.
\newblock {\em DG/0510660 (submitted)}

\bibitem[Sz6]{sz6}
Z. I. Szab\'o:
\newblock Pauli-Dirac operators and anomalous zones on Zeeman manifolds.
\newblock submitted.


\bibitem[Sze]{sze}
G. Szeg\"o:
\newblock Orthogonal Polynomials.
\newblock AMS Providence, Rhode Island, vol(23), 1939.

\bibitem[Ta]{ta}
M. E. Taylor:
\newblock Partial Differential Equations. 
\newblock Springer, 1996.

\bibitem[To]{to}
S. Tomonaga:
\newblock The Story of Spin. 
\newblock Univ. of Chicago Press, 1997.

\bibitem[Ton1]{ton1}
A. Tonomura, N. T. Matsuda, R. Suzuki, A. Fukuhara, N. Osakabe,
H. Umezaki, J. Endo, K. Shinagawa, Y. Sugita, and H. Fujiwara:
\newblock Observation of Aharanov-Bohm effect by electron holography. 
\newblock {\em Phys. Rev. Lett.}, 48:1443--1446, 1982.

\bibitem[Ton2]{ton2}
A. Tonomura, N. Osakabe, T. Matsuba, T. Kawasaki, J. Endo, S. Yano,
and H. Yamada:
\newblock Evidence for Aharanov-Bohm effect with magnetic field 
completely shielded from electron wave. 
\newblock {\em Phys. Rev. Lett.}, 56:792-795, 1986.

\bibitem[V]{v}
M. Veltman:
\newblock Facts and Mysteries in Elementary Particle Physics. 
\newblock World Scientific, 2003.

\bibitem[W]{w}
S. Weinberg:
\newblock The Quantum Theory of Fields I, II. 
\newblock Cambridge Univ. Press, 1995.

\bibitem[We]{we}
H. Weyl:
\newblock The Theory of Groups and Quantum Mechanics. 
\newblock Dover, 1950.




\end{thebibliography}
\end{document}